\documentclass{amsart}

\usepackage[utf8]{inputenc}
\usepackage[english]{babel}
\usepackage{microtype}
\usepackage{calligra}
\usepackage{amsthm}
\usepackage{amsmath}
\usepackage{amssymb}

\usepackage{graphicx}
\usepackage[
  algo2e,
  algosection,
  linesnumbered,
  ruled
]{algorithm2e}
\usepackage{caption}
\usepackage{subcaption}

\usepackage{tikz}
\usepackage{pgfplots}
\usetikzlibrary{plotmarks}
\usepackage{pgfplotstable}
\pgfplotsset{compat=1.12}
\usetikzlibrary{external}
\usetikzlibrary{spy}
\usepackage{wrapfig}

\usepackage{todonotes}
\usepackage{color}
\usepackage{setspace}
\usepackage{hyperref}

\definecolor{sixclassRdYlBu1}{rgb}{0.84,0.19,0.15}
\definecolor{sixclassRdYlBu2}{rgb}{0.99,0.55,0.35}
\definecolor{sixclassRdYlBu3}{rgb}{1.0,0.88,0.56}
\definecolor{sixclassRdYlBu4}{rgb}{0.88,0.95,0.97}
\definecolor{sixclassRdYlBu5}{rgb}{0.57,0.75,0.86}
\definecolor{sixclassRdYlBu6}{rgb}{0.27,0.46,0.71}
\definecolor{gelb}{rgb}{1,0.647,0}
\definecolor{oran}{rgb}{1,0.3137,0.035}
\definecolor{lila}{rgb}{0.4235,0,0.894}
\definecolor{gruen}{rgb}{0,0.7216,0}
\definecolor{dunkelgruen}{rgb}{0,0.47,0}
\definecolor{warm_red}{rgb}{0.835,0.0627,0.0627}
\definecolor{grau}{rgb}{0.6078,0.6078,0.549}

\newcommand{\muv}{\boldsymbol{\mu}}

\newcommand{\Params}{\mathcal{P}}
\newcommand{\R}{\mathbb{R}}
\newcommand{\N}{\mathbb{N}}
\newcommand{\Grid}{\mathcal{T}_H}
\newcommand{\EstGrid}{\tilde{\mathcal{T}}_H}
\newcommand{\grid}{\tau_h}
\DeclareMathOperator{\diam}{diam}
\newcommand{\faces}{\tau_h^\gamma}
\newcommand{\sigmafine}{\sigma}

\newcommand{\param}{\boldsymbol{\mu}}
\newcommand{\paramFixed}{\boldsymbol{\bar \mu}}
\newcommand{\paramHat}{\boldsymbol{\hat{\mu}}}
\DeclareMathOperator{\integralend}{d}
\newcommand{\intend}[1]{\integralend\hspace{-1.25pt}#1}
\newcommand{\dx}{\intend{x}}
\newcommand{\ds}{\intend{s}}

\newcommand{\norm}[1]{{\left\|{#1}\right\|}}
\newcommand{\energynorm}[2]{{|||{#1}|||_{#2}}}
\newcommand{\dualpair}[2]{{\left\langle #1 , #2 \right\rangle }}
\newcommand{\SP}[2]{(#1, #2)}
\newcommand{\etalchar}[1]{$^{#1}$}

\makeatletter
\newcommand{\enorm}{\@ifstar\@enorms\@enorm}
\newcommand{\@enorms}[1]{%
  \left|\mkern-1.5mu\left|\mkern-1.5mu\left|
   #1
  \right|\mkern-1.5mu\right|\mkern-1.5mu\right|
}
\newcommand{\@enorm}[2][]{%
  \mathopen{#1|\mkern-1.5mu#1|\mkern-1.5mu#1|}
  #2
  \mathclose{#1|\mkern-1.5mu#1|\mkern-1.5mu#1|}
}
\makeatother

\newcommand{\Ell}{\mathcal{R}_{\rm ell}}

\newcommand{\hnS}{\hspace{-1.25pt}}
\newcommand{\Hdiv}{H_\text{div}}
\newcommand{\mean}[1]{\big<#1\big>}
\newcommand{\jump}[1]{[#1]}
\newcommand{\divergence}{\nabla\hnS\hnS\cdot\hnS}
\newcommand{\gradient}{\nabla\hnS}
\newcommand{\gradienth}{\nabla\hnS_h}
\newcommand{\mydot}{\hnS\cdot\hnS}

\newcommand{\spanlin}{\operatorname{span}} 
 
\newcommand{\dist}{\operatorname{dist}} 
\newcommand{\cest}{{c_\mathrm{est}(n_t, \varepsilon_{\mathrm{testfail}})}}
\newcommand{\ceff}{{c_\mathrm{eff}(n_t, \varepsilon_{\mathrm{testfail}})}}
\newcommand{\range}{\operatorname{range}}
 
\newcommand{\subsod}{O}

\newcommand{\puconstant}{c_N}
\newcommand{\ovlpconstant}{c_\mathrm{ovlp}}
\newcommand{\pufuncconstant}{{c_{\mathrm{pu}}}}

\usepackage[section]{algorithm}
\theoremstyle{plain}
\newtheorem{theorem}[algorithm]{Theorem}
\newtheorem{lemma}[algorithm]{Lemma}
\newtheorem{proposition}[algorithm]{Proposition}
\newtheorem{remark}[algorithm]{Remark}
\newtheorem{definition}[algorithm]{Definition}
\newtheorem{example}[algorithm]{Example}
\newtheorem{corollary}[algorithm]{Corollary}

\title[Localized model reduction]{Localized model reduction for parameterized problems}

\author[A Buhr]{Andreas Buhr}
\address{Mathematics M\"unster, Einsteinstr. 62, D-48149 M\"unster, Germany}
\author[L Iapichino]{Laura Iapichino}
\address{Department of Mathematics and Computer Science, TU Eindhoven, Eindhoven, The Netherlands}
\author[M Ohlberger]{Mario Ohlberger}
\address{Mathematics M\"unster, Einsteinstr. 62, D-48149 M\"unster, Germany}
\author[S Rave]{Stephan Rave}
\address{Mathematics M\"unster, Einsteinstr. 62, D-48149 M\"unster, Germany}
\author[F Schindler]{Felix Schindler}
\address{Mathematics M\"unster, Einsteinstr. 62, D-48149 M\"unster, Germany}
\author[K Smetana]{Kathrin Smetana}
\address{University of Twente, Faculty of Electrical Engineering, Mathematics \& Computer Science, Zilverling, P.O. Box 217, 7500 AE Enschede, The Netherlands}

\date{\today}
\thanks{The authors from Mathematics M\"unster are funded by the Deutsche Forschungsgemeinschaft (DFG, German Research Foundation) 
under Germany 's Excellence Strategy – EXC 2044 – 390685587, Mathematics M\"unster: Dynamics – Geometry - Structure.
F.~Schindler acknowledges funding by the Deutsche Forschungsgemeinschaft (DFG, German Research Foundation) under contract SCHI 1493/1-1.}



\begin{document}



\begin{abstract}
In this contribution we present a survey of concepts in localized model order reduction methods for parameterized partial differential equations. The key concept of localized model order reduction is to construct local reduced spaces that have only support on part of the domain and compute a global approximation by a suitable coupling of the local spaces. In detail, we show how optimal local approximation spaces 
can be constructed and approximated by random sampling. 
An overview of possible conforming and non-conforming couplings of the local spaces is provided and corresponding localized 
a posteriori error estimates are derived. 
We introduce concepts of local basis enrichment, which includes a discussion of adaptivity. 
Implementational aspects of localized model reduction methods are addressed. 
Finally, we illustrate the presented concepts for multiscale, linear elasticity and fluid-flow problems, providing several numerical experiments. \\
This work has been accepted as a chapter in P. Benner, S. Grivet-Talocia, A. Quarteroni, G. Rozza, W.H.A. Schilders, L.M. Sileira. Handbook on Model Order Reduction.  Walter De Gruyter GmbH, Berlin, 2019+.
\end{abstract}

\maketitle


\section{Introduction}\label{sec:intro}
Projection based model order reduction has become a mature technique for simulations of large 
classes of parameterized systems; 
for an introduction, we refer to the text books and survey \cite{BCOW17,HesRozSta2016,MR3379913,BGW15}. However, especially for large-scale and multi-scale problems the ``standard'' model order reduction approach exhibits several limitations: Curse of parameter dimensionality in the sense that many parameters require prohibitively large reduced spaces, no topological flexibility, and possibly high computational costs and storage requirements in the offline stage for instance due to large computational domains. Localized model order reduction methods, which combine approaches 
from model order reduction, multiscale methods and/or domain decomposition techniques, overcome or significantly mitigate those limitations. As an further advantage, they allow using reduced spaces of different dimensions in different parts of the computational domain and accommodate (local) changes of the geometry and the partial differential equation (PDE) in the online stage. The key idea of localized model order reduction is to construct local reduced spaces on (unions of) subdomains of the decomposed computational domain and couple the local reduced spaces across interfaces either in a conforming or non-conforming manner. In this chapter we investigate localized model order reduction for linear coercive elliptic parameterized problems; inf-sup stable problems have for instance been considered in \cite{HuKnPa13b} and parabolic and nonlinear problems will be briefly dicussed at the end of this chapter.

We discuss both conforming and non-conforming localized approximations. Prominent examples for a conforming localization for non-parametric PDEs are the partition of unity method \cite{BabMel97},
the generalized finite element method (GFEM) \cite{BaCaOs94,BaBaOs04,BabMel97,BabLip11} as well as 
component mode synthesis (CMS) \cite{Hur65,BamCra68,Bou92},  \cite{HetLeh10,JaBeLa11}. 

A combination of domain decomposition and reduced basis (RB) methods has first been considered in the reduced basis element method (RBEM) \cite{MadRon02,MadRon04,LovMadRon06}, where the local RB approximations are coupled by Lagrange multipliers in a non-conforming manner. The reduced basis hybrid method \cite{IaQuRo12} extends the RBEM by additionally considering a coarse FE discretization on the whole domain to account for continuity of normal stresses in the context of Stokes equations.  Alternatively, a non-conforming coupling can be realized say by penalization as in the local reduced basis discontinuous {G}alerkin approach \cite{KOH2011}, 
the localized reduced basis multiscale method (LRBMS) \cite{AlbrechtHaasdonkEtAl2012,OS2015,ORS2017}, the 
discontinuous Galerkin reduced basis element method \cite{APQ2016}, or the generalized multiscale discontinuous Galerkin method \cite{CEL2017}. The static condensation reduced basis element (scRBE) method  \cite{HuKnPa13,HuKnPa13b,EftPat13a,SmePat16} combines intra-element RB approximations similar to the RBEM with coupling techniques from CMS resulting in a conforming approximation. A similar approach is pursued by the ArbiLoMod \cite{Betal17} that also allows for arbitrary (non-parametric) local changes of the underlying equations and/or the geometry. 

In the context of the proper generalized decomposition (PGD) method (for a review see for instance \cite{ChiAmmCue2010,ChLaCu11,ChKeLe14}) a domain decomoposition strategy has been proposed in \cite{HNC:18} and in \cite{PerErnVen10} hierarchical model reduction \cite{VogBab81,PerErnVen10, OhlSme14, SmeOhl17} has been combined with an iterative substructuring method.

Concerning the generation of local approximation spaces we focus on empirical training (see for instance \cite{EftPat13a,BabLip11,SmePat16}), i.e. local reduced spaces generated from local solutions of the PDE, and adaptive basis enrichment. In detail, we present local approximation spaces that are optimal in the sense of Kolmogorov and can be constructed by solving a local so-called transfer eigenvalue problem on the space of local solutions of the PDE. Optimal local approximation spaces for subdomains have first been proposed in \cite{BabLip11} and for interfaces and parametrized PDEs in \cite{SmePat16}. We will also show how those optimal approximation spaces can approximated by random sampling \cite{BuhSme18}.

A localizable a posteriori error estimator is crucial for an adaptive enrichment of the local reduced spaces where the reduced approximation is not accurate enough. Such an adaptive basis enrichment is one way to approach ``optimal'' computational complexity within outer-loop applications such as optimal control, inverse problems or Monte Carlo methods. With this respect, we will also
present a framework for localized residual based error control \cite{Betal17,Sme15} as well as localized a posteriori error estimation based on 
flux reconstruction \cite{ESV2010,OS2015}.

Naturally, the presented methods for localized model reduction share a lot of features with domain decomposition techniques 
and multiscale methods. 
We particularly refer to domain decomposition and preconditioning techniques with multiscale coarse spaces
such as \cite{AarHou02,GrLeSc07,GalEfe10} or the more recent contributions \cite{SDHNPS2014,GL2017,HKKR18}. 
In the context of the {FETI}-{DP} iterative substructuring method we refer to \cite{ManSou07,KlRaRh16}.
For multiscale problems there has been a tremendous development of suitable numerical methods in the last two decades
including the multiscale finite element method (MsFEM)  \cite{HouWu97,Efendiev:Hou:Ginting:2004,Efendiev:Hou:2009,HOS14}, 
the heterogeneous multiscale method (HMM) \cite{E:Engquist:2003,E:Engquist:2005,Ohlberger:2005,Abdulle2005}, the variational multiscale method (VMM) \cite{Hughes:1995,Hughes:et-al:1998,Larson:Malqvist:2005}  or the more recent local orthogonal decomposition (LOD) \cite{MalPet:14,HMP:14}. 
Model reduction can be used to accelerate the solution of localized problems which occur in multiscale methods, see e.g. \cite{Abdulle20127014,AH2015}.
Similar to the methods presented in this chapter the generalized multiscale finite element method (GMsFEM) \cite{Efendiev2013,Chung2014,Chung2014b} relies on a Galerkin projection on local subspaces, but in contrast uses ideas from multiscale methods to construct the local bases.  
A connection between multiscale methods and domain decomposition has recently been investigated in \cite{KoPeYs18,KoPoYs17,KorYse16}.

This chapter is organized as follows. In Sec. \ref{sec:setting} we introduce the problem setting and basic notation for localized model order reduction of coercive 
variational problems. Concepts for conforming and non-conforming coupling of approximation spaces are presented in Sec. \ref{sec:coupling}. 
Sec. \ref{sec:prep} deals with the preparation of local approximation spaces. Particularly, the construction of optimal local approximation spaces and their approximation via random sampling is presented and illustrated with numerical experiments for linear elasticity. In Sec. \ref{sec:a_posteriori} we present two abstract frameworks for localized a posteriori error estimation and give 
exemplifications for conforming and non-conforming localized model reduction approaches. Localized a posteriori error estimators are the key ingredient for 
basis enrichment strategies and online adaptivity that are presented in Sec. \ref{sec:adaptivity}. 
Computational aspects are discussed in Sec.  \ref{sec:computation} and numerical experiments for multiscale problems and fluid flow are presented in Sec. 
  \ref{sec:examples}.
We conclude by showing possible extensions to parabolic and nonlinear problems in Sec. \ref{sec:perspectives}.


\section{Parameterized partial differential equations and localization} \label{sec:setting}

Let $\Omega \subset \mathbb{R}^{d}$, $d=1,2,3$, be a large, bounded domain with Lipschitz boundary. Let us further introduce a Hilbert space $V$ such that  $[H^{1}_{0}(\Omega)]^{z} \subset V \subset [H^{1}(\Omega)]^{z}$, $z=1,\hdots,d$ and denote by $V'$ the dual space of $V$. Moreover, we introduce the compact set of admissible parameters $\mathcal{P}  \subset \R^p$, $p \in \N$. We consider the following variational problem.

\begin{definition}[Parameterized coercive problem in variational form]
  \label{def:parametrized_problem}
For any parameter $\param \in \Params$ find $u(\param) \in V$, such that
\begin{align}
  a\big(u(\param), v; \param\big) &= f(v; \param) &&\text{for all } v \in V.
  \label{eq:weak_solution}
\end{align}
\end{definition}

Here, $f(\cdot; \param) \in V'$ and $a(\cdot, \cdot; \param): V \times V \rightarrow \mathbb{R}$ denote parametric linear and bilinear forms, the latter being continuous and coercive w.r.t.\ the norm $\|\cdot\|_V$ induced by the inner product $(\cdot, \cdot)_V: V \times V \to \R$. That is, there exist constants $0 < \alpha \leq \alpha(\param) \leq \gamma(\param) \leq \gamma$, such that for any $\param \in \Params$,
\begin{align}
  a(v, w; \param) &\leq \gamma(\param)\; \|v\|_V\; \|w\|_V &&\text{for all } v, w \in V,\notag\\
  a(v, v; \param) &\geq  \alpha(\param)\; \|v\|_V^2 &&\text{for all } v \in V.\notag
\end{align}
Let us denote the energy norm of $u$ for parameter $\param$ as $\energynorm{u}{\param} := a(u, u; \param)^{1/2}$.
Problem \eqref{eq:weak_solution} thus admits a unique solution for all $\param \in \Params$ owing to the Lax-Milgram theorem. 
Examples for \eqref{eq:weak_solution} include elliptic multiscale problems, incompressible fluid flow or linear elasticity 
as detailed in the following.  
We will consider Neumann boundary conditions on $\Gamma_{N}$ and Dirichlet boundary conditions on $\Gamma_{D}$, 
where $\Gamma_{N}, \Gamma_{D}$ are non-overlapping and $\Gamma_{N} \cup \Gamma_{D} = \partial \Omega$.
To simplify notations, homogenous boundary conditions on $\partial \Omega$ will be prescribed in most places.

\begin{example}[Parametric elliptic multiscale problems]
  \label{ex:multiscale_problem}
  With $V = H^1_0(\Omega)$, the pressure equation in the context of two-phase flow in porous media (obtained through Darcy's law) reads: given a collection of sources and sinks $q \in L^2(\Omega)$, a parametric and possibly highly heterogeneous permeability field $\kappa: \Params \to L^\infty(\Omega)^{d \times d}$, find for each $\param \in \Params$ the global pressure $u(\param) \in V$, such that
  \begin{align}
    \label{eq:parametric-multiscale}
    -\nabla\cdot\big( \kappa(\param) \gradient{u(\param)}\big) = q&&\text{in a weak sense in }V'.
  \end{align}
  If the smallest eigenvalue of $\kappa(\param)$ is bounded from below away from zero for all $\param \in \Params$, we can consider this to be an example of Definition \ref{def:parametrized_problem} by setting $a(u, v; \param) := \int_\Omega (\kappa(\param)\gradient{u})\cdot \nabla v \dx$ and $f(v; \param) := \int_\Omega q v \dx$.
  In the context of instationary two-phase flow, \eqref{eq:parametric-multiscale} needs to be solved in each time step for varying total mobilities (modelled by the parametric nature of $\kappa$), while the permeability field $\kappa$ typically resolves fine geological structures and thus requires a very fine computational grid compared to the size of $\Omega$ (see \cite{OS2015} and the references therein).
\end{example}

\begin{example}[Incompressible fluid flow]
The Stokes and Navier-Stokes equations represent a model of the flow motion for a viscous Newtonian incompressible fluid. In the steady case it can be formulated as follows:
 \begin{equation}\label{NSequation}
\left\{
\begin{array}{ll}
-\nu \Delta {\bf y} + \delta ({\bf y} \cdot \nabla) {\bf y}+  \nabla p  = \mathbf{f}  & \quad \textrm{in }\Omega \vspace{0.05cm} \\
\nabla \cdot {\bf y} =0 & \quad  \textrm{in } \Omega \vspace{0.15cm} \\
{\bf y}= \mathbf{g}_D & \quad  \textrm{on } \Gamma_{D}  \\
-p \mathbf{n} + \nu \displaystyle \frac{\partial {\bf y}}{\partial \mathbf{n}}  =  \mathbf{g}_N & \quad  \textrm{on } \Gamma_{N}, 
\end{array}
\right.
\end{equation}
\noindent where $({\bf y},p)$ are the velocity and the pressure fields defined on the computational domain  $\Omega$.
The first equation expresses the linear momentum  conservation,  the second one the mass conservation, which is also called the continuity equation. Here 
$\mathbf{f}$ 
denotes a forcing term per unit mass, $\mathbf{g}_D$ and $\mathbf{g}_N$ are the functions addressing the Dirichlet and Neumann boundary conditions respectively on $\Gamma_{D}$ and $\Gamma_{N}$. The parameter $ \nu = \sigma/\rho$ denotes the kinematic viscosity, being $\rho$ the  density and $\sigma$ the  viscosity of the fluid.
Navier-Stokes equations correspond to the case $\delta = 1$, here we consider only $\delta = 0$, the convective term is neglected, obtaining the steady Stokes equations, which provide a model in the case of slow motion of fluids with very high viscosity.\\
We denote the functional spaces for velocity and pressure fields by 
$X = (H_{0,\Gamma_D}^1(\Omega))^d$, $Q = L^2(\Omega)$, 
respectively, where
$H_{0,\Gamma_D}^1(\Omega) = \{y \in H^1(\Omega) : y \vert_{\Gamma_D} = 0 \}$. Moreover, for simplicity, we assume that $\mathbf{g}_D=0$ (otherwise the lift function is required). The corresponding weak form of  the Stokes equations \eqref{NSequation} reads:
 find $({\bf y}, p) \in X \times Q$ such that 
 \begin{equation*}
\label{weakNS} 
\begin{array}{rl}
  \displaystyle   \nu\int_{\Omega}   \nabla \mathbf{y} : \nabla
    \mathbf{w} \ d\Omega   - \displaystyle    \int_{\Omega} p \, \nabla\cdot  \mathbf{w} \ d\Omega  &=  \displaystyle    \int_{\Omega}    \mathbf{f} \cdot
  \mathbf{w} \ d\Omega +    \int_{\Gamma_N}   \mathbf{g}_N \cdot
  \mathbf{w} \ d\Gamma ,   \hspace{.5cm} \forall \, \mathbf{w} \in X \\[12pt]
  \displaystyle \hspace{1.85cm} \int_{\Omega}    q \, \nabla \cdot \mathbf{y} \ d\Omega    &=   0, \hspace{1.55cm} \forall \, q \in
  Q.
\end{array} 
\end{equation*}
In a parameterized setting, the
input-parameter vector $\param$ 
may characterize either  the geometrical configuration  or physical properties, boundary data and sources of the problem. \\
 Denoting by $V $ the product space given by $V = X \times Q$, defining by $u(\param)  = ({\bf y}(\param),p(\param)) \in V$ and  $v=({\bf w},q)$, the parametrized abstract formulation \eqref{weakNS} can be rewritten in the following  form: find $u(\param) = ({\bf y}(\param),p(\param)) \in V$ s.t.
\begin{equation}\label{compactNS}
a(u(\param), v; \param) = {f}(v; \param), \quad \forall \, v \in V \vspace{-0.1cm}
\end{equation}
where \vspace{-0.1cm}
\begin{equation}\label{lhs_NS}
a(u,v; \param) = \nu\int_{\Omega}   \nabla \mathbf{y} : \nabla
    \mathbf{w} \ d\Omega   - \displaystyle    \int_{\Omega} p \, \nabla\cdot  \mathbf{w} \ d\Omega  - \displaystyle    \int_{\Omega} q \, \nabla\cdot  \mathbf{y} \ d\Omega, 
\end{equation}
\begin{equation}\label{rhs_NS}
f(v; \param) =   \displaystyle    \int_{\Omega}    \mathbf{f} \cdot
  \mathbf{w} \ d\Omega +    \int_{\Gamma_N}   \mathbf{g}_N \cdot
  \mathbf{w} \ d\Gamma. \vspace{-0.1cm}
\end{equation}

\end{example}

\begin{example}[Linear elasticity]
We assume that $\Omega \subset \mathbb{R}^{3}$ represents an isotropic homogeneous
material and we consider the following linear elastic boundary value problem: Find the displacement vector $u(\param)$ and the Cauchy stress tensor $\boldsymbol{\sigma}(\boldsymbol{u}(\param))$ such that
\begin{alignat}{2}\label{strong form}
\nonumber - \nabla\cdot \boldsymbol{\sigma}(\boldsymbol{u}(\param)) &= \mathbf{G}(\param) \quad &\text{in} \enspace \Omega,\\
\boldsymbol{\sigma}(\boldsymbol{u}(\param)) \cdot \boldsymbol{n} &= 0  \quad &\text{on} \enspace \Gamma_{N},\\
\nonumber \boldsymbol{u}(\param) &= \boldsymbol{g}_{D} \quad &\text{on} \enspace \Gamma_{D},
\end{alignat}
where the body force $\boldsymbol{G}: \mathcal{P} \rightarrow \mathbb{R}^{3}$ accounts for gravity. We can express for a linear elastic material the Cauchy stress tensor as $\boldsymbol{\sigma}(\boldsymbol{u}(\param)) = E(\param)\, \boldsymbol{C} : \boldsymbol{\varepsilon}(\boldsymbol{u}(\param))$, where $\boldsymbol{C}$ is the fourth-order stiffness tensor, $\boldsymbol{\varepsilon}(\boldsymbol{u}(\param)) = 0.5 (\nabla \boldsymbol{u}(\param) + (\nabla \boldsymbol{u}(\param))^{T})$ is the infinitesimal strain tensor, and the colon operator $:$ is defined as $\boldsymbol{C} : \boldsymbol{\varepsilon}(\boldsymbol{u}(\param)) = \sum_{k,l=1}^{3}\boldsymbol{C}_{ijkl}\boldsymbol{\varepsilon}_{kl}(\boldsymbol{u}(\param))$. Moreover, $E:\mathcal{P} \rightarrow L^{\infty}(\Omega)$ denotes Young's modulus, which is assumed to be piecewise constant on $\Omega$ and satisfy $E(\param) \geq E_{0} >0$ for a constant $E_{0} \in \mathbb{R}^{+}$. Therefore, the stiffness tensor can be written as 
\begin{equation*}
\boldsymbol{C}_{ijkl} = \frac{\nu}{(1+\nu)(1 -2\nu)}\delta_{ij}\delta_{kl} + \frac{1}{2(1 + \nu)}(\delta_{ik}\delta_{jl} + \delta_{il}\delta_{jk}), \quad 1 \leq i,j,k,l \leq 3,
\end{equation*}
where $\delta_{ij}$ denotes the Kronecker delta; we choose Poisson's ratio $\nu =0.3$. The corresponding variational formulation of \eqref{strong form} then reads as follows: For any $\param \in \mathcal{P}$ find $u(\param) \in V = \{\boldsymbol{v} \in [H^{1}(\Omega)]^{3}\enspace : \enspace \boldsymbol{v} = 0 \enspace \text{on} \enspace \Gamma_{D}\}$ such that
\begin{equation}\label{exact}
a(\boldsymbol{u}(\param),\boldsymbol{v};\param) = f(\boldsymbol{v};\param) \quad \forall \boldsymbol{v} \in V.
\end{equation}
Here, the bilinear and linear forms $a(\cdot, \cdot; \param): [H^{1}(\Omega)]^{3} \times  [H^{1}(\Omega)]^{3} \rightarrow \mathbb{R}$ and $f(\cdot;\param): [H^{1}(\Omega)]^{3} \rightarrow\mathbb{R}$ are defined as
\begin{equation*}
a(w,v;\param) := \int_{\Omega} E(\param) \,\frac{\partial \boldsymbol{w}^{i}}{\partial x_{j}} \boldsymbol{C}_{ijkl} \frac{\partial \boldsymbol{v}^{k}}{\partial x_{l}}  \enspace \text{and} \enspace f(\boldsymbol{v};\param) := \int_{\Omega} \boldsymbol{G}(\param)\,\cdot v \, - a(\widehat{\boldsymbol{G}}(\param),\boldsymbol{v};\param),
\end{equation*}
where $\widehat{\boldsymbol{G}}(\param) \in  [H^{1}(\Omega)]^{3}$ denotes a suitable lifting function of the possibly non-homogeneous Dirichlet boundary conditions.
\end{example}

To obtain approximate solutions of \eqref{eq:weak_solution} we presume we have an appropriate grid-based numerical method at hand (the full order model, FOM), yielding a high (but finite) dimensional approximation space $V_h$.
We consider conforming continuous Galerkin finite elements (FE), where $V_h \subset V$, and nonconforming discontinuous Galerkin (DG) or finite volume (FV) schemes, where $V_h \not\subset V$ (in which case we require broken Sobolev spaces for our analysis, see Section \ref{sec:ip-nonconforming}).
As a starting point for localized model reduction we require the FOM space to be decomposable into ``local'' spaces, which we will make more precise shortly.
While a localizing space decomposition could in general stem from any clustering of the degrees of freedom (DoF) of $V_h$ (see for instance \cite{Car2015}), we are particularly interested in local approximation spaces which are associated with a domain decomposition of the physical domain.

\begin{definition}[Non overlapping domain decomposition]\label{def:decompoisition}
  We call a finite collection of $M \in \N$ open polygonal subdomains $\Grid := \big\{\Omega_1, \Omega_2, \dots, \Omega_{M}\big\}$ a non overlapping domain decomposition of the physical domain $\Omega$, if $\bigcup_{m=1}^{M}\overline{\Omega_m} = \overline{\Omega}$ and ${\Omega}_m \cap {\Omega}_{m'} = \emptyset$ for $1 \leq m, m' \leq M$, $m \neq m'$. 
  We collect in $\Grid^v$, $\Grid^e$ and $\Grid^{\gamma}$, the set of all vertices, edges and facets (which we will denote interfaces from now on)\footnote{Note that to simplify notation we denote both the upper bound of the continuity constant and the local interfaces with $\gamma$, expecting that the respective meaning will be clear from the context.}, respectively, associated with the partition $\Grid$ and define $H := \max_{m = 1}^{M} \diam \Omega_m$. Moreover, we denote by $\Gamma:= \left(\bigcup_{m=1}^{M} \partial\Omega_{m}\right)\setminus \partial\Omega$ the whole interface of the decomposition $\Grid$.
  Note that $\Grid^e = \emptyset$ for $d = 2$ and $\Grid^e = \Grid^{\gamma} = \emptyset$ for $d = 1$.
  Each of the sets $\Grid$, $\Grid^v$, $\Grid^e$ and $\Grid^{\gamma}$ can be decomposed into elements associated with the domain boundary and inner elements, and we collect the latter in $\mathring{\Grid}$, $\mathring{\Grid^v}$, $\mathring{\Grid^e}$ and $\mathring{\Grid^\gamma}$, respectively.
  For instance, for each two adjacent subdomains $\Omega_m, \Omega_{m'} \in \Grid$, there exists a shared interface $\Gamma_{m, m'} \in \mathring{\Grid^\gamma}$, while for all boundary subdomains $\Omega_m \in \mathring{\Grid}$ there exists at least one boundary interface $\Gamma_{m, \partial\Omega} \in \Grid^{\gamma}  \backslash \mathring{\Grid^{\gamma}}$. 
\end{definition}

We can thus think of the domain decomposition as a usual grid, but without the requirements of $\Grid$ to actually resolve any data functions of the PDE.
Given such a domain decomposition, we can abstractly define a localizing space decomposition.

\begin{definition}[Localizing space decomposition]
  \label{def:localizing_space_decomposition}
  Let the FOM space $V_h$ be a finite dimensional Hilbert space with inner product and induced norm $\|\cdot\|_{V_h}^2 = (\cdot, \cdot)_{V_h}$.
  We call the direct sum decomposition of $V_h$ into subdomain spaces, interface spaces, edge spaces and vertex spaces,
  \begin{align}
    V_h = \bigoplus_{m=1}^{M} V_h^{m} \;\;\oplus\;\; \bigoplus_{\gamma \in \Grid^{\gamma}} V_h^\gamma  \;\;\oplus\;\; \bigoplus_{e \in \Grid^e} V_h^e  \;\;\oplus\;\; \bigoplus_{v \in \Grid^v} V_h^v,
    \label{eq:localizing_space_decomposition}
  \end{align}
  a localizing space decomposition.
\end{definition}

Note that such a decomposition is not unique and can always be found. Since the reduced space shall inherit this localizing decomposition, its purpose will be threefold: (i) offline, it allows for an independent and localized generation of the local reduced approximation spaces (compare Section \ref{sec:prep}), (ii) it allows to define and alter the physical domain $\Omega$ online, given that local approximation spaces for certain reference subdomains have been prepared offline, and (iii) it allows to adapt a local approximation space online (by adding basis functions or changing the local grid), while only requiring an update of local and neighboring prepared quantities (compare Section \ref{sec:adaptivity}).
For actual examples of space decompositions we refer to Section \ref{sec:coupling}.

Abstractly, we do not impose any further assumptions on the FOM as well as the reduced order model (ROM).
However, given the (bi-)linearity of $a$ and $f$, the computational benefits of the localizing space decomposition are apparent (and are made more precise throughout the rest of this chapter).
Since we allow for non conforming approximations, in general we need to consider discrete counterparts of $a$ and $f$ which are only defined on the FOM space $V_h$ and not necessarily on $V$, where we again refer to the following sections for examples.

\begin{definition}[Locally decomposed full order model (FOM)]
  \label{def:fom}
  Let $V_h$ be locally decomposable as in definition \eqref{eq:localizing_space_decomposition}, and let $a_h(\cdot,\cdot;\param): V_h \times V_h \to \R$ and $f_h(\cdot;\param) \in V_h'$ denote discrete variants of $a$ and $f$, respectively, which are continuous and coercive w.r.t.~the inner product of $V_h$.
  For each $\param \in \Params$, find $u_h(\param) \in V_h$ such that
  \begin{align}
    a_h\big(u_h(\param), v_h; \param\big) = f_h(v_h; \param) &&\text{for all } v_h \in V_h.
    \label{eq:fom_solution}
  \end{align}
\end{definition}

The idea of projection-based localized model order reduction is to consider a local reduced approximation space for each element of the localizing space decomposition \eqref{eq:localizing_space_decomposition}, in order to obtain a similarly decomposed reduced space $V_N \subset V_h$:
\begin{align}
  V_N = \bigoplus_{m=1}^M V_N^{m} \;\;\oplus\;\; \bigoplus_{\gamma \in \Grid^{\gamma}} V_N^\gamma  \;\;\oplus\;\; \bigoplus_{e \in \Grid^e} V_N^e  \;\;\oplus\;\; \bigoplus_{v \in \Grid^v} V_N^v,
  \label{eq:localized_V_N}
\end{align}
with reduced subdomain spaces $V_N^{m} \subset V_h^{m}$, reduced interface spaces $V_N^\gamma \subset V_h^\gamma$, reduced edge spaces $V_N^e \subset V_h^e$ and reduced vertex spaces $V_N^v \subset V_N^v$.
Similar to standard projection based model order reduction, we obtain the ROM simply by Galerkin projection of the locally decomposed FOM \eqref{eq:fom_solution} onto this locally decomposed reduced space.

\begin{definition}[Locally decomposed reduced order model (ROM)]
  \label{def:rom}
  Given a locally decomposed reduced space as in \eqref{eq:localized_V_N}, for each $\param \in \Params$, find $u_N(\param) \in V_N$ such that
  \begin{align}
    a_h\big(u_N(\param), v_N; \param\big) = f_h(v_N; \param) &&\text{for all } v_N \in V_N.
    \label{eq:rom_solution}
  \end{align}
\end{definition}

The main questions remain: (i) how to choose good local reduced approximation spaces to guarantee accurate and at the same time efficient reduced order approximations, (ii) how to benefit from the localization of $V_N$, that is how to obtain an offline/online decomposed scheme and in particular how to couple these local reduced approximation spaces, and (iii) how to adaptively enrich these local reduced approximation spaces online, if required.
These topics will be answered throughout the remainder of this chapter, starting with examples of how to obtain localized FOMs from standard discretization schemes and how to couple the resulting local reduced spaces.

Therefore, we introduce local grids $\grid(\Omega_m)$ on each subdomain $\Omega_m \subset \Grid$, which we presume to resolve all data functions of the underlying PDE.
As an analytical tool, we also define the global fine grid by $\grid = \cup_{\Omega_m \in \Grid} \grid(\Omega_m)$, which is usually not required in practical computations.
For simplicity, we require the local grids of two subdomains $\Omega_m, \Omega_m' \in \Grid$ to match along the shared interface $\gamma_{m, m'} \in \mathring{\Grid^\gamma}$ and denote by $\faces(\gamma_{m, m'})$ the corresponding set of all facets of $\grid$ which lie on $\gamma_{m, m'}$. Finally, we require that $\Gamma$ does not cut any grid cells.


\section{Coupling local approximation spaces }\label{sec:coupling}

\subsection{Conforming approach} \label{subsec:conforming}

There are various ways to couple local reduced spaces such that we obtain a conforming approximation, such as the partition of unity method \cite{BabMel97} or the generalized finite element method (GFEM) \cite{BaCaOs94,BaBaOs04,BabMel97,BabLip11}. However, in this section we focus on the decomposition into interface spaces and intra-element spaces, where the coupling is performed via the coupling or interface modes spanning the interface space.

\subsubsection{The multidomain problem and the Steklov-Poincar\'{e} interface equation}\label{subsub:steklovpoincare}

First, we introduce local Hilbert spaces $H^{1}_{0}(\Omega_{m}) \subset V^{m} \subset H^{1}(\Omega_{m})$, $m=1,\hdots,M$, which are supposed to respect the boundary conditions on $\partial\Omega$, the local spaces $V_{0}^{m}:=\{v \in V^{m} \, : \, v|_{\partial\Omega_{m}\setminus \partial \Omega} = 0 \}$, and the trace space $\Lambda$ associated with $\Gamma$. Moreover, we introduce local parameter-dependent bilinear and linear forms $a_{m}(\cdot, \cdot;\param): V^{m} \times V^{m} \rightarrow \mathbb{R}$ and $f_{m}(\cdot;\param) \in V^{m\prime}$, $\param \in \mathcal{P}$, $m=1,\hdots,M$, and the inner product $(\cdot,\cdot)_{V^{m}}: V^{m} \times V^{m} \rightarrow \mathbb{R}$. We may then state the variational form \eqref{eq:weak_solution} equivalently as follows (see for instance \cite{QuaVal05}): For any $\param \in \mathcal{P}$ find $u_{m}(\param) \in V^{m}$, $m =1,\hdots,M$ such that
\begin{subequations}\label{eq:variational_form_multi}
\begin{alignat}{2}
\label{eq:weak_form_loc} a_{m}(u_{m}(\param), v;\param) &= f_{m}(v;\param) \quad &\forall v \in V_{0}^{m},\\
\label{eq:weak_form_loc2} u_{m}(\param) &= u_{m'}(\param) \quad &\text{on} \enspace \Gamma_{m,m'},\\
\label{eq:weak_flux_continuity} \sum_{m=1}^{M} a_{m}(u_{m}(\param),\mathcal{E}_{m}\zeta;\param) &= \sum_{m=1}^{M} f_{m}(\mathcal{E}_{m}\zeta;\param) \quad &\forall \zeta \in \Lambda,
\end{alignat}
\end{subequations}
where $\mathcal{E}_{m}:\Lambda \rightarrow V^{m}$, $m=1,\hdots,M$ are linear and continuous extension operators. 

The formulation \eqref{eq:variational_form_multi} can then be used to derive an equation that solely acts on functions on the interface but nevertheless uniquely defines the solution $u(\param)$ of \eqref{eq:weak_solution}. To that end, we introduce a parameter-dependent lifting operator $\mathcal{E}_{\Gamma \rightarrow \Omega}(\param):\Lambda \rightarrow V$, where $\mathcal{E}_{\Gamma \rightarrow \Omega}(\param)\zeta$ is defined as the minimizer of $\inf_{v(\param) \in V} a(v(\param),v(\param);\param)$ subject to $v(\param)|_{\Gamma}= \zeta$. Note that we then also have
\begin{equation}\label{eq:def_local_lifting}
a_{m}(\mathcal{E}_{\Gamma \rightarrow \Omega}(\param) \zeta, v; \param) = 0 \quad \forall v \in V^{m}_{0} \quad \text{and} \quad \mathcal{E}(\param)_{\Gamma \rightarrow \Omega}\zeta = \zeta \quad \text{on} \enspace \Gamma \cap \partial \Omega_{m}.
\end{equation}
Then, we can rewrite the solution $u(\param)$ as 
\begin{equation}\label{eq:alternative_rep_solution}
u(\param) = \mathcal{E}_{\Gamma \rightarrow \Omega}(\param)(u(\param)|_{\Gamma}) + \sum_{m=1}^{M} u_{m}^{f}(\param),
\end{equation}
where $u_{m}^{f}(\param) \in V_{0}^{m}$ solves
\begin{equation}
a_{m}(u_{m}^{f}(\param), v;\param) = f_{m}(v;\param) \quad \forall v \in V_{0}^{m}, \enspace m=1,\hdots,M.
\end{equation}
Inserting \eqref{eq:alternative_rep_solution} into \eqref{eq:weak_flux_continuity} and choosing $\mathcal{E}_{m}=\mathcal{E}_{\Gamma \rightarrow \Omega_{m}}(\param)$
yields the Steklov-Poincar\'{e} interface equation: Find $u(\param)|_{\Gamma} \in \Lambda$ such that 
\begin{align}
\nonumber &\sum_{m=1}^{M} a_{m}(\mathcal{E}_{\Gamma \rightarrow \Omega_{m}}(\param)(u(\param)|_{\Gamma}), \mathcal{E}_{\Gamma \rightarrow \Omega_{m}}(\param)\zeta;\param) \\[-1.5ex]
&\label{eq:Steklov-Poincare}\\[-1.5ex]
\nonumber &\quad =\sum_{m=1}^{M} \left[ f_{m}( \mathcal{E}_{\Gamma \rightarrow \Omega_{m}}(\param)\zeta;\param) - a_{m}(u_{m}^{f}(\param),  \mathcal{E}_{\Gamma \rightarrow \Omega_{m}}(\param)\zeta;\param) \right] \quad \forall \zeta \in \Lambda.
\end{align}
Let us notice that the Steklov-Poincar\'{e} interface equation and its discrete, algebraic analogon, the Schur complement system, are at the base of iterative substructuring methods (see \cite{QuaVal05,TosWid05}), which have been combined with the reduced basis method in \cite {MaiHaa14}.  

We may finally define a space associated with the interface $V^{\Gamma}$ as
$V^{\Gamma}=\{ \mathcal{E}_{\Gamma \rightarrow \Omega}(\param)\zeta \in V\, : \, \zeta \in \Lambda\}$ and obtain the decomposition $V = \left( \bigoplus_{m=1}^{M}V_{0}^{m} \right) \oplus V^{\Gamma}$. This decomposition is $a$-orthogonal thanks to \eqref{eq:def_local_lifting}.

While the computation of the (harmonic) lifting operators is inherently local (see \eqref{eq:def_local_lifting}), the Steklov-Poincar\'{e} interface equation is posed on the whole interface $\Gamma$. To localize the latter we decompose $V^{\Gamma}$ as we will describe next.

\subsubsection{A conforming, localized reduced order approximation}

First, we determine basis functions associated with the vertices $v \in \mathcal{T}_{H}^{v}$. One common approach \cite{HetLeh10,JaBeLa11,Betal17} is to require that a basis function $\psi^{v} \in V_{h}\cap [H^{1}_{0}(\underset{v \subset \overline{\Omega}_{m}}{\bigcup_{m}} \overline{\Omega}_{m})]^{z}$, $z=1,\hdots,d$ associated with some vertex $v$ of the coarse mesh $\mathcal{T}_{H}$ satisfies for all $\Omega_{m}$, $m=1,\hdots,M$:
\begin{equation}\label{eq:def_vertice_basis_function}
(\psi^{v}, w)_{V^{m}} = 0 \quad \forall w \in V^{m}_{h;0} \quad \text{and} \quad \psi^{v}(\mathbf{x}^{v}) = 1, \enspace \psi^{v}(\mathbf{x}^{v'}) = 0, v \neq v'.
\end{equation}
Here, $\mathbf{x}^{v}$ are the (global) coordinates of the vertex $v$ and $V^{m}_{h;0}:=\{v \in V^{m}_{h}\,:\, v=0 \text{ on } \partial\Omega_{m}\setminus\Gamma_{N}\}$ . To uniquely define $\psi^{v}$ we need to prescribe the respective values on $\Gamma$. We may for instance require $\psi^{v}$ to be linear on the respective edges or bilinear on the respective interfaces (see e.~g. \cite{Betal17}). For multiscale problems in two space dimensions with a permeability $\kappa(\mathbf{x}_{1},\mathbf{x}_{2};\paramFixed)$ it has been suggested in \cite{HouWu97} to prescribe 
\begin{equation}
\psi^{v}(\mathbf{x}_{1},\mathbf{x}_{2}^{v}) := \left(\int_{\mathbf{x}_{1}^{v'}}^{\mathbf{x}_{1}} \frac{ds}{\kappa(s,\mathbf{x}_{2}^{v};\paramFixed)} \right) \bigl / \left( \int_{\mathbf{x}_{1}^{v'}}^{\mathbf{x}_{1}^{v}} \frac{ds}{\kappa(s,\mathbf{x}_{2}^{v};\paramFixed)} \right)
\end{equation} 
on a horizontal edge $[\mathbf{x}_{1}^{v'},\mathbf{x}_{1}^{v}] \times \{\mathbf{x}_{2}^{v}\}$  in a uniform rectangular coarse grid $\mathcal{T}_{H}$.

Next, we assume that we have given sets of discrete edge basis functions $\{\chi^{e}_{k}\}_{k=1}^{N^{e}_{h;0}} \in V_{h}|_{e}$ and discrete interface basis functions $\{\chi^{\gamma}_{k}\}_{k=1}^{N^{\gamma}_{h;0}} \in V_{h}|_{\gamma}$ defined on the respective edge $e \in \mathcal{T}_{H}^{e}$ or interface $\gamma \in \mathcal{T}_{H}^{\gamma}$. Here, we set $N^{e}_{h;0}:=\dim(V_{h}|_{e\setminus \partial e})$ and $N^{\gamma}_{h;0}:=\dim(V_{h}|_{\gamma\setminus \partial \gamma})$ as we require that $\chi^{e}_{k}$ and $\chi^{\gamma}_{k}$ are zero on the boundary of the edge and interface, respectively. Furthermore, we define $\Lambda_{N^{e}_{h;0}}^{e}:=\spanlin\{\chi^{e}_{1},\hdots,\chi^{e}_{N^{e}_{h;0}}\}$ and $\Lambda_{N^{\gamma}}^{\gamma}:=\spanlin\{\chi^{\gamma}_{1},\hdots,\chi^{\gamma}_{N^{\gamma}_{h;0}}\}$. 

Similarly as for the vertices we may then define associated basis functions that have support on the union of subdomains that share the respective edge or interface: Find $\psi^{\gamma}_{k} \in V_{h} \cap [H^{1}_{0}(\underset{\gamma \subset \overline{\Omega}_{m}}{\bigcup_{m}}\overline{\Omega}_{m})]^{z}$, $z=1,\hdots,d$, $\gamma \in \mathcal{T}_{H}^{\gamma}$, $k=1,\hdots,N^{\gamma}_{h;0}$ such that
\begin{equation}\label{eq:def_face_basis_function}
(\psi^{\gamma}_{k}, w)_{V^{m}} = 0 \quad \forall w \in V^{m}_{h;0} \quad \text{and} \quad \psi^{\gamma}_{k}|_{\gamma} = \chi^{\gamma}_{k}.
\end{equation}
Likewise, we find $\psi^{e}_{k} \in V_{h}\cap [H^{1}_{0}(\underset{e \subset \overline{\Omega}_{m}}{\bigcup_{m}} \overline{\Omega}_{m})]^{z}$, $z=1,\hdots,d$, $e \in \mathcal{T}_{H}^{e}$, $k=1,\hdots,N^{e}_{h;0}$ such that
\begin{equation}\label{eq:def_edge_basis_function}
(\psi^{e}_{k}, w)_{V^{m}} = 0 \quad \forall w \in V^{m}_{h;0} \quad \text{and} \quad \psi^{e}_{k}|_{e} = \chi^{e}_{k}.
\end{equation}
Again, we need to provide the value of $\psi^{e}_{k}$ on the interfaces sharing the edge $e \in \mathcal{T}_{H}^{e}$ in order to uniquely define $\psi^{e}_{k}$. Similarly to above we may require that $\psi^{e}_{k}$ is linear on the respective interfaces as suggested for instance in \cite{Betal17} or define a function which takes into account also the coefficient function. 

Note that if the interfaces are mutually disjoint, which is for instance the case if we associate the subdomains $\Omega_{m}$, $m=1,\hdots,M$ with the components of a structure, only the basis functions $\psi^{\gamma}_{k}$, $k=1,\hdots,N^{\gamma}_{h;0}$ ($d=3$) or $\psi^{e}_{k}$, $k=1,\hdots,N^{e}_{h;0}$ ($d=2$) are needed. Here, the values of the basis functions on the boundary of the interfaces or edges are determined by the boundary conditions on $\partial\Omega$ (see for instance \cite{HuKnPa13,HuKnPa13b,EftPat13a,SmePat16}).

For $N^{\gamma} \ll N^{\gamma}_{h;0}$ and $N^{e} \ll N^{e}_{h;0}$ may now define the reduced space associated with $\Gamma$ as follows:
\begin{align}\label{eq:reduced interface space}
V^{\Gamma}_{N} &:= \bigoplus_{v \in \mathcal{T}_{H}^{v}} \spanlin\{\psi^{v}\} \oplus \bigoplus_{e \in \mathcal{T}_{H}^{e}} \spanlin\{\psi^{e}_{1},\hdots,\psi^{e}_{N^{e}}\} \oplus \bigoplus_{\gamma \in \mathcal{T}_{H}^{\gamma}} \spanlin\{\psi^{\gamma}_{1},\hdots,\psi^{\gamma}_{N^{\gamma}}\}.
\end{align}
Such reduced interface spaces are for instance employed in (adaptive) component mode synthesis (CMS) \cite{HetLeh10,JaBeLa11}, the static condensation reduced basis element (scRBE) method  \cite{HuKnPa13,HuKnPa13b,EftPat13a,SmePat16} for mutually disjoint interfaces, or in the ArbiLoMod \cite{Betal17}. Subspaces of $V_{\Gamma}^{N}$ are considered in certain multiscale methods. For example in the  MsFEM of Hou and Wu  \cite{HouWu97} the reduced space is spanned by the basis functions $\psi^{v}$, $v \in \mathcal{T}_{H}^{v}$. For further relations between CMS, MsFEM and GFEM we refer e.g. to \cite{HetLeh10}.

Recall that the basis functions associated with the vertices, edges, and interfaces have all been computed w.r.t. an inner product that does not depend on the parameter (see \eqref{eq:def_vertice_basis_function}, \eqref{eq:def_face_basis_function}, \eqref{eq:def_edge_basis_function}). Therefore, we finally assume that we have also given reduced spaces $V_{N;0}^{m}:= \spanlin\{\zeta_{1}^{m},\hdots,\zeta_{N^{m}}^{m}\} \subset V^{m}_{h;0}$, $m=1,\hdots,M$, that will account for parameter variations. In detail, we obtain approximations $\tilde{\psi}^{*}_{k}(\param)$, $*=v,e,\gamma$ by solving
\begin{equation}\label{eq:def_intra_element_correction}
\text{find } \tilde{b}_{k}^{*}(\param) \in V^{m}_{N;0}:  \quad a_{m}(\psi^{*}_{k} + \tilde{b}_{k}^{*}(\param), w; \param) = 0 \quad \forall w \in V^{m}_{N;0}
\end{equation}
and setting $\tilde{\psi}^{*}_{k}(\param)= \psi^{*}_{k} + \tilde{b}_{k}^{*}(\param)$, $*=v,e,f$. Finally, we define $\tilde{b}^{m}(\param) \in V^{m}_{N;0}$ as the solution of
\begin{equation}\label{eq:def_intra_element_correction_source}
\text{find } \tilde{b}^{m}(\param) \in V^{m}_{N;0}:  \quad a_{m}(\tilde{b}^{m}(\param), w; \param) = f(w,\param) \quad \forall w \in V^{m}_{N;0}.
\end{equation}
Note, that both $\tilde{b}_{k}^{*}(\param)$, $*=v,e,\gamma$ and $\tilde{b}^{m}(\param)$ can be interpreted as intra-element RB approximations. The corresponding reduced spaces $V_{N;0}^{m}$, $m=1,\hdots,M$ can for instance be constructed from solutions $b_{k}^{*}(\param),b^{m}(\param) \in V^{m}_{h;0}$, $*=v,e,f$ of
\begin{equation}\label{eq:def_intra_element_snapshot}
\quad a_{m}(\psi^{*}_{k} + b_{k}^{*}(\param), w; \param) = 0 \quad \forall w \in V^{m}_{h;0}
\end{equation}
and 
\begin{equation}\label{eq:def_intra_element_snapshot_source}
\quad a_{m}(b^{m}(\param), w; \param) = f(w,\param) \quad \forall w \in V^{m}_{h;0},
\end{equation}
respectively, via a standard greedy algorithm or a POD.\footnote{Note that in actual practice one would construct the reduced bases only on a certain number $<M$ of reference domains; see for instance \cite{HuKnPa13}.} Let us also remark that for instance in the scRBE method for the approximation of each basis function $\psi^{*}_{k}$, $*=v,e,\gamma$ a different RB space is considered, to further reduce the size of problems \eqref{eq:def_intra_element_correction}, \eqref{eq:def_intra_element_correction_source}. Finally, we define the reduced spaces
\begin{equation}\label{eq:reduced space conform}
V_{N}=\bigoplus_{m=1}^{M} V^{m}_{N;0} \oplus V^{\Gamma}_{N}
\end{equation}
and
\begin{align}
\nonumber V^{\Gamma}_{N}(\param) &:=\bigoplus_{v \in \mathcal{T}_{H}^{v}} \spanlin\{\tilde{\psi}^{v}(\param)\} \oplus \bigoplus_{e \in \mathcal{T}_{H}^{e}} \spanlin\{\tilde{\psi}^{e}_{1}(\param),\hdots,\tilde{\psi}^{e}_{N^{e}}(\param)\} \\[-1.5ex]
 &\label{eq:reduced interface space parameter}\\[-1.5ex]
\nonumber &\qquad\qquad \oplus \bigoplus_{\gamma \in \mathcal{T}_{H}^{\gamma}} \spanlin\{\tilde{\psi}^{\gamma}_{1}(\param),\hdots,\tilde{\psi}^{\gamma}_{N^{\gamma}}(\param)\}.
\end{align}
The global reduced approximation $u_{N}(\param)$ can then be computed by performing a Galerkin projection onto the reduced space $V^{\Gamma}_{N}(\param)$ or a Petrov-Galerkin approximation using $V^{\Gamma}_{N}(\param)$ as a trial and $V^{\Gamma}_{N}$ as a test space (see e.g. \cite{EftPat13b,Sme15}). Instead of eliminating the volume degrees of freedom via \eqref{eq:def_intra_element_correction}, \eqref{eq:def_intra_element_correction_source}, $u_{N}(\param)$ can also directly be determined by performing a Galerkin projection onto $V_N$ (see for instance \cite{Betal17}). Similarly, for CMS and a fixed parameter a Galerkin projection onto $V_{N}$ may be performed to compute the reduced solution; here, the reduced space $V_{N;0}^{m}$ is constructed from an eigenvalue problem and does not account for parameter variations (see e.g. \cite{HetLeh10}). Finally, in the reduced basis---domain decomposition---finite element (RDF) method \cite{IaQuRo16} the reduced space $V_{N}$ is chosen as a direct sum of $\bigoplus_{m=1}^{M} V^{m}_{N;0}$ and standard FE spaces defined on the interface or on a (small) area around the interface. Here, the intra-element reduced spaces $V^{m}_{N;0}$ are constructed via a greedy algorithm considering a parametrized linear combination of standard Lagrange basis functions or Fourier modes as boundary conditions. Then, a Galerkin projection on $V_{N}$ is performed to compute $u_{N}(\param)$.

\subsection{Non-conforming approach}\label{subsec:non-conforming}
With the term {\em non-conforming approach} we want to classify a set of alternative techniques to solve the reduced problem on the global computational domain.
A first approach consists in considering a global system of equations given by local parametrized problems and additional equations ensuring the matching between the different subdomains through the use of Lagrange multipliers. This approach has been used for solving elliptic equation in \cite{MadRon02,MadRon04} and Stokes equations in \cite{LovMadRon06,IaQuRo12}.

Another approach consists in coupling local FOM spaces by interior penalty bilinear forms, inspired by discontinuous Galerkin FEM.
Here, we refer to the discontinuous Galerkin reduced basis element method \cite{CHM11,APQ2016,PaGeQu16} and 
the local reduced basis discontinuous {G}alerkin approach \cite{KOH2011}. A discontinuous Galerkin approach with local POD 
modes has been presented in \cite{FeIoLa18}.
In the context of multiscale problems (cf. Example \ref{ex:multiscale_problem}), the generalized multiscale discontinuous Galerkin 
method has been proposed in \cite{CEL2017,CEM2018} and used for solving the heat problem with phase change in \cite{SVV2018}.
In this chapter we are going to present the localized reduced basis method (LRBMS) in Subsection \ref{sec:ip-nonconforming} below.
LRBMS has been introduced in \cite{AlbrechtHaasdonkEtAl2012} and analyzed in \cite{OS2014, OS2015} for elliptic and 
in \cite{ORS2017} for parabolic problems. Applications to the simulation of two phase flow in porous media have been addressed in 
\cite{KFHLO2015} and to battery simulation with resolved electrodes in \cite{OR17}. 

\subsubsection{Non-conforming coupling based on Lagrange multipliers}\label{subsub:lagrange}

We want to reformulate the problem \eqref{eq:variational_form_multi}, with the idea that exact coincidence of the traces of the discrete functions (equation \eqref{eq:weak_form_loc2}) is generally too stringent, and may, in fact, lead to imposing $u_{m}=0$ on the internal interfaces; thus, the gluing process can be done in a dual way through Lagrange multipliers.
We assume that local basis functions are computed in each subdomain $\Omega_m, m=1,\dots, M$  by solving local parametrized variational problems coming from the original problem \eqref{eq:weak_solution} with proper boundary conditions along the boundaries which correspond to internal ones in the original domain. The choice of the boundary conditions is strongly related to the problem aimed to be solved. Thus, local reduced spaces are defined via these local solutions and denoted by $V_{N}^{m}, m=1,\dots, M$. Possible ways to construct $V_{N}^{m}$ are presented in section \ref{sec:prep}. If two or more subdomains are characterized by the same type of parameter and the same type of boundary conditions, the same local reduced space can be associated to those subdomains. For simplicity we consider different spaces for each different subdomain. \\
We define the following operator:
\begin{align}\label{LM}
\mathcal L^{m,m'}(u(\muv),\psi)= \int_{\Gamma_{m,m'}} (u(\boldsymbol{\mu})|_{\Omega_m}-u(\boldsymbol{\mu})|_{\Omega_{m'}}) \psi ds=0, \quad \forall \psi \in W_{m,m'},
\end{align}
where $ m,m' \in \{1,\cdots, M\}$, $\Gamma_{m,m'}$ is the interface between two adjacent sub-domains denoted with the indices $m$ and $m'$ respectively, and $ W_{m,m'}$ is the Lagrange multiplier space defined on this interface. Typical choices for the latter are low-order polynomial spaces \cite{MadRon02,IaQuRo12} or spaces constructed from snapshots (and their derivatives) \cite{MadRon04}. 

A basis for $ W_{m,m'}$ can then for instance be provided by the characteristic Lagrange polynomials $\psi_q$, $q=1,\cdots,Q_{m,m'}$ associated with the $Q_{m,m'}$ nodes of $\Gamma_{m,m'}$.\\
{If we suppose that $\Omega$ has $M-1$ internal interfaces, $\Gamma_{m,m+1}, m=1,\cdots,M-1$, the reduced global problem of this approach reads: find $u_N(\boldsymbol \muv)\in V_{N}^{1}\times\cdots\times V_{N}^{M},\lambda_N \in W_{m, m+1}, m=1,\cdots,M-1,$ such that
\begin {equation}\label{LMproblem}
\begin{cases} 
a( u_N(\muv),v_N,\muv)+\displaystyle\sum_{i=1}^{M-1}\mathcal L^{m, m+1}(v_N,  \lambda_N)=f( w,\muv)\forall v_N\in V_{N}^{1}\times\cdots\times V_{N}^{M},\\
\mathcal L^{m,m+1}(u_N(\muv),  \psi)=0\hspace{2.5cm}  m=1,\ldots,M-1, \forall   \psi \in W_{m,m+1}.
\end{cases}
\end{equation}

\subsubsection{Non-conforming coupling based on interior penalties}\label{subsec:ipdg}
\label{sec:ip-nonconforming}

We demonstrate how to obtain a localized FOM by applying ideas from interior penalty (IP) DG schemes w.r.t.~the domain decomposition $\Grid$ in the context of the parametric multiscale Example \ref{ex:multiscale_problem}.
To define the localized FOM, we presume we are given a discretization on the full global grid $\grid$ (which is not used in actual computations), which we make precise by specifying the approximation space with an associated inner product and discrete variants of $a$ and $f$.
As a common ground for the analysis of conforming as well as non conforming schemes, we introduce the broken Sobolev space $H^s(\grid(\omega)) := \big\{ v \in L^2(\omega) \;\big|\; v|_t \in H^s(t)\;\;\forall t \in \grid(\omega)\big\}$, for a given grid $\grid(\omega)$ of some domain $\omega \subseteq \Omega$ and $s \geq 1$, and associated broken gradient operator $\gradienth : H^1(\grid(\omega)) \to L^2(\omega)^d$ by $(\gradienth v)|_t := \gradient(v|_t)$ on all $t \in \grid$ for $v \in H^1(\grid(\omega))$.

\begin{example}[Continuous Galerkin (CG) FEM]
  The CG FEM scheme for the conforming approximation of Example \ref{ex:multiscale_problem} w.r.t.~the full global grid $\grid$ is given by
  the conforming approximation space of order $k \geq 1$,
    \begin{align}
      V_h^\text{CG}(\grid) := \big\{ v \in V \;\big|\; v|_t \in \mathbb{P}_k(t) \;\; \forall t \in \grid \big\} \;\subset\; V \;\subset\; H^1(\grid),
      \notag
    \end{align}
    where $\mathbb{P}_k(\omega)$ for any $\omega \subset \Omega$ denotes the space of all polynomials defined on $\omega$ of degree at most $k \geq 0$;
    the bilinear form $(\cdot, \cdot)^\text{CG}: H^1(\grid) \times H^1(\grid) \to \R$, given by $(u, v)^\text{CG} := \int_\Omega \gradienth u \cdot \gradienth v \dx$, as the inner product on $V_h^\text{CG}(\grid)$ (where we note that its restriction to $V \subset H^1(\grid)$ coincides with the $V$-inner product); and
    the discrete bilinear form $a_h^\text{CG}(\cdot,\cdot;\param): H^1(\grid) \times H^1(\grid) \to \R$ and linear functional $f_h^\text{CG}: H^1(\grid) \to \R$, given by
    \begin{align}
      a_h^\text{CG}(u, v; \param) := \int_\Omega \big(\kappa(\param) \gradienth u\big)\cdot \gradienth v \dx &&\text{and}&& f_h^\text{CG}(v) := \int_\Omega q v \dx
    \notag
    \end{align}
    (again noting that their respective restriction to $V$ coincide with $a$ and $f$).
\end{example}

The definition of the non conforming scheme is more involved.
We denote the set of faces of $\grid$ by $\faces$ and to each face $\sigmafine \in \faces$, we assign a unique normal $n_\sigmafine \in \R^d$ pointing away from $t^+$, where the face may be either an inner face $\sigmafine \in \faces$, given by the intersection of two grid elements $t^+, t^- \in \grid$, $\sigmafine = \overline{t^+ \cap t^-}$, or a boundary face $\sigmafine \in \faces$, given by $\sigmafine = \overline{t^+ \cap \partial\Omega}$ for some $t^+ \in \grid$.
Since functions in the broken Sobolev space are two-valued on grid faces, we introduce the mean $\mean{\cdot}$ and jump $\jump{\cdot}$ on a boundary face by $\mean{v} := \jump{v} := v|_{t^+}$ and by $\mean{v} := \tfrac{1}{2}(v|_{t^+} + v|_{t^-})$ and $\jump{v} := v|_{t^+} - v|_{t^-}$, respectively, on any other face.

Considering the family of interior penalty discontinuous Galerkin (DG) schemes, we present the symmetric variant for ease of notation, and refer to the symmetric weighted variant \cite{ESZ2009}, which is particularly well suited for multi-scale problems with highly varying or high-contrast coefficients.

\begin{example}[Interior penalty (IP) discontinuous Galerkin (DG) FEM]
  \label{ex:ipdg}
  The symmetric IPDG FEM scheme for the non conforming approximation of Example \ref{ex:multiscale_problem} w.r.t.~the full global grid $\grid$ is given by
  the non conforming approximation space of order $k \geq 1$,
    \begin{align}
      V_h^\text{DG}(\grid) := \big\{ v \in L^2(\Omega) \;\big|\; v|_t \in \mathbb{P}_k(t) \;\; \forall t \in \grid \big\} \;\subset\; H^1(\grid);
      \notag
    \end{align}
    the bilinear form $(\cdot, \cdot)^\text{DG}: H^1(\grid) \times H^1(\grid) \to \R$, given by
    \begin{align}
      (u, v)^\text{DG} := (u, v)^\text{CG} + \sum_{\sigmafine \in \faces} (u, v)_\sigmafine^p &&\text{with}&& (u, v)_\sigmafine^p := \int_\sigmafine h_\sigmafine^{-1} \jump{u}\jump{v} \ds,
    \notag
    \end{align}
    as inner product on $V_h^\text{DG}(\grid)$, where $h_\sigmafine$ is a positive number associated with each face $\sigmafine \in \faces$, e.g., $h_\sigmafine := \diam(\sigmafine)$ for $d \geq 2$ and $h_\sigmafine := \min\{\diam(t^+), \diam(t^-)\}$ for $d = 1$; and
    the linear functional $f_h^\text{DG}: H^1(\grid) \to \R$ given by $f_h^\text{DG}(v) := f_h^\text{CG}(v)$ and the discrete bilinear form $a_h^\text{DG}(\cdot, \cdot; \param): H^2(\grid) \times H^2(\grid) \to \R$, given by
    \begin{align}
      a_h^\text{DG}(u, v; \param) := a_h^\text{CG}(u, v; \param) + \sum_{\sigmafine \in \faces} a_\sigmafine(u, v; \param)
    \notag
    \end{align}
    with the face bilinear form $a_\sigmafine$ for any $\sigmafine \in \faces$ given by
    \begin{align}
      a_\sigmafine(v, u; \param) := a_\sigmafine^c(v, u; \param) + a_\sigmafine^c(u, v; \param) + (u, v)_\sigmafine^p\, w_\sigmafine
    \notag
    \end{align}
    with $a_\sigmafine^c(u, v; \param) := \int_\sigmafine-\mean{\big(\kappa(\param)\gradienth v\big)\cdot n_\sigmafine}\jump{u} \ds$ and a user-dependent penalty weight $w_\sigmafine > 0$, such that $a_h^\text{DG}$ is continuous and coercive w.r.t.~the above inner product.
\end{example}

The main idea of an \emph{IP localized FOM} is to consider the restriction of either of the above discretization schemes to each subdomain of the domain decomposition, and to again couple those with IP techniques along the interfaces of the subdomain.
We thus choose $* \in \{\text{CG}, \text{DG}\}$ and obtain the localized FOM space in the sense of Definition \ref{def:localizing_space_decomposition} as a direct sum of subdomain spaces (with empty interface, edge and vertex spaces)
\begin{align}
  V_h := \bigoplus_{m=1}^M V_h^{m}, &&\text{with}&& V_h^{m} := \big\{ v|_{\Omega_m} \;|\; v \in V^* \big\},
\notag
\end{align}
with associated inner product $(\cdot, \cdot)_{V_h}: V_h \times V_h \to \R$ given by
\begin{align}
  (u, v)_{V_h} := \sum_{m=1}^M (u|_{\Omega_m}, v|_{\Omega_m})^* + \sum_{\Gamma' \in \mathring{\Grid^\gamma}} \sum_{\sigmafine \in \faces(\Gamma')} (u, v)_\sigmafine^p.
\notag
\end{align}
We also define the linear functional $f_h: V_h \to \R$ by $f_h := f_h^*$ and, in a similar manner as above, the non conforming bilinear form $a_h(\cdot,\cdot;\param): V_h \times V_h \to \R$ by
\begin{align}
  a_h(u, v; \param) := \sum_{m=1}^M a_h^*(u|_{\Omega_m}, v|_{\Omega_m}; \param) + \sum_{\Gamma' \in \mathring{\Grid^\gamma}} \sum_{\sigmafine \in \faces(\Gamma')} a_\sigmafine(u, v; \param).
  \notag
\end{align}
We have thus fully specified a localized FOM in the sense of Definition \ref{def:fom} and comment on two special cases: for $* = \text{CG}$ and a trivial domain decomposition of a single subdomain, $\Grid = \{\Omega\}$, we obtain the above standard CG FEM while for $* = \text{DG}$ the resulting FOM coincides with the above standard symmetric IPDG FEM.

To make the coupling more precise, we may rearrange the above terms to obtain a localization of $a_h$ w.r.t.~the domain decomposition in the sense of
\begin{align}
  a_h(u, v; \param) &= \sum_{m = 1}^M a_h^{m}(u, v; \param) + \sum_{\Gamma' \in \mathring{\Grid^\gamma}} a_h^{\Gamma'}(u, v; \param),
\notag\\
\intertext{%
  with the subdomain and interface bilinear forms
}
  a_h^{m}(u, v; \param) :&= a_h^*(u|_{\Omega_m}, v|_{\Omega_m}; \param) + \sum_{\Gamma' \in \mathring{\Grid^\gamma}\cap\Omega_m} \sum_{\sigmafine \in \faces(\Gamma')} a_\sigmafine(u|_{\Omega_m}, v|_{\Omega_m}; \param),
\notag\\
  a_h^{\Gamma'}(u, v; \param) :&= \sum_{\sigmafine \in \faces(\Gamma')} \big\{a_\sigmafine(u|_{\Omega^+}, v|_{\Omega^-}; \param) + a_\sigmafine(u|_{\Omega^-}, v|_{\Omega^+}; \param)\big\},
\notag
\end{align}
respectively, for all $1 \leq m \leq M$ and all $\Gamma' \in \mathring{\Grid^\gamma}$, with the subdomains $\Omega^+, \Omega^- \in \Grid$ sharing the interface $\Gamma'$.

Now given a local reduced space $V_N^{m} \subset V_h^{m}$ for each subdomain we obtain the locally decomposed broken reduced space in the sense of \eqref{eq:localized_V_N} by
\begin{align}
  V_N = \bigoplus_{m=1}^M V_N^{m} \;\;\subset V_h.
\notag
\end{align}
Using the above decomposition of $a_h$ into subdomain and interface contributions, we can readily observe that the locally decomposed ROM can be offline/online decomposed by local computations: namely by projection of the subdomain bilinear forms $a_h^{m}(\cdot,\cdot;\param)$ onto $V_N^{m}\times V_N^{m}$ and the interface bilinear forms $a_h^{\Gamma'}(\cdot,\cdot;\param)$ onto $V_N^{m} \times V_N^{n}$, with $1 \leq m, n \leq M$, such that $\Omega^+ = \Omega_m$ and $\Omega^- = \Omega_n$, respectively.

We thus obtain a sparse matrix representation of the resulting reduced system, with a sparsity pattern which coincides with the one from standard IPDG schemes.


\section{Preparation of local approximation spaces}\label{sec:prep}

Both, couplings that yield a conforming and non-conforming approximation require either reduced spaces $\Lambda_{N^{\gamma}}^{\gamma} \subset V_{h}|_{\gamma}$ for interfaces and/or edges $\Lambda_{N^{e}}^{e} \subset V_{h}|_{e}$ (see subsection \ref{subsec:conforming}) or reduced spaces $V_{N}^{m}$  (see subsection \ref{subsec:non-conforming}) or both. As the generation of edge basis functions can be done analogously to the construction of interface basis functions we restrict ourselves to the latter in order to simplify notation. To fix the setting we thus consider the task of finding a suitable reduced space either on a subdomain $\Omega_{m} \subsetneq \Omega_{out} \subset \Omega$ with $\dist(\Gamma_{out},\partial \Omega_{m}) \geq \rho > 0$, $\Gamma_{out}:=\partial\Omega_{out}\setminus\partial \Omega$ or an interface $\Gamma_{m,m'} \subset \partial \Omega_{m}$, where $\dist(\Gamma_{out},\Gamma_{m,m'})\geq \rho > 0$. Possible geometric configurations of the oversampling domain $\Omega_{out}$ are illustrated in Fig. \ref{fig:illustration geometry}.

\begin{figure}[t]
\begin{center}
              \includegraphics[height=0.15\textwidth]{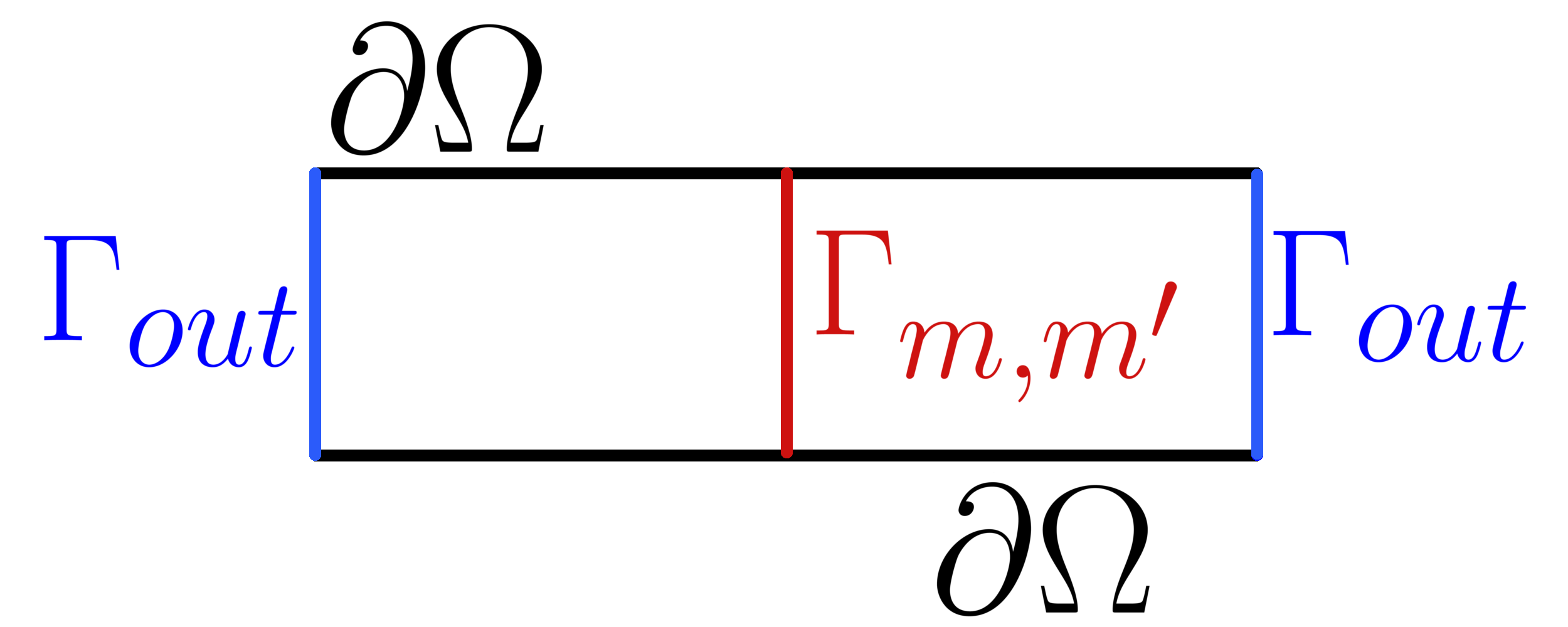} \quad
              \includegraphics[height=0.15\textwidth]{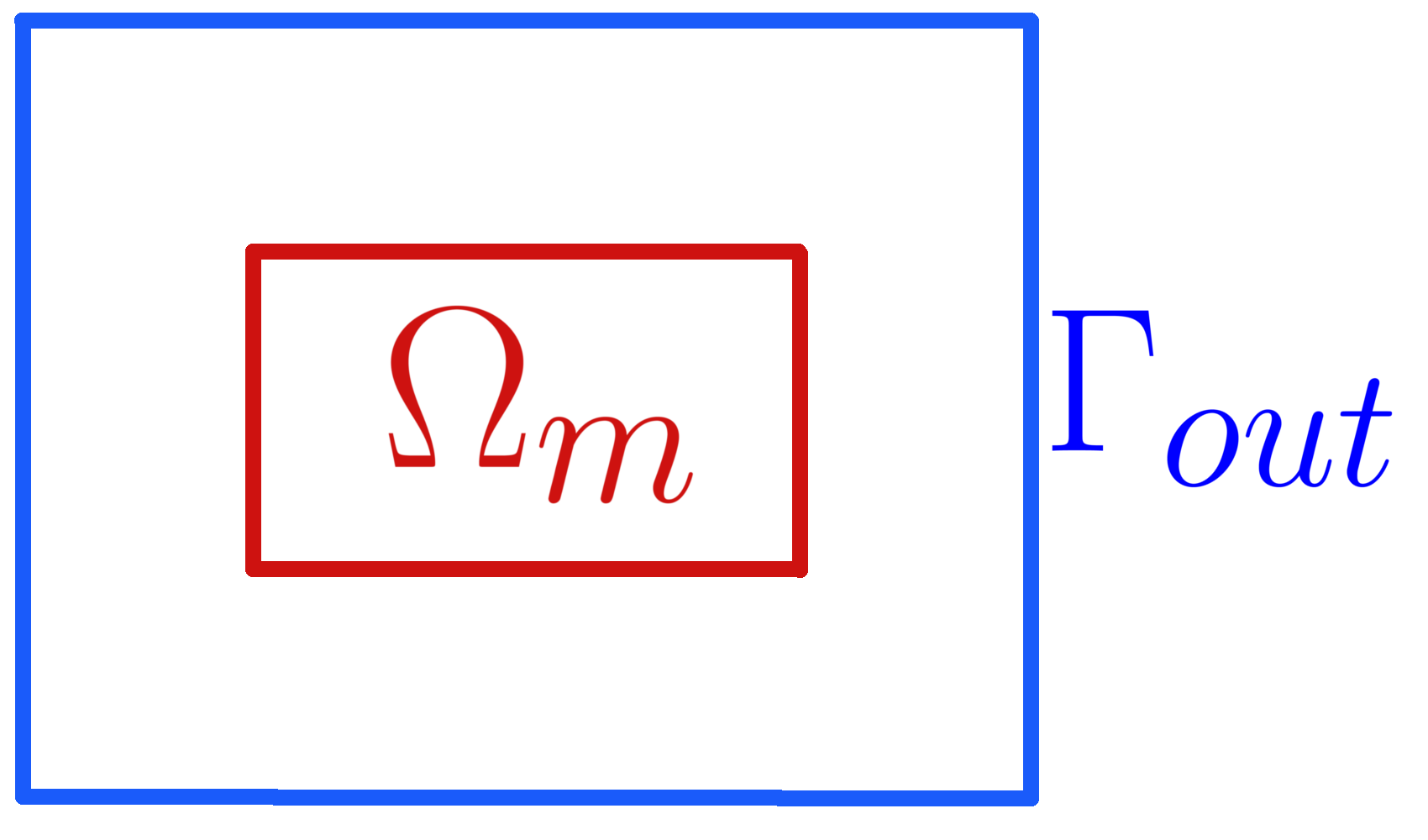}
              \caption{\footnotesize Illustration of possible decompositions of $\Omega_{out}$ with respect to $\Gamma_{m,m'}$ or $\Omega_{m}$.}\label{fig:illustration geometry}
              \end{center}
\vspace{-15pt}
\end{figure}

We will first briefly discuss in subsection \ref{subsec:bas_gen_polynomial} reduced spaces that are spanned by \emph{polynomials or solutions of ``standard'' eigenvalue problems} and are thus related to the spectral element method or hp-FEM. Subsequently, in subsection \ref{subsec:bas_gen_empirical} we will present reduced spaces that are generated from local solutions of the PDE, are thus of \emph{empirical} nature, and are optimal in the sense of Kolmogorov. We will also show how those optimal basis functions can be efficiently and accurately approximated by means of random sampling. 

\subsection{Polynomial-based local approximation spaces}\label{subsec:bas_gen_polynomial} 

Component mode synthesis (CMS) as introduced in \cite{Hur65,BamCra68} relies on free vibration modes or eigenmodes of local, constrained eigenvalue problems \cite{Hur65,BamCra68,Bou92,HetLeh10,JaBeLa11,HetKla14} for the approximation within subdomains. To couple the modes at the interfaces a reduced interface space spanned by eigenmodes is employed \cite{Hur65,BamCra68,Bou92,HetLeh10,JaBeLa11,HetKla14}. 

A combination of domain decomposition and RB methods has first been considered in the reduced basis element method (RBEM) \cite{MadRon02}. Here, inspired by the the mortar spectral element method \cite{BeMaPa94}, the Lagrange multiplier space $W_{m,m'}$ as defined in Subsection \ref{subsec:non-conforming} is chosen as a low-order polynomial space. The reduced basis hybrid method \cite{IaQuRo12} extends the RBEM by additionally considering a coarse FE discretization on the whole domain to account for continuity of normal stresses and also employs a low-order polynomial Lagrange multiplier space on the interface. For the scRBE method a reduced interface space spanned by the eigenvectors of a discrete generalized eigenvalue problem based on the Laplacian has been suggested in 
\cite{HuKnPa13,HuKnPa13b} and eigenmodes of a singular Sturm–Liouville eigenproblem
have been used in \cite{EftPat13a}. Finally, in the RDF method \cite{IaQuRo16} a standard FE space is considered on the interface or on a (small) area around the interface.

\subsection{Local approximation spaces based on empirical training}\label{subsec:bas_gen_empirical} 

In this subsection we are concerned with local approximation spaces that are constructed from local solutions of the PDE; those approaches are often called \emph{empirical}. Basis functions on the interfaces selected from local snapshots are for instance suggested in \cite{EftPat13a}, where an empirical pairwise training procedure for interface reduction within the scRBE context is developed, and within a heterogeneous domain decomposition method in \cite{MarRozHaa14}. Local approximation spaces that are optimal in the sense of Kolmogorov have been introduced for subdomains within the GFEM in \cite{BabLip11} for parameter-independent PDEs and for interfaces within static condensation procedures \cite{SmePat16} for parametrized PDEs. While the authors of \cite{SmePat16} introduce and analyze a spectral greedy algorithm to deal with parameter variations, \cite{TadPat18} suggests using a POD making use of the hierarchical approximate POD \cite{HiLeRa18}. Those optimal local spaces both allow for a rigorous {\it{a priori}} theory and yield a rapidly (and often exponentially) convergent approximation; in certain cases the superalgebraic convergence can be proved \cite{BabLip11}. Recently, in \cite{TadPat18,Tad16} the results in \cite{BabLip11,BaHuLi14,SmePat16} have been generalized from linear differential operators whose associated bilinear form is coercive to elliptic, inf-sup stable ones. In \cite{BuhSme18} it has been shown that those optimal local approximation spaces can be efficiently approximated by transferring methods from randomized numerical linear algebra \cite{HaMaTr11}; the local approximation approximation spaces in \cite{BuhSme18} are constructed from local solutions of the PDE with random boundary conditions. Local reduced spaces generated from random snapshots have also been suggested in \cite{Betal17,EftPat13a} and methods from randomized linear algebra have been exploited in the FETI-2$\lambda$
domain decomposition method in \cite{WanVou15} and in \cite{CEGL16} for the generalized multiscale finite element method. 

We will first present the optimal local approximation spaces as introduced in \cite{BabLip11,SmePat16} for a fixed parameter $\mathcal{P}=\{\paramFixed\}$, subsequently discuss their approximation via random sampling, and conclude this subsection with the discussion of the general case $\mathcal{P}\neq \{\paramFixed\}$. To simplify notation we will omit $\paramFixed$ as long as it is fixed.  

\subsubsection{Optimal local approximation spaces for $\mathcal{P}=\{\paramFixed\}$}\label{subsubsec:bas gen optimal}

To enable maximum flexibility regarding the shape of $\Omega$ on the user's side, we assume that we do not have any a priori knowledge of the shape of $\Omega$ when constructing the reduced order model. We thus know that the global solution $u$ satisfies the considered PDE locally on $\Omega_{out}$ but suppose that the trace of $u$ on $\partial\Omega_{out}$ is \emph{unknown} to us. Therefore, we aim at approximating all local solutions $u_{loc}$ of 
\begin{equation}\label{eq:local_PDE}
a_{loc}(u_{loc}, v) = f_{loc}(v) \quad \forall v \in V_{loc},
\end{equation}
with \emph{arbitrary Dirichlet boundary conditions on $\Gamma_{out}$}. Here, the Hilbert space $V_{loc}$ is defined such that $[H^{1}_{0}(\Omega_{out})]^{z} \subset V_{loc} \subset [H^{1}(\Omega_{out})]^{z}$, $z =1,\hdots,d$, respecting the boundary conditions on $\partial\Omega$, and $a_{loc}: [H^{1}(\Omega_{out})]^{z} \times [H^{1}(\Omega_{out})]^{z} \rightarrow \mathbb{R}$, $f_{loc}: V_{loc} \rightarrow \mathbb{R}$ are local bilinear and linear forms. We will first restrict ourselves to the case $f_{loc}=0$, $g_{D}=0$, and $\partial\Omega_{m}\cap\Gamma_{D}=\emptyset$; the general case will be dealt with at the end of this subsubsection. We may then define the space of all local solutions of the PDE as 
\begin{equation}\label{eq:space of local sol PDE}
\mathcal{H}:= \{ w \in [H^{1}(\Omega_{out})]^{z}\, : \, w \text{ solves } \eqref{eq:local_PDE}, w = 0 \text{ on }\Gamma_{D}\cap\partial\Omega_{out}\}, \enspace z=1,\hdots,d.
\end{equation}
As suggested in \cite{BabLip11,SmePat16} we introduce a transfer operator $\mathcal{T}: \mathcal{S} \rightarrow \mathcal{R}$ for Hilbert spaces $\mathcal{S}$ and $\mathcal{R}$, where $\mathcal{S}= \{ w|_{\Gamma_{out}} \, : \, w \in \mathcal{H}\}$.  We define $\mathcal{T}$ for interfaces or subdomains, respectively, for $w \in \mathcal{H}$ as
\begin{equation}\label{eq:transfer_op_gamma}
\mathcal{T}(w|_{\Gamma_{out}})=\left(w-P_{\Omega_{out}}(w)\right)|_{\Gamma_{m,m'}}
\enspace \text{or} \enspace
\mathcal{T}(w|_{\Gamma_{out}})=\left(w - P_{\Omega_{m}}(w)\right)|_{\Omega_{m}}
\end{equation}
and set
$\mathcal{R}= \{ v|_{\Gamma_{m,m'}} \, : \, v = w - P_{\Omega_{out}}(w),  w \in \mathcal{H}\}$
or
$\mathcal{R}= \{ \left(w - P_{\Omega_{m}}w\right)|_{\Omega_{m}} \, : \, w \in \mathcal{H}\}$. Here, $P_{D}$, $D \subset \Omega_{out}$, denotes an orthogonal projection onto the kernel of the bilinear form; for further details see \cite{BuhSme18,SmePat16}. In the case of heat conduction we would for instance subtract the mean value of the respective function on $D$. Note that subtracting this projection is necessary to prove compactness of the transfer operator $\mathcal{T}$. The key argument to show compactness of $\mathcal{T}$ is Caccioppoli's inequality, which estimates the energy norm of a function in $\mathcal{H}$ on $\Omega_{m}$ in terms of the $L^{2}$-norm on $\Omega_{out}$ of the respective function. Using the Hilbert-Schmidt theorem and Theorem 2.2 in  \cite[Chapter 4]{Pinkus85} it can then be shown that certain eigenfunctions of $\mathcal{T}^{*}\mathcal{T}$ span the optimal local approximation space, where $\mathcal{T}^{*}:\mathcal{R} \rightarrow \mathcal{S}$ denotes the adjoint operator of $\mathcal{T}$. As we aim at approximating $\mathcal{H}$ and thus a whole set of functions, the concept of optimality of Kolmogorov \cite{Kol36} is used: A subspace $\mathcal{R}_{n} \subset \mathcal{R}$ of dimension at most $n$ for which holds
\begin{equation*}
d_{n}(\mathcal{T}(\mathcal{S});\mathcal{R}) = \sup_{\psi \in \mathcal{S}} \inf_{\zeta \in \mathcal{R}_{n}} \frac{\|\mathcal{T}\psi - \zeta\|_{\mathcal{R}}}{\|\psi\|_{\mathcal{S}}}
\end{equation*}
is called an optimal subspace for $d_{n}(\mathcal{\mathcal{T}}(\mathcal{S});\mathcal{R})$, where the Kolmogorov $n$-width $d_{n}(\mathcal{T}(\mathcal{S});\mathcal{R})$ is defined as 
\begin{equation*}
d_{n}(\mathcal{T}(\mathcal{S});\mathcal{R}) := \underset{\dim(\mathcal{R}_{n})=n}{\inf_{\mathcal{R}_{n}\subset \mathcal{R}}} \sup_{\psi \in \mathcal{S}} \inf_{\zeta \in \mathcal{R}_{n}} \frac{\|\mathcal{T}\psi - \zeta\|_{\mathcal{R}}}{\|\psi\|_{\mathcal{S}}}.
\end{equation*}

We summarize the findings about the optimal local approximation spaces in the following theorem.

\begin{theorem}[Optimal local approximation spaces \cite{BabLip11,SmePat16}]\label{theorem:eigenvalue problem}
The optimal approximation space for $d_{n}(\mathcal{T}(\mathcal{S});\mathcal{R})$ is given by 
\begin{equation}\label{eq:optimal space}
\mathcal{R}_{n}:= \spanlin \{\chi_{1}^{sp},...,\chi_{n}^{sp}\}, \qquad \text{where} \enspace \chi_{j}^{sp}=\mathcal{T}\phi_{j},\quad j=1,...,n,
\end{equation}
and 
$\lambda_j$ are the largest $n$ eigenvalues and $\phi_j$ the corresponding eigenfunctions
that satisfy the transfer eigenvalue problem: Find $(\phi_{j},\lambda_{j}) \in (\mathcal{S},\mathbb{R}^{+})$ such that
\begin{align}\label{eq:transfer eigenvalue problem}
(\,\mathcal{T}\phi_{j}\,,\,\mathcal{T}w\,)_{\mathcal{R}} = \lambda_{j} (\,\phi_{j}\,,\,w\,)_{\mathcal{S}} \quad \forall w \in \mathcal{S}.
\end{align}
Moreover, we have: 
\begin{equation}\label{eq:n-width_equals_eigenvalue}
d_{n}(\mathcal{T}(\mathcal{S});\mathcal{R}) = \sup_{\xi \in \mathcal{S}} \inf_{\zeta \in \mathcal{R}_{n}} \frac{\| \mathcal{T}\xi - \zeta\|_{\mathcal{R}}}{\|\xi\|_{\mathcal{S}}} = \sqrt{\lambda_{n+1}}.
\end{equation}
\end{theorem} 
\begin{remark}
We emphasize that the optimal space $\mathcal{R}_{n}$ is optimal in the sense of Kolmogorov for the approximation of the range of $\mathcal{T}$ and \emph{not} necessarily for the approximation of $u(\param)$. Moreover, we remark that $\chi_{i}^{sp}$ are the left singular vectors and $\sqrt{\lambda_{i}}$ the singular values of $\mathcal{T}$.
\end{remark}

Next, for $f_{loc}\neq 0$ but still $g_{D}=0$ we solve the problem: Find $u_{loc}^{f} \in V_{loc}$ such that 
$a_{loc}(u_{loc}^{f}, v) = f_{loc}(v)$ for all $v \in V_{loc}$ and augment the space $\mathcal{R}_{n}$ either with $u_{loc}^{f}|_{\Omega_{m}}$ or $u_{loc}^{f}|_{\Gamma_{m,m'}}$. To take non-homogeneous Dirichlet boundary conditions into account 
one can proceed for instance with a standard lifting approach, adjusting $f_{loc}$ accordingly. Note that for homogeneous boundary conditions we proceed very similar to above, prescribing ``arbitrary'' boundary conditions on $\Gamma_{out}$ and homogeneous boundary conditions on $\partial\Omega \cap \partial \Omega_{out}$. The optimal local approximation space for subdomains are then defined as 
\begin{equation}\label{eq:optimal space subdomain}
\mathcal{R}_{n}^{+}:=\spanlin \{ \chi_{1}^{sp},...,\chi_{n}^{sp}, u^{f}_{loc}|_{\Omega_{m}}\} \oplus \ker(a_{m}(\cdot,v))
\end{equation}
and similarly for interfaces as
\begin{equation}\label{eq:optimal space interface}
\mathcal{R}_{n}^{+}:=\spanlin \{ \chi_{1}^{sp},...,\chi_{n}^{sp}, u^{f}_{loc}|_{\Gamma_{m,m'}}\} \oplus \ker(a_{m}(\cdot,v))|_{\Gamma_{m,m'}},
\end{equation}
respectively. Here, $\ker(a_{m}(\cdot,v))$ denotes the kernel of the mapping $a_{m}(\cdot,v): [H^{1}(\Omega_{m})]^{z} \rightarrow \mathbb{R}$, $z=1,\hdots,d$, $v \in V_{0}^{m}$ for the bilinear form $a_{m}$ defined in subsection \ref{subsec:conforming}. In case $\partial\Omega_{m}\cap\Gamma_{D}\neq\emptyset$ all modifications in this subsubsection involving the kernel of the bilinear form are waived. 

The result in \eqref{eq:n-width_equals_eigenvalue} can be exploited to derive an {\it a priori} error bound for the approximation error between the solution $u(\paramFixed)$ of \eqref{eq:weak_solution} still for a fixed reference parameter $\paramFixed$ and the optimal static condensation approximation $u^{n}(\paramFixed)$ as stated in the following proposition:

\begin{proposition}[{\it A priori} error bound \cite{SmePat16}]\label{prop:a priori bound}
Assume that the interfaces $\gamma \in \mathcal{T}_{H}^{\gamma}$ are mutually disjoint, that all interfaces have the same geometry, and that each $\Omega_{m}$, $m=1,\hdots,M$ has exactly two interfaces. Let $u(\paramFixed)$ be the (exact) solution of \eqref{eq:weak_solution} for a fixed parameter $\paramFixed$. Moreover, let $u_{n_{+}}(\paramFixed)$ be the static condensation approximation defined in subsection \ref{subsec:conforming}, where we employ the optimal interface space $\mathcal{R}_{n}^{+}$ for each $\gamma \in \mathcal{T}_{H}^{\gamma}$ and assume that the error due to the intra-element RB approximation is zero. Then, we have the following {\it a priori} error bound:
\begin{equation}\label{eq a priori bound}
\frac{\energynorm{u(\paramFixed) - u_{n_{+}}(\paramFixed)}{\paramFixed}}{\energynorm{u(\paramFixed)}{\paramFixed}} \leq \#\gamma \max_{\gamma \in \mathcal{T}_{H}^{\gamma}} \left( C_{\gamma} \,\sqrt{\lambda_{n+1}^{\gamma}}\right),
\end{equation}
where $\#\gamma$ denotes the number of interfaces in $\mathcal{T}_{H}^{\gamma}$ and $\lambda_{n+1}^{\gamma}$ is the $n+1$-th eigenvalue of \eqref{eq:transfer eigenvalue problem} for the interface $\gamma \in \mathcal{T}_{H}^{\gamma}$. The constant $C_{\gamma}$ depends only on the subdomains that share the interface $\gamma$ and neither on $\Omega$ nor on $u(\paramFixed)$.
\end{proposition}

To define reduced interface spaces $\Lambda_{N^{\gamma}}^{\gamma}$, $\gamma \in \mathcal{T}_{H}^{\gamma}$ and reduced spaces $V_{N}^{m}$, $m=1,\hdots,M$ we approximate \eqref{eq:transfer eigenvalue problem} with finite elements. To that end, we introduce a conforming FE space $V_{h;loc} \subset V_{loc}$, the FE source space $S:=\{ v|_{\Gamma_{out}}\, : \, v \in V_{h}\}$ of dimension $N_{S}$, and the FE range space $R:=\{ (v - P_{\Omega_{out}}(v))|_{\Gamma_{m,m'}}\, : \, v \in V_{h}\}$ or $R:=\{ (v - P_{\Omega_{m}})|_{\Omega_{m}}\, : \, v \in V_{h}\}$ with $\dim(R)=N_{R}$. We may then define the discrete transfer operator $T:S \rightarrow R$ for $w \in \mathcal{H}_{h}=\{ w \in V_{h}|_{\Omega_{out}}\, : \, a_{loc}(w,\varphi) = 0 \,\forall \varphi \in V_{h;loc}, \enspace w = 0 \enspace \text{on} \enspace \Gamma_{D}\cap\partial\Omega_{out}\}$ as 
\begin{equation}\label{eq: discrete transfer operator}
T(w|_{\Gamma_{out}})=\left(w-P_{\Omega_{out}}(w)\right)|_{\Gamma_{m,m'}} \enspace \text{or} \enspace
T(w|_{\Gamma_{out}})=\left(w - P_{\Omega_{m}}(w)\right)|_{\Omega_{m}}.
\end{equation}
In order to define a matrix form of the transfer operator we introduce DOF mappings $\mathbb{B}_{S\rightarrow V_{h}|_{\Omega_{out}}} \in \mathbb{R}^{\dim(V_{h}|_{\Omega_{out}}) \times N_{S}}$ and $\mathbb{B}_{V_{h}|_{\Omega_{out}}\rightarrow R} \in \mathbb{R}^{N_{R}\times \dim(V_{h}|_{\Omega_{out}})}$ that map the DOFs of $S$ to the DOFs of $V_{h}|_{\Omega_{out}}$ and the DOFs of $V_{h}|_{\Omega_{out}}$ to the DOFs of $R$, respectively. Moreover, we introduce the stiffness matrix $\mathbb{A}_{loc}$ obtained from the FE discretization of \eqref{eq:local_PDE}, where we assume that in the rows associated with the Dirichlet DOFs the non-diagonal entries are zero and the diagonal entries equal one. By denoting by $\boldsymbol{\zeta}$ the FE coefficients of $\zeta \in S$ and by defining $\mathbb{P}_{D}$ as the matrix of the orthogonal projection on the kernel of the bilinear form on $D \subset \Omega_{out}$,
we obtain the following matrix representation $\mathbb{T} \in \mathbb{R}^{N_{R} \times N_{S}}$ of the transfer operator for subdomains
\begin{eqnarray}\label{eq:matrix_form_transfer_operator}
\mathbb{T}\,\boldsymbol{\zeta} = 
\left(1 - \mathbb{P}_{\Omega_{m}} \right)
\mathbb{B}_{V_{h}|_{\Omega_{out}}\rightarrow R} \, \mathbb{A}^{-1} \mathbb{B}_{S\rightarrow V_{h}|_{\Omega_{out}}} \, \boldsymbol{\zeta}
\end{eqnarray}
and interfaces
\begin{eqnarray}\label{eq:matrix_form_transfer_operator2}
\mathbb{T}\,\boldsymbol{\zeta} = 
\mathbb{B}_{V_{h}|_{\Omega_{out}}\rightarrow R} \, 
\left(1 -  \mathbb{P}_{\Omega_{out}} \right) \,
\mathbb{A}^{-1} \mathbb{B}_{S\rightarrow V_{h}|_{\Omega_{out}}} \, \boldsymbol{\zeta}.
\end{eqnarray}
Finally, we denote by $\mathbb{M}_S$ the inner product matrix of $S$ and by $\mathbb{M}_R$ the inner product matrix of $R$. Then, the FE approximation of the transfer eigenvalue problem reads as follows: Find the eigenvectors $\boldsymbol{\zeta}_{j} \in \mathbb{R}^{N_{S}}$ and the eigenvalues $\lambda_{j} \in \mathbb{R}^{+}_0$ such that
\begin{equation}\label{eq:matrix_version_transfer_eigenvalue_problem}
\mathbb{T}^{t}\mathbb{M}_{R}\mathbb{T} \,\boldsymbol{\zeta}_{j} = \lambda_{j}\, \mathbb{M}_{S} \,\boldsymbol{\zeta}_{j}.
\end{equation}
The coefficients of the FE approximation of the basis functions $\{\chi_{h,1}^{sp},...,\chi_{h,n}^{sp}\}$ of the discrete optimal local approximation space 
\begin{equation}\label{eq:discrete optimal space}
R_{n}:=\spanlin \{\chi_{h,1}^{sp},...,\chi_{h,n}^{sp}\}
\end{equation}
are then given by
$
\boldsymbol{\chi}_{h,j}^{sp} = \mathbb{T}\,\boldsymbol{\zeta}_{j},$ $j =1,\hdots,n.
$
Adding the representation of the right-hand side, the boundary conditions, and a basis of the kernel of the bilinear form yields the optimal spaces $\Lambda_{N^{\gamma}}^{\gamma}$ and $V_{N}^{m}$. 

Note that in actual practice we would not assemble the matrix $\mathbb{T}$. Instead one may solve the PDE locally $N_{S}$ times prescribing the basis functions of $S$ as Dirichlet boundary conditions on $\Gamma_{out}$ and subsequently assemble and solve the transfer eigenvalue problem. Alternatively, one may pass $\mathbb{T}$ implicitly to the Lanczos method. For instance, the implicitly restarted Lanczos method as implemented in ARPACK \cite{LeSoYa98} requires $\mathcal{O}(n)$ local solutions of the PDE in each iteration and applications of the adjoint $T^{*}$. In the next subsubsection we will show how methods from randomized linear algebra \cite{HaMaTr11,DriMah16,Mah11,MahDri09} can be used to compute an approximation of the optimal local approximation spaces. However, beforehand, we conclude this subsubsection with some numerical experiments on the transfer eigenvalues and thus via 
Proposition \ref{prop:a priori bound} on the convergence behavior of the relative approximation error.

\begin{figure}[t]
\begin{center}
\includegraphics[scale = 0.2]{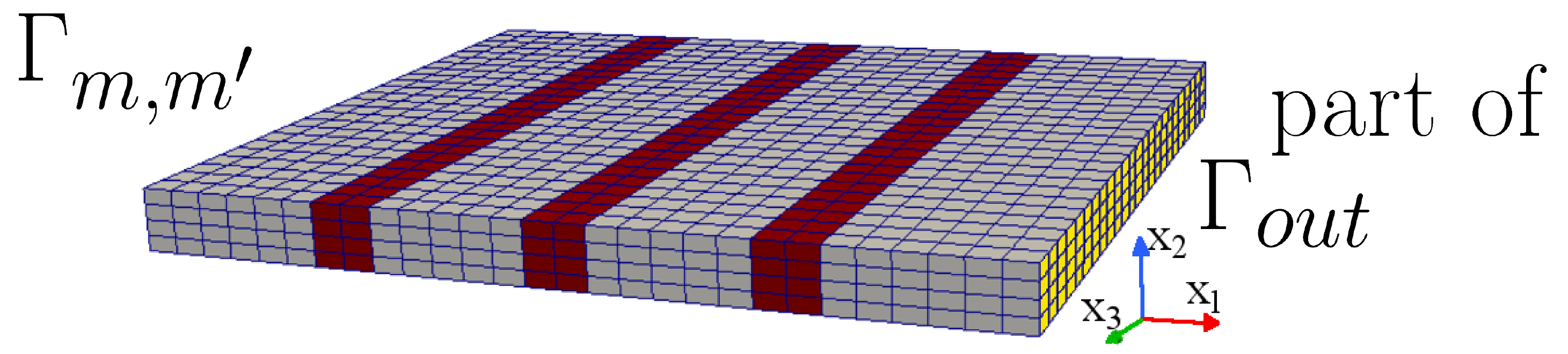}
\caption{Mesh in the subdomain $\Omega_{2}$ for the ship stiffener. The part of $\Gamma_{out}$ in $\Omega_{2}$ is indicated in yellow and on the opposite we have the interface $\Gamma_{m,m'}$. In the red shaded areas Young's modulus may be varied between $1$ and $20$ and in the gray areas we consider $E(\param)\equiv 1$.} \label{fig boat mesh}
\end{center}
\end{figure}

To this end, we present the simplified model for a ship stiffener from \cite{SmePat16}: We consider $\overline{\Omega}_{out} = \bar{\Omega}_{1} \cup \bar{\Omega}_{2}$ and $\Gamma_{m,m'}=\Gamma_{1,2}=\bar{\Omega}_{1} \cap \bar{\Omega}_{2}$, where $\Omega_{2}$ is depicted in Fig.~\ref{fig boat mesh}, $\Omega_{1}$ is just a shifted version of $\Omega_{2}$, and the part of $\Gamma_{out}$ in $\Omega_{2}$ is indicated in yellow in Fig.~\ref{fig boat mesh}. We allow $E(\param)$ to vary in the red areas of the subdomains between $1$ and $20$ and prescribe $E(\param)\equiv 1$ in the gray areas; we choose $\boldsymbol{G}(\param)=(0,0,0)^{T}$.

In detail, we consider $\Omega_{1}=(-0.7,0.7)\times (-0.05,0.05)\times (-0.6,0.6)$, $\Omega_{2}=(0.7,2.1)\times (-0.05,0.05)\times (-0.6,0.6)$ and $\Gamma_{out}=\{-0.7\}\times(-0.05,0.05)\times (-0.6,0.6) \cup \{2.1\}\times(-0.05,0.05)\times (-0.6,0.6)$. We employ a conforming linear FE space associated with the mesh depicted in Fig.~\ref{fig boat mesh}, resulting in $N = 13125$ degrees of freedom per subdomain and an FE interface space of dimension $N_{\Gamma}=375$. Finally, we equip both $S$ and $R$ with a lifting inner product based on the lifting operator $\mathcal{E}_{\Gamma \rightarrow \Omega_{m}}(\paramFixed)$ defined in subsubsection \ref{subsub:steklovpoincare}; for further details we refer to \cite{SmePat16}.

We consider different values for Young's modulus (ratios) $E_{i}^{r}$, $i=1,2$ in the red areas of the subdomains and observe in Fig.~\ref{fig eigenvalues boat} for the ship stiffener application an exponential convergence of order $\approx e^{-n}$ of the eigenvalues $\lambda_{n}(\param)$ and thus the static condensation approximation. 
We emphasize that we observe in Fig.~\ref{fig eigenvalues boat} that the eigenvalues associated with the stiffened plate ($E_{1}^{r}=E_{2}^{r}=20$) decay fastest, while we see the slowest decay for the non-stiffened plate ($E_{1}^{r}=E_{2}^{r}=1$). This is consistent with the expectation that stiffening the plate decreases the deflection of the plate, eliminating the higher eigenmodes. Moreover, an inspection of the optimal interface modes reveals many ``classical'' mode shapes such as bending or torsional modes of beams and demonstrates again the physical significance of the optimal modes. Also for beams of different shapes, including an I-beam with a crack and thus an irregular domain, an exponential convergence of the transfer eigenvalues and the physical significance of the transfer eigenmodes can be observed; for further details see \cite{SmePat16}.
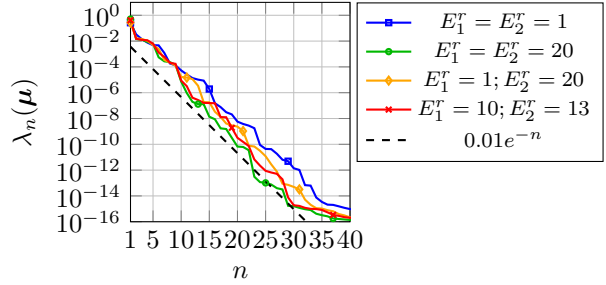
\begin{figure}[t]
\centering
\begin{tikzpicture}
\begin{semilogyaxis}[
    width=4.5cm,
    height=4.5cm,
    xmin=1,
    xmax=40,
    ymin=1e-16,
    ymax=10,
    xlabel=$n$,
    ylabel=$\lambda_n(\param)$,
    grid=both,
    grid style={line width=.1pt, draw=gray!20},
    major grid style={line width=.2pt,draw=gray!50},
    xtick={1,5,10,15,20,25,30,35,40},
    ytick={1e-16, 1e-14, 1e-12, 1e-10, 1e-8, 1e-6, 1e-4, 1e-2, 1e0},
    legend pos=outer north east,
  ]
  \addplot+[solid, blue, thick, mark=square, mark size=1pt, mark repeat=14] table[x index=0,y index=1] {eigenvalues.dat};
    \addplot+[solid, gruen, thick, mark=none, mark = o, mark size = 1.0pt, mark repeat = 12] table[x index=0,y index=2] {eigenvalues.dat};
  \addplot+[solid, gelb, thick, mark=none,mark=diamond, mark size =1.5pt, mark repeat=10] table[x index=0,y index=3] {eigenvalues.dat};
  \addplot+[solid, red, thick, mark=x, mark size=1.5pt, mark repeat =18] table[x index=0,y index=4] {eigenvalues.dat};
  \addplot+[dashed, black, thick, mark=none] table[x index=0,y index=5] {eigenvalues.dat};
  \legend{
  {\footnotesize $E^{r}_1=E^{r}_2=1$},
  {\footnotesize $E^{r}_1=E^{r}_2=20$}, 
  {\footnotesize $E^{r}_1=1; E^{r}_2=20$}, 
  {\footnotesize $E^{r}_1=10; E^{r}_2=13$},  
  {\footnotesize $0.01 e^{-n}$}
  };
\end{semilogyaxis}
\end{tikzpicture}
\caption{Eigenvalues $\lambda_{n}(\param)$ for different Young's modulus ratios $E_{i}^{r}$ in $\Omega_{i}$, $i=1,2$.\label{fig eigenvalues boat}}
\end{figure}

\subsubsection{Randomized Training}\label{subsubsec:bas gen random}

\begin{algorithm2e}[t]
\DontPrintSemicolon
\SetAlgoVlined
\SetKwFunction{AdaptiveRandomizedRangeApproximation}{AdaptiveRandomizedRangeApproximation}
\SetKwInOut{Input}{Input}
\SetKwInOut{Output}{Output}
  \Input{Operator $T$, target accuracy $\normalfont{\texttt{tol}}$, number of test vectors $n_t$, maximum failure probability $\varepsilon_\mathrm{algofail}$}
  \Output{space $R_{n}^\mathrm{rand}$ with property $P\left(\norm{T - P_{R_{n}^\mathrm{rand}}T} \leq \normalfont{\texttt{tol}} \right) > \left(1 - \varepsilon_\mathrm{algofail}\right)$}
  \textbf{Initialize:} $B \leftarrow \emptyset$, $M \leftarrow \{T D_S^{-1} \boldsymbol{r}_1, \ \dots, \  T D_S^{-1} \boldsymbol{r}_{n_t} \} $ \label{algo:line:init}\;
  \textbf{Compute error estimator factors:}\\ $\varepsilon_\mathrm{testfail} \leftarrow \varepsilon_\mathrm{algofail} / N_T$; \quad $c_\mathrm{est} \leftarrow \left[ \sqrt{2 \lambda^{\mathbb{M}_S}_{min}} \ \mathrm{erf}^{-1} \left( \sqrt[n_t]{\varepsilon_{\mathrm{testfail}} } \right)
\right]^{-1}$ \label{algo:line:constant_init} \;
  \While{$\left( \max_{t \in M} \norm{t}_R \right) \cdot c_\mathrm{est} > \normalfont{\texttt{tol}}$ \label{algo:line:convergence}}{
    $B \leftarrow B \cup (T D_S^{-1} \boldsymbol{r}) $ \label{algo:line:basis_extension}\;
    $B \leftarrow \mathrm{orthonormalize}(B)$ \label{algo:line:basis_extension2}\;
    \textbf{orthogonalize test vectors:}
    $M \leftarrow \left\{t - P_{\mathrm{span}\{B\}} t \ \Big| \ t \in M  \right\}$ \label{algo:line:test_vector_update}\;
  }
  \Return $R_{n}^\mathrm{rand} = \spanlin\{B\}$\;
\caption{Adaptive Randomized Range Approximation}
\label{algo:adaptive_range_approximation}
\end{algorithm2e}
In order to compute an efficient approximation $R_{n}^\mathrm{rand}$ of $R_{n}$ the adaptive randomized range approximation algorithm \ref{algo:adaptive_range_approximation} as suggested in \cite{BuhSme18} iteratively enhances the reduced space with applications of $T$ to a random function until a certain convergence criterion is satisfied. 

In detail, in each iteration in line \ref{algo:line:basis_extension} we draw a new random vector $\boldsymbol{r} \in \mathbb{R}^{N_{S}}$ whose entries are independent and identically distributed random variables with standard normal distribution. Then, we employ the mapping $D_{S}^{-1}:\mathbb{R}^{N_{S}} \rightarrow S$ to define a unique FE function in $S$ whose coefficients are the components of $\boldsymbol{r}$. Subsequently, we apply the transfer operator $T$ to $D_{S}^{-1}\boldsymbol{r}$, meaning that we solve the PDE locally on $\Omega_{out}$ with random boundary conditions and restrict the solution to $\Omega_{m}$ or $\Gamma_{m,m'}$; the resulting function is added to the set of basis functions $B$. Finally, the basis $B$ is orthonormalized. Note that the orthonormalization is numerically challenging, 
as the basis functions are nearly linear dependent when
$\mathrm{span}\{B\}$ is already a good approximation of the range of $T$; in \cite{BuhSme18} using the numerically stable Gram-Schmidt with adaptive re-iteration from \cite{BEOR14a} is suggested. The main loop of the algorithm is terminated when the following a posteriori norm estimator is smaller than the desired tolerance $\normalfont{\texttt{tol}}$.
\begin{proposition}[A probabilistic a posteriori norm estimator \cite{BuhSme18}]\label{thm:errest_reliable}
Let $\boldsymbol{r}_i$, $i=1,\hdots,n_{t}$ be  $n_{t}$ random normal test vectors and $\lambda^{\mathbb{M}_S}_{min}$ and  $\lambda^{\mathbb{M}_S}_{max}$ the smallest and largest eigenvalues of the matrix of the inner product in $S$. Then, the a posteriori norm estimator 
\begin{equation}\label{eq:rand a posteriori error estimator}
\Delta(n_t, \varepsilon_\mathrm{testfail}) 
:= \cest \max_{i \in 1, \dots, n_t} \norm{ (T - P_{R_n^\mathrm{rand}}T) \ D_S^{-1} \ \boldsymbol{r}_i }_R
\end{equation}
satisfies
\begin{equation}\label{eq:rand err_est_upper bound}
P\left\{\|T - P_{R_n^\mathrm{rand}}T\| \leq  \Delta(n_t, \varepsilon_\mathrm{testfail})\right\} \geq (1 - \varepsilon_\mathrm{testfail}),
\end{equation}
where $\cest := 
1 / [ (2 \lambda^{\mathbb{M}_S}_{min})^{1/2} \ \mathrm{erf}^{-1} ( \sqrt[n_t]{\varepsilon_{\mathrm{testfail}} } ) ].
$ Additionally, there holds
$$
P\left\{ \frac{
\Delta(n_t, \varepsilon_\mathrm{testfail}) 
}{
\norm{T - P_{R_n^\mathrm{rand}}T}
} \leq \ceff \right\} \geq 1 - \varepsilon_\mathrm{testfail}
,
$$
where the constant $\ceff$ is defined as
$$
\ceff := 
\left[
Q^{-1}\left(\frac{N_T}{2}, \frac{\varepsilon_\mathrm{testfail}}{n_t}\right)
\frac{\lambda^{\mathbb{M}_S}_{max}}{\lambda^{\mathbb{M}_S}_{min}}
\left( \mathrm{erf}^{-1} \left( \sqrt[n_t]{\varepsilon_\mathrm{testfail}} \right) \right)^{-2}
\right]^{1/2}
$$
and $Q^{-1}$ is the inverse of the upper normalized incomplete gamma function.
\end{proposition}

The constant $\cest$ is calculated in line \ref{algo:line:constant_init} using $N_T$, which denotes the rank of operator $T$. In practice $N_{T}$ is unknown
and an upper bound for $N_T$ such as $\min(N_S, N_R)$ can be used instead. Note that the term
$
\left( \max_{t \in M} \norm{t}_R \right) \cdot \cest
$
is the norm estimator \eqref{eq:rand a posteriori error estimator}.
The test vectors are reused for all iterations.

To finally analyze the failure probability of Algorithm \ref{algo:adaptive_range_approximation} we first note that after $N_{T}$ steps we have $R_{n}^\mathrm{rand}=\range(T)$ and thus $\norm{T- P_{R_{n}^\mathrm{rand}}T}=0$, yielding the termination of Algorithm \ref{algo:adaptive_range_approximation}. Using the fact that the a posteriori error estimator defined in \eqref{eq:rand a posteriori error estimator} is therefore executed at most $N_T$ times combined with the probability for one estimate to fail in \eqref{eq:rand err_est_upper bound} and an union bound argument we infer that the failure probability for the whole algorithm is
$
\varepsilon_\mathrm{algofail} 
\leq
N_T
 \ \varepsilon_\mathrm{testfail}
.
$

Remarkably, the convergence behavior of the reduced space $R_{n}^\mathrm{rand}$ is only slightly worse than the rate $\sqrt{\lambda_{n+1}}$, which is achieved by the optimal local approximation spaces defined in Theorem \ref{theorem:eigenvalue problem}:
\begin{proposition}[A priori error bound \cite{BuhSme18}]\label{prop:a priori}
\label{thm:convergence_rate_main}
Let $\lambda^{\mathbb{M}_R}_{max}$ and $\lambda^{\mathbb{M}_R}_{min}$ denote the largest and smallest eigenvalues of the inner product matrix $\mathbb{M}_R$ and let $R_{n}^\mathrm{rand}$ be the outcome of Algorithm \ref{algo:adaptive_range_approximation}. Then, for $n\geq 4$ there holds
\begin{equation}\label{eq:a priori mean}
\mathbb{E}
\norm{T - P_{R_{n}^\mathrm{rand}} T}
\leq C_{R,S}\vspace{-5pt}
\min_{\overset{k+p=n}{k\geq 2, p\geq 2}}
\left[
\left(1+\sqrt{\frac{k}{p-1}} \right) \vspace{-1pt} \sqrt{\lambda_{k+1}} +
\frac{ e \sqrt{n}}{p} \left( \sum_{j > k} \lambda_j \right)^{\frac{1}{2}}
\right]
,
\end{equation}
where $C_{R,S}=(\lambda^{\mathbb{M}_{R}}_{max}/\lambda^{\mathbb{M}_{R}}_{min})^{1/2}(\lambda^{\mathbb{M}_{S}}_{max}\lambda^{\mathbb{M}_{S}}_{min})^{1/2}$.
\end{proposition}

It can be observed in numerical experiments that the a priori bound in Proposition \ref{prop:a priori} is sharp in terms of the predicted convergence behavior as we will show now for a test case from \cite{BuhSme18}. Moreover, we will investigate the performance of Algorithm \ref{algo:adaptive_range_approximation} also for a test case from \cite{BuhSme18}. 
To that end, let $\widehat{\Omega}_{m}=(-0.5,0.5)\times (-0.25,0.25) \times (-0.5,0.5)$ and $\Omega_{m} = (-0.5,0.5) \times (-0.5,0.5) \times (-0.5,0.5)$ be the subdomains on which we aim to construct a local approximation space, $\widehat{\Omega}_{out}=(-2,2)\times (-0.25,0.25) \times (-2,2)$  and $\Omega_{out} = (-2,2) \times (-0.5,0.5) \times (-2,2)$ the corresponding oversampling domains and $\widehat{\Gamma}_{out}=\{-2,2\} \times (-0.25,0.25) \times (-2,2) \cup (-2,2) \times (-0.25,0.25) \times \{-2,2\}$ and $\Gamma_{out}=\{-2,2\} \times (-0.5,0.5) \times (-2,2) \cup (-2,2) \times (-0.5,0.5) \times \{-2,2\}$ the respective outer boundaries. On $\partial\widehat{\Omega}_{out} \setminus \widehat{\Gamma}_{out}$ and $\partial\Omega_{out} \setminus \Gamma_{out}$ we prescribe homogeneous Neumann boundary conditions and we suppose that $\widehat{\Omega}_{out}$ and  $\Omega_{out}$ do not border the Dirichlet boundary of $\Omega$. For the FE discretization we use a regular mesh with hexahedral elements and a mesh size $h=0.1$ in each space direction and a corresponding conforming FE space with linear FE resulting in $\dim(V_{h}|_{\widehat{\Omega}_{out}})=30258$, $\dim(R)=N_{R}=2172$, $\dim(S)=N_{S}=2880$ for $\widehat{\Omega}_{out}$ and $V_{h}|_{\Omega_{out}}=55473$, $N_{R}=3987$, and $N_{S}=5280$ for $\Omega_{out}$.\footnote{Note that although in theory we should subtract the orthogonal projection on the six rigid body motions from the FE basis functions, in actual practice we avoid that by subtracting the orthogonal projection from the harmonic extensions only.} We equip the source space $S$ with the $L^{2}$-inner product and the range space $R$ with the energy inner product. Finally, for all results in this subsubsection we computed the statistics over $1000$ samples.

Analyzing the convergence behavior of $\mathbb{E}(\|T - P_{R_{k+p}^\mathrm{rand}}T\|)$ on $\widehat{\Omega}_{out}$ for a growing number of randomly generated basis functions $k$ and a (fixed) oversampling parameter $p=2$ in Fig.~\ref{fig: linear elasticity a priori 2} we see that until $k\approx 75$ the a priori bound reproduces the convergence behavior of $\mathbb{E}(\|T - P_{R_{k+p}^\mathrm{rand}}T\|)$ perfectly. We may thus conclude that the a priori bound in \eqref{eq:a priori mean} seems to be sharp regarding the convergence behavior of $\mathbb{E}(\|T-P_{R_{k+p}^\mathrm{rand}}T\|)$ in the basis size $k$. We also observe that the a priori bound is rather pessimistic as it overestimates $\mathbb{E}(\|T - P_{R_{k+p}^\mathrm{rand}}T\|)$ by a factor of more than $100$; this is mainly due to the square root of the conditions of the inner product matrices. 

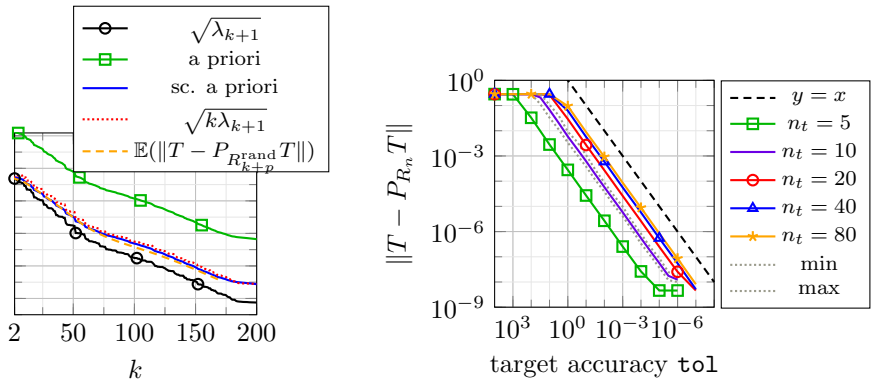
\begin{figure}
\centering
\begin{subfigure}[t]{0.4\textwidth}
\centering
\begin{tikzpicture}
\begin{semilogyaxis}[
    width=4.8cm,
    height=4.0cm,
    xmin=2,
    xmax=200,
    ymin=1e-9,
    ymax=150,
    xlabel=$k$,
    yticklabels={,,},
    grid=both,
    grid style={line width=.1pt, draw=gray!20},
    major grid style={line width=.2pt,draw=gray!50},
    minor xtick={2,25,50,75,100,125,150,175,200},
    minor ytick={1e-9, 1e-8, 1e-7, 1e-6, 1e-5, 1e-4, 1e-3, 1e-2, 1e-1, 1e0, 1e1, 1e2},
    xtick={2,50,100,150,200},
    ytick={1e2, 1e0, 1e-2, 1e-4, 1e-6, 1e-8},
    max space between ticks=25,
   legend style={at={(1.3,1.7)},anchor=north east},
  ]
 \addplot+[solid, black, thick, mark=o,mark repeat=50] table[x index=0,y index=2] {elasticity_plate_lambda.dat};
   \addplot+[solid,          gruen, thick, mark=square, mark repeat=50] table[x index=0,y index=1] {elasticity_plate_a_priori_bound.dat};
  \addplot+[solid, blue, thick, mark=none] table[x index=0,y index=2] {elasticity_plate_a_priori_bound.dat};
  \addplot+[densely dotted, red, thick, mark=none] table[x index=0,y index=1] {elasticity_plate_lambda.dat};
  \addplot+[densely dashed, gelb, thick, mark=none] table[x index=0,y index=1] {elasticity_plate_error.dat};
  \legend{
  {\footnotesize $\sqrt{\lambda_{k+1}}$}, 
  {\footnotesize a priori}, 
  {\footnotesize sc. a priori}, 
  {\footnotesize $\sqrt{k\lambda_{k+1}}$},
  {\footnotesize $\mathbb{E}(\|T-P_{R_{k+p}^\mathrm{rand}}T\|)$}
  };
\end{semilogyaxis}
\end{tikzpicture}
\subcaption{\footnotesize Convergence behavior on $\widehat{\Omega}_{out}$.}
\label{fig: linear elasticity a priori 2}
\end{subfigure}
\begin{subfigure}[t]{0.55\textwidth}
\centering
\begin{tikzpicture}
\begin{loglogaxis}[
    width=4.5cm,
    height=4.6cm,
    xmin=1e-8,
    xmax=1e4,
    x dir=reverse,
    ymin=1e-9,
    ymax=1e0,
    xlabel=target accuracy \texttt{tol},
    ylabel=$\|T-P_{R_{n}}T\|$,
    grid=both,
    grid style={line width=.1pt, draw=gray!20},
    major grid style={line width=.2pt,draw=gray!50},
    minor xtick={1e-15, 1e-14, 1e-13, 1e-12, 1e-11, 1e-10, 1e-9, 1e-8, 1e-7, 1e-6, 1e-5, 1e-4, 1e-3, 1e-2, 1e-1, 1e0, 1e1, 1e2, 1e3},
    minor ytick={1e-15, 1e-14, 1e-13, 1e-12, 1e-11, 1e-10, 1e-9, 1e-8, 1e-7, 1e-6, 1e-5, 1e-4, 1e-3, 1e-2, 1e-1, 1e0, 1e1, 1e2, 1e3},
    xtick={1e3,1e0, 1e-3, 1e-6},
    ytick={1e3, 1e0, 1e-3, 1e-6, 1e-9, 1e-12},
    legend pos=outer north east,
  ]
  \addplot+[densely dashed,black,thick,mark=none] table[x index=0,y index=0] {elasticity_test_vectors_y_x.dat};
  \addplot+[solid,gruen,thick,mark=square] table[x index=0,y index=1] {elasticity_test_vectors_r_fuenf.dat};
  \addplot+[solid,lila,thick,mark=none] table[x index=0,y index=1] {elasticity_test_vectors_r_zehn.dat};
 \addplot+[solid,red,thick,mark=o, mark repeat=5] table[x index=0,y index=1] {elasticity_test_vectors_r_rest.dat};
  \addplot+[solid,blue,thick,mark=triangle, mark repeat=3] table[x index=0,y index=2] {elasticity_test_vectors_r_rest.dat};
  \addplot+[solid,gelb,thick,mark=star, mark repeat=2] table[x index=0,y index=3] {elasticity_test_vectors_r_rest.dat};
    \addplot+[densely dotted, grau, thick, mark=none] table[x index=0,y index=1] {elasticity_adaptive_prctile.dat};
  \addplot+[densely dotted, grau, thick, mark=none] table[x index=0,y index=5] {elasticity_adaptive_prctile.dat};
  \legend{
    {\footnotesize $y=x$},
    {\footnotesize $n_t=5$},
    {\footnotesize $n_t=10$},
    {\footnotesize $n_t=20$},
    {\footnotesize $n_t=40$},
    {\footnotesize $n_t=80$},
    {\footnotesize min},
    {\footnotesize max},
  };
\end{loglogaxis}
\end{tikzpicture}
\subcaption{\footnotesize Median for varying $n_{t}$ on $\Omega_{out}$}
\label{fig: linear elasticity adaptive algorithm 2}
\end{subfigure}
\caption{Comparison of the convergence behavior of $\sqrt{\lambda_{k+1}}$, $\sqrt{k \lambda_{k+1}}$, $\mathbb{E}(\|T-P_{R^{k+p}}T\|)$, the a priori error bound \eqref{eq:a priori mean}, and the a priori error bound of \eqref{eq:a priori mean} scaled with a constant such that its value for $k=2$ equals the one of $\mathbb{E}(\|T-P_{R_{k+p}^\mathrm{rand}}T\|)$ (sc. a priori) for increasing $k$ for and $p=2$ for the oversampling domain $\widehat{\Omega}_{out}$ (a).
Median of the projection error $\|T-P_{R_{n}^\mathrm{rand}}T\|$ for a decreasing target accuracy \texttt{tol} for a varying number of test vectors $n_{t}$ and the minimal and maximal values for $n_t=10$ on $\Omega_{out}$ (b).}
\end{figure}

Regarding the performance of Algorithm \ref{algo:adaptive_range_approximation} on $\Omega_{out}$ we first observe in Fig.~\ref{fig: linear elasticity adaptive algorithm 2} that the actual error $\|T-P_{R_{n}^\mathrm{rand}}T\|$ lies below the target tolerance \texttt{tol} for all $1000$ samples for $n_t=10$; which holds also true for all other considered values of $n_t$. Here, we prescribe $\varepsilon_\mathrm{algofail} = 10^{-10}$ and use $3993$ as an upper bound for $N_{T}$. We see in Fig.~\ref{fig: linear elasticity adaptive algorithm 2} that increasing the number of test vectors $n_{t}$ from $5$ to $10$ or from $10$ to $20$ increases the ratio between the median of the actual error $\|T-P_{R_{n}^\mathrm{rand}}T\|$ and the target accuracy $\normalfont{\texttt{tol}}$ significantly --- for the former by more than one magnitude --- while an increase from $n_{t}=40$ to $n_{t}=80$ has hardly any influence; similar results have been obtained in \cite{BuhSme18} for heat conduction and a Helmholtz problem. This can be explained by the scaling of the effectivity of the employed a posteriori error estimator, which is of the order of $1000$ for $n_{t}=5$ and of the order of $10$ for $n_{t} \geq 20$. Regarding the choice of $n_t$ it seems that for the present test case a value of about $20$ is in the sweet spot. We thus infer that for the present test case only very few local solutions in addition to the optimal amount are required, demonstrating that Algorithm \ref{algo:adaptive_range_approximation} performs nearly optimally in terms of computational complexity for the current problem. 

\subsubsection{The general setting $\mathcal{P}\neq\{\paramFixed\}$}\label{subsubsec:bas gen para}

The processes in subsubsection \ref{subsubsec:bas gen optimal} and \ref{subsubsec:bas gen random} yield for every $\param \in \mathcal{P}$ the local approximation space $R_{n}^{+}(\param)$ for this specific parameter $\param \in \mathcal{P}$. $R_{n}^{+}(\param)$ can also be generated by some other process, where we require that there holds
\begin{equation}\label{eq:requirement local space}
\| T(\param) - P_{R_{n}(\param)}T(\param)\| \leq \frac{\varepsilon}{2 C_{1}(\mathcal{T}_{H},\param)} 
\end{equation}
possibly only at high probability and that $R_{n}^{+}(\param)$ is defined as the direct sum of $R_{n}(\param)$, the kernel of the bilinear form, and representations of non-homogeneous Dirichlet boundary conditions and the right-hand side. We abuse notation in this subsubsection by omitting henceforth the remark that the estimate may only hold in a probabilistic sense. The constant $C_{1}(\mathcal{T}_{H},\param)$ has to be chosen in such a manner that if one uses the parameter-dependent spaces $R_{n}^{+}(\param)$ to define $u^{n}(\param)$, we have 
\begin{equation}
\frac{\energynorm{ u(\param) - u_{n_{+}}(\param) }{\param}}{\energynorm{ u(\param)}{\param}} \leq \frac{\varepsilon}{2}.
\end{equation}

The spectral greedy algorithm as introduced in \cite{SmePat16}\footnote{For a generalization to a setting where the discrete parameter set describes different geometries such as a beam with or without a crack we refer to \cite{Sme18}.} constructs \emph{one} (quasi-optimal) parameter-independent approximation space $R_{N}$ which approximates those parameter-dependent spaces $R_{n}^{+}(\param)$ with a given accuracy on a finite dimensional training set $\Xi \subset \mathcal{P}$. In the spectral greedy algorithm we exploit the fact that we expect that the local spaces $R_{n}^{+}(\param)$, and in particular the spectral modes that correspond to the largest eigenvalues, are not affected too much by a variation in the parameter thanks to the expected very rapid decay of the higher eigenfunctions in the interior of $\Omega_{out}$. 

The spectral greedy as described in Algorithm \ref{algo:spectral greedy} then proceeds as follows. After the initialization we compute for all $\param \in \Xi$ the parameter-dependent spaces $R_{n}^{+}(\param)$ such that we have \eqref{eq:requirement local space}. Note that for a decomposition $\mathcal{T}_{H}$ with mutually disjoint interfaces (also called ports), where each $\Omega_{m}$, $m=1,\hdots,M$ has exactly two interfaces and all interfaces have the same geometry we have the following {\it a priori} error bound \cite{SmePat16} for the error between $u(\param)$ and the continuous port-reduced static condensation approximation $u_{n_{+}}(\param)$ corresponding to the \emph{parameter-dependent} optimal interface space $\mathcal{R}_{n}^{+}(\param)$:
\begin{equation}\label{eq: a priori bound system}
\frac{\energynorm{ u(\param) - u_{n_{+}}(\param) }{\param}}{\energynorm{ u(\param)}{\param}} \leq \#\gamma c_{1}(\param)c_{2}(\param) \, \max_{\gamma \in \mathcal{T}_{H}^{\gamma}} \left(  C_{\gamma,1}(\Omega_{\gamma},\param) \,\sqrt{\lambda_{\gamma,n+1}(\param)}\right).
\end{equation}
Here, the constant $C_{\gamma,1}(\Omega_{\gamma},\param)$ depends only on the subdomains that share $\gamma$ and not on $\Omega$ or on $u(\param)$. Moreover, $c_{1}(\param)$ and $c_{2}(\param)$ are chosen such that we have $c_{1}(\param) \energynorm{\cdot}{\paramFixed} \leq \energynorm{\cdot}{\param} \leq c_{2}(\param) \energynorm{\cdot}{\paramFixed}$ for all $\param \in \mathcal{P}$ and a fixed reference parameter $\bar{\param} \in \mathcal{P}$. Choosing $C_{1}(\mathcal{T}_{H},\param)=\#\gamma c_{1}(\param)c_{2}(\param) \, \max_{\gamma \in \mathcal{T}_{H}^{\gamma}} C_{\gamma,1}(\Omega_{\gamma},\param)$ and $\sqrt{\lambda_{\gamma,n+1}(\param)} \leq \varepsilon/2$ yields a reduced space $R_{n}^{+}(\param)$ that satisfies the requirements stated in the beginning for every $\param \in \Xi$. Although precise estimates for $C_{\gamma,1}(\Omega_{\gamma},\param)$ can be obtained, setting $C_{\gamma,1}(\Omega_{\gamma},\param)=1$ yields in general good results as another value would just result in rescaling $\varepsilon$; for further details see \cite{SmePat16}. After having collected all functions on $\Gamma_{m,m'}$ or $\Omega_{m}$ that are essential to obtain a good approximation for all local solutions $u_{loc}(\param)$ of the PDE evaluated on $\Gamma_{m,m'}$ or $\Omega_{m}$, $\param \in \Xi$, we must select a suitable basis from those functions. This is realized in an iterative manner in Lines 5-14. 

\begin{algorithm2e}[t]
\caption{spectral greedy \cite{SmePat16}\label{algo:spectral greedy}}
\SetAlgoVlined
\SetKwInOut{Input}{Input}\SetKwInOut{Output}{Output}
\Input{train sample $\Xi \subset \mathcal{P}$, tolerance $\varepsilon$}
\Output{set of chosen parameters $\Xi_{N}$, local approximation space $R_N$}
\BlankLine
{\bf Initialize} $N\leftarrow \dim(\mathrm{ker}(a_{m}(\cdot,v)))$, $\Xi_{N}\leftarrow\emptyset, R_{N}\leftarrow \ker(a_{m}(\cdot,v))\text{ or } R_{N}\leftarrow \ker(a_{m}(\cdot,v))|_{\Gamma_{m,m'}}$ \\
\ForEach{$\param \in \Xi$}{
Compute $R_{n}^{+}(\param)$ such that $\| T(\param) - P_{R_{n}^{+}(\param)}T(\param)\| \leq \frac{\varepsilon}{2 C_{1}(\mathcal{T_{H}},\param)}$.\label{alg compute Lambda}
}
\BlankLine
\While{true}{
\If{$\max_{\param \in \Xi} E(S(R_{n}^{+}(\param)),R_N) \leq \varepsilon/(\varepsilon + 2C_{2}(\mathcal{T}_{H},\param)c_{1}(\param)c_{2}(\param))$ \label{alg thres}}{
{\bf return}\\
}
$\param^{*}\leftarrow \arg \max_{\param \in \Xi} E(S(R_{n}^{+}(\param)),R_N)$\label{alg mu max}\\
$\Xi_{N+1} \leftarrow \Xi_{N} \cup \param^{*}$\\
$\kappa \leftarrow \arg \sup_{\rho \in S(R_{n}^{+}(\param^{*}))} \inf_{\zeta \in R_N} \|\rho - \zeta\|_{R}$\label{kappa}\\
$R_{N+1} \leftarrow R_N + \mathrm{span}\{\kappa\}$ \label{alg extend}\\
$ N \leftarrow N+1$}
 \Return $\Xi_{N}$, $R_N$
\end{algorithm2e}

In each iteration we first identify in Line \ref{alg mu max} the reduced space $R_{n}^{+}(\param^{*})$ that maximizes the deviation 
\begin{equation*}
E(S(R_{n}^{+}(\param)),R_N):= \sup_{\xi \in S(R_{n}^{+}(\param))} \inf_{\zeta \in R_{N}} \|\xi - \zeta\|_{R}, \quad \param \in \Xi,
\end{equation*}
where possible choices of $S(R_{n}^{+}(\param))\subset R_{n}^{+}(\param)$ will be discussed below. Subsequently, we determine in Line \ref{kappa} the function $\kappa \in S(R_{n}^{+}(\param^{*}))$ that is worst approximated by the space $R_{N}$ and enhance $R_{N}$ with the span of $\kappa$. The spectral greedy algorithm terminates if for all $\param \in \Xi$ we have 
\begin{equation}\label{eq:stop_greedy}
\max_{\param \in \Xi} E(S(R_{n}^{+}(\param)),R_{N}) \leq \varepsilon/(\varepsilon + 2C_{2}(\mathcal{T}_{H},\param)c_{1}(\param)c_{2}(\param))
\end{equation}
for a constant $C_{2}(\mathcal{T}_{H},\param)$, which can in general be chosen equal to one. We emphasize, that both $C_{1}(\mathcal{T}_{H},\param)$ and $C_{2}(\mathcal{T}_{H},\param)$ do in general only depend on the \emph{number} of faces or subspaces on which the respective reduced space $R_{N}$ is used and \emph{not on the precise decomposition} of $\Omega$; see \eqref{eq: a priori bound system}. A slight modification of the stopping criterion \eqref{eq:stop_greedy} and a different scaling of $\varepsilon$ in the threshold for the a priori error bound in Line \ref{alg compute Lambda} allows to prove that after termination of the spectral greedy for a decomposition $\mathcal{T}_{H}$ with mutually disjoint interfaces, where each $\Omega_{m}$, $m=1,\hdots,M$ has exactly two interfaces and all interfaces have the same geometry we have \cite{SmePat16}
\begin{equation}
\energynorm{u(\param) - u_{N}(\param)}{\param}/\energynorm{u(\param)}{\param} \leq \varepsilon.
\end{equation}
Here, $u_{N}(\param)$ is the continuous port-reduced static condensation approximation corresponding to $\mathcal{R}_{N}$; $\mathcal{R}_{N}$ being the continuous outcome of the spectral greedy. \\

{\bf Choice of the subset $S(R_{n}^{+}(\param))$} First, we emphasize that in contrast to the standard greedy as introduced in \cite{VePrRoPa03} we have an ordering of the basis functions in $R_{n}^{+}(\param)$ in terms of their approximation properties thanks to the transfer eigenvalue problem; the sorting of the basis functions in terms of their approximation properties is implicitly saved in their norms as $\|\chi_{j}^{sp}(\param)\|_{R}^{2} = \lambda_{j}(\param)$, $j=1,\hdots,n$. To obtain local approximation spaces $R_{N}$ that yield a (very) good approximation $u^{N}(\param)$ already for moderate $N$ it is therefore desirable that the spectral greedy algorithm selects the lower eigenmodes sooner rather than later during the {\tt while}-loop. As suggested in \cite{SmePat16} we thus propose to consider 
\begin{align}\label{choice of subset in greedy}
S(R_{n}^{+}(\param))&:= \{ \zeta(\param) \in R_{n}^{+}(\param)\,:\,  \|\zeta(\param)\|_{R_{n}^{+}(\param)} \leq 1\} \\
\nonumber\enspace \text{with} \enspace \|\zeta(\param)\|_{R_{n}^{+}(\param)}&:= \left(\sum_{i=1}^{n_{+}} (\boldsymbol{\zeta}_{i}(\param))^{2}\right)^{1/2}
\end{align}
where $\zeta(\param) = \sum_{i=1}^{n_{+}} \boldsymbol{\zeta}_{i}(\param) \chi_{i}(\param)$, $n_{+}:=\dim(R_{n}^{+}(\param))$ and here and henceforth $\{\chi_{i}(\param)\}_{i=1}^{n_{+}}$ denotes the orthonormal basis of $R_{n}^{+}(\param)$.
Note that we are therefore considering a weighted norm in $R_{n}^{+}(\param)$. The deviation $E(S(R_{n}^{+}(\param)),R_N)$ can then be computed by solving the eigenvalue problem: Find $(\boldsymbol{\varrho}_{j}(\param),\sigma_{j}(\param)) \in (\mathbb{R}^{n_{+}},\mathbb{R}^{+})$ such that
\begin{align*}
\mathbb{Z}(\param)\boldsymbol{\varrho}_{j}(\param) &= \sigma_{j}(\param) \boldsymbol{\varrho}_{j}(\param),\\
\text{where} \enspace \mathbb{Z}_{i,l}(\param)&:= (\chi_{l}(\param) - \sum_{k=1}^{N}(\chi_{l}(\param),\chi_{k})_{R}\chi_{k},\chi_{i}(\param) - \sum_{k=1}^{N}(\chi_{i}(\param),\chi_{k})_{R}\chi_{k})_{R}
\end{align*}
and $\chi_{k}$ denotes the orthonormal basis of $R_{N}$. We thus obtain 
$
E(S(R_{n}^{+}(\param)),R_N) = \sqrt{\sigma_{1}(\param)}
$, for all $\param \in \Xi$,
and $\kappa = \sum_{i=1}^{n_{+}} \boldsymbol{\varrho}_{1}(\param^{*})\chi_{i}(\param^*)$ at each iteration. 

Note that were we to consider the norm $\| \cdot \|_{R}$ in \eqref{choice of subset in greedy} the sorting of the spectral basis $\chi_{i}(\param)$ of $R_{n}^{+}(\param)$ in terms of approximation properties is neglected in the {\tt while} loop of Algorithm \ref{algo:spectral greedy}; for further explanations see \cite{SmePat16}.

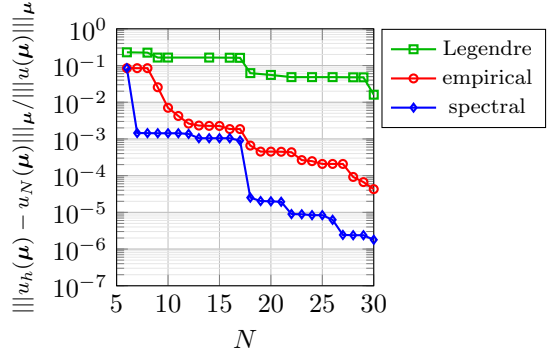
\begin{figure}[t]
\centering
\begin{tikzpicture}
\begin{semilogyaxis}[
    width=5cm,
    height=5cm,
    xmin=5,
    xmax=30,
    ymin=1e-7,
    ymax=1,
    xlabel=$N$,
    ylabel={\footnotesize $\energynorm{u_{h}(\param) - u_{N}(\param)}{\param}/\energynorm{u(\param)}{\param}$},
    grid=both,
    grid style={line width=.1pt, draw=gray!20},
    major grid style={line width=.2pt,draw=gray!50},
    xtick={1,5,10,15,20,25,30},
    ytick={1e-7, 1e-6,1e-5, 1e-4, 1e-3, 1e-2, 1e-1, 1e0},
    legend pos=outer north east,
  ]
    \addplot+[solid, gruen, thick, mark=square, mark size=1.5pt] table[x index=0,y index=1] {beam_legendre.dat};
    \addplot+[solid, red, thick, mark=o, mark size=1.5pt] table[x index=0,y index=1] {beam_empirical.dat};
    \addplot+[solid, blue, thick, mark=diamond, mark size=1.5pt] table[x index=0,y index=2] {beam_empirical.dat};
  \legend{
  {\footnotesize Legendre}, 
  {\footnotesize empirical}, 
  {\footnotesize spectral} };
\end{semilogyaxis}
\end{tikzpicture}
\caption{$\energynorm{u_{h}(\param) - u^{N}(\param)}{\param}/\energynorm{u(\param)}{\param}$ for the Legendre, empirical, and spectral interface basis functions for the solid beam. \label{fig beam error}}
\end{figure}

Finally, we compare in Fig. \ref{fig beam error} the spectral modes generated by the spectral greedy algorithm \ref{algo:spectral greedy} numerically with other interface modes, demonstrating the superior convergence of the former. In detail, we compare the relative error of the port-reduced static condensation approximation for interface spaces comprising ``Legendre polynomial''-type functions\footnote{Note that each component of the displacement is the solution of a scalar singular Sturm-Liouville eigenproblem.} \cite{EftPat13a}, empirical port modes constructed by a pairwise training algorithm\footnote{Following the notation in \cite{EftPat13b} we have chosen $N_{\text{samples}} =500$ and $\gamma = 3$ in the pairwise training algorithm.} \cite{EftPat13a,EftPat13b}, and the spectral modes. To that end, we consider a domain $\Omega$ which consists of two identical solid beams, each of whom is associated with a subdomain $\Omega_{i}$, $i=1,2$. Here, we choose $\Omega_{1}=(-0.5,0.5)\times (-0.5,0.5) \times (0,5)$, $\Omega_{2}=(-0.5,0.5)\times (-0.5,0.5) \times (5,10)$ and $\Gamma_{out} = \Gamma_{1} \cup \Gamma_{2}$, with $\Gamma_{1}=(-0.5,0.5)\times (-0.5,0.5) \times \{0\}$ and $\Gamma_{2}=(-0.5,0.5)\times (-0.5,0.5) \times \{10\}$. The underlying FE discretization has $N=3348$ degrees of freedom per subdomain and $N_{\Gamma}=108$ degrees of freedom per interface. We require $E(\param)$ to be uniform within each subdomain, the constant varying in $[1,10]$ and choose for $\boldsymbol{G}(\param) \in \mathbb{R}^{3}$ the admissible set of parameters to be $[-1,1]\times [-1,1]\times [-1,1]$. Finally, we equip both $S$ and $R$ again with a lifting inner product. Within the spectral greedy we have considered $200$ parameter values sampled from the uniform distribution over $\mathcal{P}$ and $\varepsilon = 1 \cdot 10^{-6}$. On average the interface spaces $R_{n}^{+}(\param)$ have had a size of $13.65$ and the resulting parameter-independent port space $R_{N}$ has a size of $56$. 

In the online stage we consider $E(\param)\equiv 1$ in both components, $\boldsymbol{G}=(0,0,0)^{T}$, and prescribe $\boldsymbol{g}_{D,1}=(0,0,0)^{T}$ at $\Gamma_{1}$ and $\boldsymbol{g}_{D,2}=(1,1,1)^{T}$ at $\Gamma_{2}$. We observe that the Legendre modes perform by far the worst, demonstrating that including information on the solution manifold in the basis construction procedure can significantly improve the approximation behavior. We remark that the Legendre modes will perform even worse in the case of less regular behaviour on the interface, which further justifies the need for problem-specific local approximation spaces in the sense of model reduction. The empirical modes and spectral modes exhibit a comparable convergence until $N=17$, but for $N>17$ the relative error in the spectral approximation is one order of magnitude smaller than that of the empirical port mode approximation. This can be explained by the fact that thanks to its conception the pairwise training algorithm is able to identify and include the most significant modes, but (in contrast to the spectral greedy algorithm) might have difficulties to detect subtle modes that affect the shape of the function at the interface $\Gamma_{m,m'}$ only slightly. Note that the temporary stagnation of the relative error for $N=7,...,17$ for the spectral modes is due to the fact that the spectral greedy prepares the interface space for all possible boundary conditions and parameter configurations. Thus, for the boundary conditions considered here some spectral modes, as say a mode related to a twisting (torsion) of the beam, are not needed for the approximation.


\section{A posteriori error estimation}\label{sec:apost}\label{sec:a_posteriori}

\subsection{Residual based a posteriori error estimation} 

A global residual based a posteriori error estimator for projection based model reduction is readily defined as
\begin{equation}\label{eq:error_global}
\Delta(u_N(\param)) := \frac{1}{\alpha(\param)} \norm{R(u_N(\param);\param)}_{V_h'}
\end{equation}
where $R(u_N(\param);\param) \in V_h^\prime$ is the global residual given as
$
	\dualpair{R(u_N(\param);\param)}{\varphi_h} = f_h({\varphi_h;\param}) - a_h(u_N(\param), \varphi_h; \param)
$ for all $\varphi_h \in V_h.$
This error estimator is known to be robust and efficient (cf. \cite[Proposition 4.4]{HesRozSta2016}), i.e. we have
\begin{equation}\label{eq:global_estimator}
\norm{u_h(\param) - u_N(\param)}_V \leq \Delta(u_N(\param)) 
\leq \frac{\gamma(\param)}{\alpha(\param)} \norm{u_h(\param) - u_N(\param)}_V.
\end{equation}

For localized model order reduction, however, we are merely interested in localized a posteriori error estimation. 
To this end, we first present abstract localized lower and upper bounds for the dual norm of a linear functional 
(see \cite{Betal17}).

\begin{theorem}[Localized lower and upper bounds for functionals]
\label{thm:a_posteriori}
Let $\subsod_i$, $1 \leq i \leq \tilde{M}$ be
a collection of linear subspaces of $V_h$, and let $P_{\subsod_i}: V_h \longrightarrow \subsod_i \subseteq V_h$
be mappings which satisfy $\sum_{i=1}^{\tilde{M}} P_{\subsod_i} = \operatorname{id}_{V_h}$.
Moreover, assume that for $J \in \mathbb{N}$ there exists a partition $\dot{\bigcup}_{j=1}^J \Upsilon_j = \{1, \ldots, \tilde{M}\}$
such that for arbitrary $1 \leq j \leq J$ and $i_1 \neq i_2 \in \Upsilon_{j}$ we have $\subsod_{i_1} \perp \subsod_{i_2}$.

Defining the stability constant of this partition modulo $V_N$ as
\begin{equation}\label{eq:def stability constant}
\puconstant := \sup_{\varphi \in V_h \setminus \{0\}} \frac{( \sum_{i = 1}^{\tilde{M}} \inf_{\tilde \varphi \in V_N \cap \subsod_i}  \norm{P_{\subsod_i}(\varphi) - \tilde \varphi}^2)^\frac{1}{2}}{\norm{\varphi}}
\end{equation}
we have for any linear functional $f \in V_h^\prime$ with $ \dualpair{f}{\varphi} = 0 \  \forall {\varphi} \in V_N$  the estimate
\begin{equation}
	\label{eq:abstract_localized_estimate}
	 \frac{1}{\sqrt{J}}  \Big( \sum_{i = 1}^{\tilde{M}} \norm{f}^2_{\subsod_i^\prime} \Big)^{\frac{1}{2}} \leq
     \norm{f}_{V_h^\prime} \leq \puconstant \cdot \Big( \sum_{i = 1}^{\tilde{M}} \norm{f}^2_{\subsod_i^\prime} \Big)^{\frac{1}{2}}.
\end{equation}
Here, $\norm{f}_{\subsod_i^\prime}$ denotes the norm of the restriction of $f$ to $\subsod_i$.
\end{theorem}

When grouping the spaces $\subsod_i$ so that in each group, all
spaces are orthogonal to each other, $J$ is the number of groups needed. Note that subtracting the projection onto $V_N$ in \eqref{eq:def stability constant} allows subtracting say the mean value of a function or the orthogonal projection onto the rigid body motions, if the respective functions are included in $V_N$. We may thus employ say Poincar\'{e}'s inequality or Korn's inequality in subdomains that do not lie at $\Gamma_{D}$. 

Applying both estimates to the residual $R(u_N(\param);\param) \in V_h^\prime$, we obtain from \eqref{eq:error_global} and Theorem \ref{thm:a_posteriori}
a robust and efficient, localized error estimate:
\begin{corollary}[Localized residual based a posteriori error estimate]
\label{thm:abstract_error_estimate}
Let the assumptions on the subspace collection $\subsod_i$ and the mappings $P_i$ from Theorem \ref{thm:a_posteriori} be satisfied. 
Then, the error estimator $\Delta_{loc}(u_N(\param))$ defined as
\begin{equation}
	\label{eq:localized_estimate}
\Delta_{loc}(u_N(\param)) := \frac{1}{\alpha(\param)}  \puconstant  \big( \sum_{i = 1}^{\tilde{M}} \norm{R(u_N(\param);\param)}^2_{\subsod_i^\prime} \big)^\frac{1}{2}
\end{equation}
is robust and efficient, i.e.
\begin{equation}
	\label{eq:efficiency_estimate}
	\norm{u_h(\param) - u_N(\param)}_V \leq \Delta_{loc}(u_N(\param))
	\leq \frac{\gamma(\param)  \sqrt{J}  \puconstant}{\alpha(\param)}   \norm{ u_h(\param) - u_N(\param)}_V.
\end{equation}
\end{corollary}

Online-offline decomposition of this error estimator can be done by 
applying the usual strategy for online-offline decomposition
used with the standard RB error estimator (see e.g. \cite[Sec.  4.2.5]{HesRozSta2016}
or a numerically more stable approach \cite{BEOR14a,CaErLe14,SoFaDy15})
to every dual norm in $\Delta_{loc}(u_N(\param))$.

The a posteriori error estimator for the ArbiLoMod derived in \cite{Betal17} and the a posteriori error estimator for the scRBE method as suggested in \cite{Sme15} both fit into the framework above as will be detailed below in Examples \ref{ex:Arbilomod} and \ref{ex:scrbe}. In contrast, for instance the error estimators proposed in \cite{HuKnPa13,HuKnPa13b} for the scRBE method exploit matrix perturbation analysis at the system level to bound the Euclidean norm of the error between the coefficients of the static condensation solution and the coefficients of the static condensation solution using an RB approximation in the interior. To estimate the error caused by interface reduction in \cite{EftPat13a} a computationally tractable non-conforming approximation to the exact error is employed. To take into account the error due to the intra-element RB approximations ideas from \cite{HuKnPa13} are used. It can also be noted that the error estimators in \cite{HuKnPa13,HuKnPa13b,EftPat13a} are only valid under certain assumptions on the accuracy of the RB approximation. In \cite{MarRozHaa14} a localized a posteriori error estimator for interface reduction and intra-element RB approximation is presented for the coupled Stokes-Darcy system. The a posteriori error estimator for the CMS method derived in \cite{JaBeLa11} employs the dual norms of residuals and eigenvalues of the eigenproblems used for the construction of the (local) basis functions. The error estimator in \cite{JaBeLa11} is however only partially local as it involves the residual for the port or interface space on the whole interface $\Gamma$.
For localized a posteriori error estimation in the context of adaptive GMsFEM we refer to \cite{Chung2014,Chung2014b,Chung2015,chung2015online}.

\begin{example}[Localized a posteriori error estimate for ArbiLoMod \cite{Betal17}]\label{ex:Arbilomod}
Let us assume $V_h \subset V = H^1_0(\Omega)$, and
choose $\subsod_i$ as subspaces of $H^1(\Omega_i)$ where $\EstGrid:=\{\tilde\Omega_1, \ldots,
\tilde\Omega_{\tilde{M}}\}$ is an arbitrary \emph{overlapping} decomposition of $\Omega$, which may be chosen independently from $\Grid$.
Assume that there is a partition of unity $p_i \in H^{1,\infty}(\tilde\Omega_i) \cap C(\tilde\Omega_i)$,
$\sum_{i=1}^{\tilde{M}} p_i = 1$, such that $\norm{p_i}_\infty \leq 1$ and
$\norm{\nabla p_i}_\infty \leq \pufuncconstant \diam(\tilde\Omega_i)^{-1}$.
The constant $\pufuncconstant$ will depend on size of the overlap of the subdomains
$\tilde\Omega_i$ with their neighbors in relation to their diameters.

Moreover, we assume that there is a linear interpolation operator $\mathcal{I}$ onto $V_h$ such that
$\mathcal{I}$ is the identity on $V_h$ with $\mathcal{I}(p_i V_h) \subseteq \subsod_i$ and
$\|\mathcal{I}(p_i v_h) - p_i v_h\|_V \leq c_I \|v_h\|_{\tilde\Omega_i,1}$ for all $v_h \in V_h$.
We then can define mappings
\begin{equation*}
  P_{O_i}(v_h) := \mathcal{I}(p_i\cdot v_h).
\end{equation*}
which satisfy the assumptions of Theorem~\ref{thm:a_posteriori}.
In case $V_h$ comes from a finite element discretization, a possible choice for $\mathcal{I}$ is Lagrange interpolation.

If we now ensure that the partition of unity $p_i$ is included in $V_N$, we can choose
$\tilde\varphi$ in the definition of $\puconstant$ as $\tilde\varphi := p_i \cdot |\tilde\Omega_i|^{-1}\int_{\tilde\Omega_i}\varphi$, which allows us to prove \cite[Proposition 5.7]{Betal17} that
$\puconstant$ can be bounded by
\begin{equation*}
\puconstant \leq \sqrt{4 + 2c_I^2 + 4(\pufuncconstant c_{\mathrm pc})^2} \cdot \sqrt{\ovlpconstant}.
\end{equation*}
In this estimate $\ovlpconstant := \max_{x \in \Omega} \#\{i \in \Upsilon_E \ | \ x \in \Omega_i\}$ 
is the maximum number of estimator domains $\tilde\Omega_i$ overlapping in any point $x$ of $\Omega$,
and $c_{\mathrm pc}$ is a Poincar\'e-inequality constant associated with $\EstGrid$. 
In particular, this result shows that the efficiency of \eqref{eq:localized_estimate} is independent
from the number of subdomains in $\Grid$, provided that the partition of unity $p_i$ is
included in $V_N$.
\end{example}

\begin{example}[scRBE method and interface reduction from \cite{Sme15}]\label{ex:scrbe}
We exemplify the a posteriori error estimator from Corollary \ref{thm:abstract_error_estimate} for the scRBE method, which is equally applicable when considering solely static condensation and no intra-element RB approximations.\footnote{The error estimator in \cite{Sme15} is derived for mutually disjoint interfaces. However, we conjecture that the estimator can be generalized to general decompositions of $\Omega$.} To simplify notations we define interface spaces $V^{\gamma}_{h}:=\spanlin\{\psi_{1}^{\gamma},\hdots,\psi_{N^{\gamma}_{h}}^{\gamma}\}$, where $N^{\gamma}_{h}=\dim(V_{h}|_{\gamma})$; for the definition of $\psi_{k}^{\gamma}$ we refer to subsection \ref{subsec:conforming}. Recall that we then have the following space decomposition of the (global) finite element space $V_{h}$
\begin{equation}\label{eq:scRBE space decomposition}
V_{h} = \bigoplus_{m=1}^{M} V_{h;0}^{m} \oplus \left(\bigoplus_{\gamma \in \mathcal{T}_{H}^{\gamma}} V^{\gamma}_{h}\right). 
\end{equation}
We may thus uniquely rewrite every $\varphi \in V_{h}$ as
\begin{equation}\label{eq:scRBE space decomp func}
\varphi = \sum_{m=1}^{M} \varphi^{m} + \sum_{\gamma \in \mathcal{T}_{H}^{\gamma}} \varphi^{\gamma}
\end{equation}
where $\varphi^{m} \in V_{h;0}^{m}$ and $\varphi^{\gamma}\in V^{\gamma}_{h}$, extending $\varphi^{m}$ and $\varphi^{\gamma}$ by zero. This allows us to define mappings $P_{V_{h;0}^{m}}: V_{h} \rightarrow V_{h;0}^{m}$, $\varphi \mapsto \varphi^{m}$ and $P_{V^{\gamma}_{h}}: V_{h} \rightarrow V^{\gamma}_{h}$, $\varphi \mapsto \varphi^{\gamma}$ as required in Theorem \ref{thm:a_posteriori}. Thanks to  \eqref{eq:def_face_basis_function} we also obtain
\begin{equation*}
V_{h;0}^{m} \perp V_{h;0}^{m'}, \enspace m\neq m' \quad \text{and} \quad V_{h;0}^{m} \perp V^{\gamma}_{h}, \enspace m=1,\hdots, \Omega, \gamma \in \mathcal{T}_{H}^{\gamma}.
\end{equation*}

It thus remains to verify that we can bound the constant $\puconstant$ with $V_{N} $ as defined in \eqref{eq:reduced space conform}. To that end, we first note that thanks to \eqref{eq:def_face_basis_function} we have the following stability result \cite[Proposition 4.1]{Sme15}:
\begin{equation}\label{eq:stable decomp}
\| \varphi \|_{V}^{2} = \sum_{m=1}^{M} \|  \varphi^{m} \|_{V}^{2} + \| \sum_{\gamma \in \mathcal{T}_{H}^{\gamma}} \varphi^{\gamma} \|^{2}_{V}. 
\end{equation}
We thus obtain
\begin{align*}
\puconstant &\leq \sup_{\varphi \in V_h \setminus \{0\}} \frac{\left(\sum_{m=1}^{M} \|  \varphi_{m} \|_{V}^{2} +  \sum_{\gamma \in \mathcal{T}_{H}^{\gamma}} \inf_{\tilde \varphi_{f} \in V_{N}^{\gamma}} \| \varphi^{\gamma} - \tilde \varphi^{\gamma} \|^{2}_{V}\right)^{1/2}}{\|\varphi\|_{V}}\\
&\overset{\eqref{eq:stable decomp}}{\leq} \sup_{\varphi \in V_h \setminus \{0\}} \frac{\left(\| \varphi \|_{V}^{2}
 +  \sum_{\gamma \in \mathcal{T}_{H}^{\gamma}} \inf_{\tilde \varphi^{\gamma} \in V_{N}^{\gamma}} \| \varphi^{\gamma} - \tilde \varphi^{\gamma} \|^{2}_{V}\right)^{1/2}}{\|\varphi\|_{V}},
\end{align*}
where $V^{\gamma}_{N}:=\spanlin\{\psi_{1}^{\gamma},\hdots,\psi_{N^{\gamma}}^{\gamma}\}$. To show $\sum_{\gamma \in \mathcal{T}_{H}^{\gamma}} \inf_{\tilde \varphi^{\gamma} \in V_{N}^{\gamma}} \| \varphi^{\gamma} - \tilde \varphi^{\gamma} \|^{2}_{V} \leq c \|\varphi\|_{V}$ for a constant $c$ we choose $\tilde \varphi^{\gamma}$ such that $(\varphi^{\gamma} - \tilde \varphi^{\gamma})|_{\gamma}$ equals the trace of $\varphi$ minus the orthogonal projection on the kernel of the bilinear form; for further details see \cite{BuhSme18,SmePat16}. Then, we can use \cite[Lemma B.4]{SmePat16} to conclude boundedness of $\puconstant$ and thus \eqref{eq:efficiency_estimate}, the latter corresponding to \cite[Proposition 4.2 and Corollary 4.6]{Sme15}.

Finally, we shortly discuss how to compute the dual norms of the residuals in \eqref{eq:abstract_localized_estimate}. The dual norms of the residuals of the intra-element RB approximations can be computed by employing Riesz representations. The dual norms of the residuals in the interface space can be computed by means of conservative fluxes \cite{HuEnMaLa00}, which have been extended to interface reduction in \cite{Sme15}. In detail, we compute the conservative flux $H^{m}_{N}(\param)$ such that 
\begin{equation}\label{eq:cons fluxes}
\sum_{\gamma \in \overline{\Omega}_{m}} (H^{m}_{N}(\param),\psi^{\gamma})_{\gamma} = f_{m}(\psi^{\gamma};\param) - a_{m}(u_{N}(\param), \psi^{\gamma}; \param) \quad \forall \psi^{\gamma} \in \bigoplus_{\gamma \in \overline{\Omega}_{m}} V^{\gamma}_{h},
\end{equation}
where $(\cdot, \cdot)_{\gamma}$ denotes a suitable inner product on the interface $\gamma$. Note that thanks to our mutual disjoint interface assumption problem \eqref{eq:cons fluxes} decouples and we may compute the conservative flux separately for each interface $\gamma$. Moreover, by orthonormalizing the interface basis functions $\chi^{\gamma}_{k}$ defined in subsection \ref{subsec:conforming} w.r.t the $(\cdot, \cdot)_{\gamma}$ inner product, the computation of $H^{m}_{N}(\param)$ reduces to the assembling of the residual in \eqref{eq:cons fluxes}. The computational costs thus scale linearly in $(N^{\gamma}_{h} - N^{\gamma})$ and $N^{\gamma}$. For further details we refer to \cite{Sme15}.
\end{example}


\subsection{Local flux reconstruction based error estimation}\label{sec:flux-reconstruction} 

Following \cite{OS2015}, we discuss local flux reconstruction based a posteriori error estimation of the full approximation error $u(\param) - u_N(\param)$ (that is: the discretization as well as the model reduction error) in the context of non conforming approximations of elliptic multiscale problems such as Example \ref{ex:multiscale_problem}.
An extension to convection--diffusion--reaction problems based on \cite{ESV2010} is straightforward.
This estimate was introduced in the IP localized non conforming setting of the LRBMS (compare Section \ref{sec:ip-nonconforming}).

Recalling the broken Sobolev space and broken gradient operator from Section \ref{sec:ip-nonconforming}, the key idea of flux reconstruction based error estimation is to observe that not only the approximate solution $u_N(\param)$ is non conforming, but also the approximate diffusive flux $-\kappa(\param)\gradienth u_N(\param)$, in the sense that it is not contained in $\Hdiv(\Omega)$ (i.e., the space of functions in $L^2(\Omega)^d$ whose divergence exists in a weak sense and lies in $L^2(\Omega)$). 

We may then obtain computable estimates by comparing these quantities with conforming reconstructions, as detailed further below. 
\footnote{%
  Note that all of the analysis holds for the FOM solution $u_h(\param)$ as well as the ROM solution $u_N(\param)$ (compare \cite{OS2015}), but we restrict the exposition to the latter.
  In particular, the presented estimates can thus also be used to steer grid adaptation of the FOM solution.
}
The respective reconstructed diffusive flux is locally conservative 
and is related to the conservative flux reconstruction to compute the dual norm of the residuals in the interface space in Example \ref{ex:scrbe}.

To begin with, we specify the parameter dependent (semi-)energy norm induced by the bilinear form $a$ for a parameter $\paramFixed \in \Params$, $\energynorm{\cdot}{\paramFixed} : H^1(\grid) \to \R$, $v \mapsto \energynorm{v}{\paramFixed} := a(v, v; \paramFixed)^{\frac{1}{2}},$ (by using the broken gradient in the definition of $a$) and note that we can compare these semi norms for two parameters by means of the affine decomposition of $a$ (compare \eqref{eq:affine_decomp}),
\begin{align}
  \underline{\Theta_a}(\param, \paramFixed)^{1/2}\; \energynorm{v}{\paramFixed}\; \leq\; \energynorm{v}{\param}\;\; \leq\;\; \overline{\Theta_a}(\param, \paramFixed)^{\frac{1}{2}}\; \energynorm{v}{\paramFixed},
  \notag
\end{align}
with the equivalence constants given by $\underline{\Theta_a}(\param, \paramFixed) := \min_{q = 1}^{Q_a} \Theta_a^q(\param)\, \Theta_a^q(\paramFixed)^{-1}$ and $\overline{\Theta_a}(\param, \paramFixed) := \max_{q = 1}^{Q_a} \Theta_a^q(\param)\, \Theta_a^q(\paramFixed)^{-1}$, respectively.
The first abstract result is the following discretization-agnostic lemma, which leaves the choice of the reconstructions, $v$ and $s$, open.
(We give estimates on the full $V_h$-norm at the end of this subsection.)

\begin{lemma}[Abstract energy norm estimate (Lemma 4.1 in \cite{OS2015})]
  \label{lem:lrbms_abstract_estimate}
  For $\param \in \Params$, let $u(\param) \in V$ denote the weak solution of \eqref{eq:weak_solution} with the data functions $\kappa$ and $q$ as in example \ref{ex:multiscale_problem}.
  It then holds for arbitrary $v_N \in H^1(\grid)$ and $\paramFixed \in \Params$, that
  \begin{align}
    & \hspace*{-1em}\energynorm{u(\param) - v_N}{\paramFixed} 
  \notag\\
    &\leq \underline{\Theta_a}(\param, \paramFixed)^{-\frac{1}{2}} \Big\{
  \overline{\Theta_a}(\param, \paramFixed)^{\frac{1}{2}} \inf_{v \in V} \energynorm{u(\param) - v}{\paramFixed} 
  \notag\\
    &\hspace{10pt}+\hnS\inf_{s \in \Hdiv}\hnS\Big(
    \hnS\sup_{\substack{\varphi\in V\\|\varphi|_{a; \param}=1}}\hnS\big\{
      (q - \divergence s, \varphi)_{L^2(\Omega)}
    \,-\,(\kappa(\param) \gradienth v_N + s, \gradient\varphi)_{L^2(\Omega)}
    \big\}
    \Big)
    \Big\}
  \notag\\
    &\leq \tfrac{\overline{\Theta_a}(\param, \paramFixed)}{\underline{\Theta_a}(\param, \paramFixed)}^{\frac{1}{2}}\, 2\, \energynorm{u(\param) - v_N}{\paramFixed}.
  \notag
  \end{align}
\end{lemma}

To obtain a fully computable localizable estimate we need to specify the conforming reconstruction of the solution ($v$ in the above lemma) and of the diffusive flux ($s$ in the above lemma).
We define both reconstructions w.r.t.~the global fine grid $\grid$ and note, that their respective computations can be localized w.r.t.~the domain decomposition to allow for offline/online decomposable localized estimates.

We reconstruct the non conforming solution $u_h(\param) \in V_h$ by means of its \emph{Oswald interpolant} $I_\text{OS}[u_h(\param)] \in V$.
We define the corresponding Oswald interpolation operator $I_\text{OS}: V_h \to V_h \cap V$ by specifying its values on each Lagrange node $\nu$ of $\grid$: given any $v_h \in V_h$, we set $I_\text{OS}[v_h](\nu) := v_h|_t(\nu)$ for any Lagrange node lying inside a grid element $t \in \grid$,
\begin{align}
  I_\text{OS}[v_h](\nu) := 0 &&\text{for all boundary nodes and}&& I_\text{OS}[v_h](\nu) := \tfrac{1}{|\grid^\nu|}\sum_{t \in \grid^\nu}v_h|_t(\nu)
\notag
\end{align}
for all nodes which are shared by multiple grid elements, which we collect in $\grid^\nu \subset \grid$.

The definition of the conforming reconstruction of the non conforming diffusive flux $-\kappa(\param)\gradienth u_h(\param) \in L^2(\Omega)^d$ is more involved.
Given $l \geq 0$, we define the $l$th order \emph{Raviart-Thomas-N\'ed\'elec} space of vector valued functions by
\begin{align}
  RTN_h^l(\grid) := \big\{ s \in \Hdiv(\Omega) \;\big|\; s|_t \in [\mathbb{P}_l(t)]^d + x \mathbb{P}_l(t) \;\;\forall t \in \grid\big\}
\notag
\end{align}
and note that the DoFs of any $s_h \in RTN_h^l(\grid)$ are uniquely defined by specifying the moments of order up to $l - 1$ of $s_h|_t$ on all elements $t \in \grid$ and the moments of order up to $l$ of $s_h|_\sigmafine \cdot n_\sigmafine$ on all faces $\sigmafine \in \faces$ (compare \cite{BF1991}).
With these preliminaries we define the \emph{diffusive flux reconstruction operator} $R_h^l: \Params \to [V_h \to RTN_h^l(\grid)]$, given some $v_h \in V_h$ and some $\param \in \Params$ by specifying the DoFs of $R_h^l[v_h; \param] \in RTN_h^l(\grid)$, such that
\begin{align}
  \big(R_h^l[v_h; \param]\cdot n_\sigmafine, r\big)_{L^2(\sigmafine)} &= a_\sigmafine^c(v_h, r; \param) + (v_h, r)_\sigmafine^p \quad\quad\text{for all } r \in \mathbb{P}_l(\sigmafine),
\label{eq::flux_recostruction::1}
\intertext{%
  on all $\sigmafine \in \faces$ and
}
  \big(R_h^l[v_h; \param], \gradient r\big)_{L^2(t)} &= -a^\text{CG}(R_h^l[v_h; \param]\big|_t, r; \param) - \sum_{\sigmafine \in \faces \cap t} a_\sigmafine^c(r, v_h; \param)
\label{eq::flux_recostruction::2}
\end{align}
for all $\gradient r \in [\mathbb{P}_{l - 1}(t)]^d$ with $r \in \mathbb{P}_l(t)$ on all $t \in \grid$.
Given a FOM space $V_h$ of polynomial order $k \geq 1$, we choose a $k - 1$st order reconstruction.
With this definition, the reconstructed diffusive flux of a given a reduced solution $u_N(\param)$ fulfills the following \emph{local conservation} property, given that the constant function $1$ is present in the local reduced spaces $V_N^{m}$:
\begin{align}
  \big(\divergence{R_h^{k - 1}[u_N(\param); \param]}, 1 \big)_{L^2(\Omega_m)} = \big(q, 1\big)_{L^2(\Omega_m)}, &&\text{for all }\Omega_m \in \Grid.
\notag
\end{align}
When inserting this diffusive flux reconstruction for $s$ in lemma \ref{lem:lrbms_abstract_estimate}, this local conservation property is key to obtaining the following estimate.

\begin{theorem}[Locally computable energy norm a posteriori estimate]
  \label{thm:flux_reconstruction}
  Let the domain decomposition $\Grid$ from Definition \ref{def:decompoisition} be such, that the Poincar\'e-inequality holds on each subdomain $\Omega_m \in \Grid$ with a constant $C_P^{m} > 0$,
  \begin{align}
    \|\varphi - \Pi_0^{m}\varphi\|_{L^2(\Omega_m)}^2 \leq C_P^{m}\, h_{m}^2\, \|\gradient\varphi\|_{L^2(\Omega_m)}^2&&\text{for all }\varphi \in H^1(\Omega_m),
  \notag
  \end{align}
  where $h_{m} := \diam(\Omega_m)$ and where $\Pi_0^{m}\varphi$ denotes the mean value of $\varphi$ over $\Omega_m$.
  Let further $u(\param) \in V$ be the weak solution of \eqref{eq:parametric-multiscale} and $u_N(\param) \in V_N$ be the IP localized ROM solution, with $1 \in V_N^{m}$ for $1 \leq m \leq M$.
  It then holds for arbitrary $\paramFixed, \paramHat \in \Params$ that
  \begin{align}
    \energynorm{u(\param) - u_N(\param)}{ \paramFixed} \leq \eta(\param; \paramFixed; \paramHat)
  \notag
  \end{align}
  with the a posterior error estimator $\eta(\param; \paramFixed; \paramHat)$ given by
  \begin{align}
    \eta(\param; \paramFixed; \paramHat) := \underline{\Theta_a}(\param, \paramFixed)^{-\frac{1}{2}}
      \Big[
        &\overline{\Theta_a}(\param, \paramFixed)^{\frac{1}{2}} \Big( \sum_{\Omega_m \in \Grid} \eta_\textnormal{nc}^{\Omega_m}(\param; \paramFixed)^2 \Big)^{\frac{1}{2}}
    \notag\\
        &+\Big(
        \sum_{\Omega_m \in \Grid}\big(\eta_\textnormal{r}^{\Omega_m}(\param) + \underline{\Theta_a}(\param, \paramHat)^{-1}\; \eta_\textnormal{df}(\param; \paramHat)\big)^2
        \Big)^{\frac{1}{2}}
      \Big],
  \notag
  \end{align}
  and the local non conformity, residual and diffusive flux indicators given by
  \begin{align}
    \eta_\textnormal{nc}^{\Omega_m}(\param; \paramFixed) &:= \big|\big(v_N(\param) - I_\textnormal{OS}[v_N(\param)\big)|_{\Omega_m}\big|_{a; \paramFixed},
  \notag\\
    \eta_\textnormal{r}^{\Omega_m}(\param) &:= \tfrac{C_{\Omega_m}^P}{\underline{\kappa}_{\Omega_m}}^{\frac{1}{2}} \big\| q - \divergence{R_h^{k - 1}[u_N(\param);\param]}\big\|_{L^2(\Omega_m)}\text{ and}
  \notag\\
  \eta_\textnormal{df}(\param; \paramHat)\big) &:= \big\|\kappa(\paramHat)^{-1}\big(\kappa(\param)\gradienth u_N(\param) + R_h^{k -1}[u_N(\param);\param]\big)\big\|_{L^2(\Omega_m)}
  \end{align}
  respectively, where $\underline{\kappa}_{\Omega_m}$ denotes the minimum eigenvalue of $\kappa$ over $\Omega_m$ and $\Params$.
\end{theorem}

We obtain an a posterior error estimate w.r.t.~the $V_h$-norm or a full energy norm, $\energynorm{\cdot}{\param} + \big(\sum_{\sigmafine \in \faces} (\cdot, \cdot)_\sigmafine^p\big)^{\frac{1}{2}}$, by noting that $\big(u(\param), u(\param)\big)_\sigmafine^p = 0$ for a weak solution $u(\param)$ of sufficient regularity.


\section{Basis enrichment and online adaptivity} \label{sec:adaptivity} 
Model order reduction is usually employed either (i) in the context of real-time decision making and embedded devices, or (ii) in the context of outer-loop applications, such as optimal control, inverse problems or Monte Carlo methods.
In (i), one is usually interested in reduced spaces $V_N$ of very low dimension to obtain ROMs as small as possible, at the possible expense of very involved offline computations.
Here, localized model order reduction may help to reduce the latter, but we can usually not expect the resulting reduced space to be smaller than the one generated using traditional global model order reduction methods.
In (ii), however, one is interested in a black-box-like approximation scheme which is queried for a huge amount of parameters, with a somehow ``optimal'' computational cost (including offline as well as online cost).
Here, one may keep high-dimensional data throughout the computational process (offline as well as online), and it is in this context that localized model order reduction techniques may truly outperform other approaches.
In the context of PDE constrained optimization this has been investigated e.g. in \cite{OS2017,OSS2018,FR2018}. 

The localized a posteriori error estimation as discussed in Section \ref{sec:apost} enables adaptive enrichment
of the local reduced approximation spaces, whenever the quality of the reduced scheme is estimated to be insufficient -- be it due to insufficient training due to lacking computational resources or due to limited knowledge about the range of possible parameters or due to other reasons altogether.

Let us thus assume that an initial (possibly empty) localized reduced approximation space $V_N$
is given, compare Section \ref{sec:prep}.
The goal of an adaptive enrichment is to enlarge the local solution spaces with additional modes that reflect non-local influences of the
true solution such as channeling effect or singularities.
Local adaptive basis enrichment can be employed both offline for the whole parameter range and/or online for a specific chosen parameter. 
Empirical training followed by offline enrichment is e.g. used in a Greedy manner for the basis construction in ArbiLoMod (cf. Example \ref{ex:Arbilomod}) in \cite{Betal17}. 
Adaptive enrichment for the GMsFEM is presented in \cite{Chung2014,Chung2014b} and online adaptive enrichment in \cite{Chung2015,chung2015online}.
For the exposition in this section, we restrict to online enrichment as introduced \cite{OS2015}, i.e. for local enrichment of the basis when a certain parameter is already chosen.

From a birds-eye perspective, we can think of an online adaptive reduced scheme as a $p$-adaptive FE scheme with problem adapted basis functions, where the local reduced bases are adapted during online enrichment.\footnote{%
  We would also like to mention the $h$-adaptive model order reduction approach from \cite{Car2015} which is based on a $k$-means clustering of the DoFs, but we restrict the exposition here to localization w.r.t.~a domain decomposition.
}
Thus, we can think of online enrichment in the usual \underline{S}olve $\to$ \underline{E}stimate $\to$ \underline{M}ark $\to$ \underline{R}efine (SEMR) manner, well known in grid-adaptive discretization schemes. In the \underline{E}stimate step we employ an a posteriori error estimate $\eta$ that is localizable w.r.t.~the domain decomposition, 
i.e. $\eta^2 \leq \sum_{m = 1}^M \eta_m^2$, with appropriate local indicators $\eta_m$. Examples are given in Section \ref{sec:a_posteriori}.
As such, most marking strategies from grid-adaptive schemes are applicable, and we give examples in Section \ref{sec:multiscale_experiments}.
In this context, refinement is locally done by enriching the local reduced spaces, that is: by adding additional basis functions to the local reduced bases on selected subdomains.
We thus presume we are given a parameter $\param \in \Params$ and a reduced solution $u_N(\param) \in V_N$, the estimated error of which is above a given tolerance.

\begin{algorithm2e}[t]
\DontPrintSemicolon
\SetAlgoVlined
\SetKwInOut{Input}{Input}
\SetKwInOut{Output}{Output}
\Input{a marking strategy \texttt{MARK}, an orthonormalization procedure \texttt{ONB}, a localizable offline/online decomposable a posteriori error estimate $\eta(\param)^2 \leq \sum_{m = 1}^M \eta_m(\param)^2$, local reduced bases $\varPhi^m$ for $1 \leq m \leq M$, $\param \in \Params$, $u_N(\param)$, $\Delta_\text{online} > 0$}
\Output{Updated reduced solution}
${\varPhi^m}^{(0)} \gets \varPhi^m$, $\forall 1 \leq m \leq M$ \;
$n \gets 0$\;
\While{$\eta(\param) > \Delta_\text{online}$}{%
  \ForAll{$1 \leq m \leq M$}{compute local error indicator $\eta_m(\param)$}
  $\widetilde{\Grid} \gets \texttt{MARK}\big( \Grid, \{\eta_m(\param)\}_{1 \leq m \leq M} \big)$\;
  \ForAll{$\Omega_m \in \widetilde{\Grid}$}{%
    Solve \eqref{eq::oversampled_problem} for $\varphi^{\tilde{\Omega}_m}$
    ${\varPhi^m}^{(n+1)} \gets \texttt{ONB}\big(\{{\varPhi^m}^{(n)}, \big(\varphi^{\tilde{\Omega}_m} + u_N(\param)\big)|_{\Omega_m}\} \big)$
  }
  update all reduced quantities (system matrices, error estimates) w.r.t.~the newly added basis elements\;
  solve \eqref{eq:rom_solution} for the reduced solution $u_N(\param)$ using the updated quantities
}
\Return $u_N(\param)$
\caption{Adaptive online enrichment in the context of the LRBMS%
.%
\label{alg:online_adaptive_LRBMS}}
\end{algorithm2e}

As an example, we detail the online enrichment procedure used in the context of the LRBMS (compare Section \ref{sec:ip-nonconforming}),
using the a posteriori error estimation techniques from Section \ref{sec:flux-reconstruction} .
Inspired by domain decomposition as well as numerical multiscale methods, we may then obtain a candidate for the next element of a local reduced basis by solving local corrector problems on a collection $\widetilde{\Grid}\subseteq\Grid$ of marked subdomains with $u_N(\param)$ as boundary values.
For each marked subdomain $\Omega_m \in \widetilde{\Grid}$, we denote by $\tilde{\Omega}_m := \big\{ \Omega_{m'} \in \Grid \;\big|\; \Omega_m \cap \Omega_{m'} \neq \emptyset \big\}$ an overlapping subdomain and by $V_h^{\tilde{\Omega}_m} := \big\{ v|_{\tilde{\Omega}_m} \;\big|\; v \in V_h, v|_{\partial\tilde{\Omega}_m} = 0\big\}$ the associated restricted FOM space, encoding zero Dirichlet boundary values.
We are then looking for a local correction $\varphi^{\tilde{\Omega}_m} \in V_h^{\tilde{\Omega}_m}$, such that
\begin{align}
  a_h(\varphi^{\tilde{\Omega}_m}, v_h; \param) = f_h(v_h; \param) - a_h(u_N(\param)|_{\tilde{\Omega}_m}, v_h; \param) &&\text{for all } v_h \in V_h^{\tilde{\Omega}_m},
  \label{eq::oversampled_problem}
\end{align}
where we understand all quantities to be implicitly extended to $\Omega$ by zero, if required, and note that $\varphi^{\tilde{\Omega}_m}$ can be computed involving only quantities associated with $\tilde{\Omega}_m$.
Using this local correction on the overlapping subdomain, we obtain the next element of the local reduced basis associated with $\Omega_m$ by an orthonormalization of $\big(\varphi^{\tilde{\Omega}_m} + u_N(\param)\big)|_{\Omega_m}$ with respect to the existing basis on $V_N^m$.

Given a marking strategy and an orthonormalization procedure, we summarize the adaptive online enrichment used in the context of the LRBMS 
in Algorithm \ref{alg:online_adaptive_LRBMS}.


\section{Computational aspects} \label{sec:computation}
In this section we discuss the computational efficiency of localized model order reduction
schemes in comparison to standard, non-localized techniques.
Imposing a localization constraint on the reduced space naturally yields sub-optimal
spaces in the sense of Kolmogorov $N$-width. However, this is mitigated by the sparse structure
of the resulting reduced system matrices.
In particular, for problems with large-dimensional parameter domains with localized influence
of each parameter component on the solution, we can expect localized ROMs to show comparable or
even better online efficiency in comparison to a standard ROM\@. In addition, localized model order reduction 
provides more flexibility to balance computational and storage requirements between the offline 
and online phase and has thus the potential to be optimized with respect to the specific needs.
This is particularly favorable for large-scale or multiscale problems, where global snapshot 
computations are extremely costly or even prohibitive. 

In the offline (and enrichment) phase of the localized schemes, only relatively small-dimensional local
problems are solved instead of the computation of global solution snapshots. 
In comparison to a global reduction approach with a parallel solver for snapshot
generation (e.g.\ a domain decomposition scheme), the preparation of the local reduced spaces via training
(Section~\ref{sec:prep}) can be performed almost communication-free, allowing the application of these schemes on
parallel compute architecture without fast interconnect such as cloud environments.
Via adaptive enrichment of the approximation spaces -- based on the solution of local correction problems 
(Section~\ref{sec:adaptivity}) -- smaller and more efficient ROMs can be obtained.
In comparison to domain-decomposition methods, where similar correction problems are solved,
these correction problems are only solved in regions of the domain where the approximation space
is insufficient.
Thus, for problems with a localized effect of the parameterization, a significant reduction of the computational
effort can be expected in the reduced basis generation process.

In the context of component-based localized model order reduction (e.g.\ CMS, scRBE, RBHM, RDF) large computational savings can be achieved
by the preparation of local approximation spaces (components) w.r.t.\ arbitrary neighboring components
(connected through so-called ports).
In addition to parametric changes of the governing equations or computational domain, this allows the (non-parametric)
recombination of components in arbitrary new configurations without requiring additional offline computations.

\subsection{Online efficiency}
In view of Definition~\ref{def:rom}, we can interpret the localized model order reduction
methods introduced in Section~\ref{sec:coupling} as standard projection based model reduction methods -- such as the reduced basis
method -- subject to the constraint that the reduced space $V_N$ admits a localizing decomposition of
the form \eqref{eq:localized_V_N}. 
As such, the usual offline/online decomposition methodology can be applied.
To this end, let us assume that the bilinear form $a(\cdot, \cdot; \param)$ and the source functional 
$f(\cdot; \param)$ admit affine decompositions
\begin{equation}\label{eq:affine_decomp}
    a(v, w; \param) = \sum_{q=1}^{Q_a} \Theta^a_q(\param)a^q(v; w), \qquad
    f(w; \param) = \sum_{q=1}^{Q_f} \Theta^f_q(\param)f^q(w),
\end{equation}
for all $v, w \in V$, $\param \in \Params$ with non-parametric bilinear forms $a^q: V \times V \to
\mathbb{R}$, functionals $f^q \in V^\prime$ and some parameter functionals $\Theta^a_q,
\Theta^f_q: \Params \to \mathbb{R}$.
If the given problem is not of the form \eqref{eq:affine_decomp}, we can employ empirical
interpolation \cite{BMNP04} to compute an approximate affine decomposition.

We begin by computing the reduced approximation space $V_N$ using the methods
outlined in Section~\ref{sec:prep}.
After that, a reduced model is assembled by computing matrix representations
$\mathbb{A}^q \in \mathbb{R}^{N\times N}$ of $a^q$ and vector representations
$\mathbb{F}^q \in \mathbb{R}^N$ of $f^q$ w.r.t.\ a given basis $\varphi_1,\ldots,\varphi_N$ of $V_N$, i.e.
\begin{equation}\label{eq:fom_matrices}
    \mathbb{A}^q_{ij}:= a^q(v_j, w_i), \qquad \mathbb{F}^q_i:=f^q(w_i).
\end{equation}
After this computationally demanding offline phase, the coordinate representation
$\mathbb{U}_N(\param) \in \mathbb{R}^N$ of the reduced solution $u_N(\param)$ of \eqref{eq:rom_solution} is quickly
obtained for arbitrary new parameters $\param$ by solving
\begin{equation}\label{eq:rom_in_matrices}
    \sum_{q=1}^{Q_a} \Theta^a_q(\param)\mathbb{A}^q \cdot \mathbb{U}_N(\param)
    = \sum_{q=1}^{Q_f} \Theta^f_q(\param)\mathbb{F}^q
\end{equation}
in the following online phase.
The computational effort to determine $u_N(\param)$ is of order
\begin{equation}\label{eq:rom_effort}
    \mathcal{O}(Q_aN^2 + Q_fN) + \mathcal{O}(N^3)
\end{equation}
for the assembly and solution of the dense equation system \eqref{eq:rom_in_matrices}.
In particular, we have obtained full offline/online splitting, i.e. the effort to obtain
$\mathbb{U}_N(\param)$ is independent of $\dim V_h$. 
From $\mathbb{U}_N(\param)$ we can then either reconstruct $u_N(\param)$ by linear combination with
the reduced basis or evaluate arbitrary linear functionals of $u_N(\param)$ by additionally
computing vector representations of these functionals in the offline phase.

RB methods aim at constructing reduced spaces $V_N$ which are near-optimal approximation spaces for
the discrete solution manifold $\{u_h(\param)\,|\,\param\in\Params\}$ in the sense of Kolmogorov, i.e. it should
hold that
\begin{equation}\label{eq:nwidth_rb}
    \sup_{\param\in\Params}\inf_{v\in V_N}\|u_h(\mu)  - v\| \approx
    d_N := \inf_{\substack{W \subseteq V_h\\\dim W = N}}\sup_{\param\in\Params}\inf_{v\in W}\|u_h(\mu) - v\|,
\end{equation}
where $d_N$ is the Kolmogorov $N$-width of the solution manifold.
Localized RB methods aim at reducing the computational effort of the offline phase by replacing the
computation of solution snapshots $u_h(\param)$ of the global discrete full order model by solutions
of smaller localized problems associated with the domain $\Grid$ (see below).
This comes at the expense of replacing the set of all $N$-dimensional subspaces of $V_h$
by the smaller set of all $N$-dimensional subspaces of $V_h$ of the form \eqref{eq:localized_V_N},
i.e. we aim at constructing $V_N$ with
\begin{equation}\label{eq:nwidth_lrb}
    \sup_{\param\in\Params}\inf_{v\in V_N}\|u_h(\mu)  - v\| \approx
    d_N^{loc} := \inf_{\substack{W \subseteq V_h\\\dim W = N\\W\text{ satisf.
    \eqref{eq:localized_V_N}}}}\sup_{\param\in\Params}\inf_{v\in W}\|u_h(\mu) - v\|.
\end{equation}

As $d_N^{loc} > d_N$, localized RB methods generally result in larger $V_N$ to satisfy a given
approximation error tolerance $\varepsilon$.
Since we can represent any basis vector of a global RB approximation of \eqref{eq:fom_solution}
w.r.t.\ the localizing space decomposition \eqref{eq:localizing_space_decomposition} as a sum of $M^{tot}
:= M + \#\Grid^\gamma + \#\Grid^e + \#\Grid^v$ local vectors, we have the a priori bound
$d_{M^{tot}\cdot N}^{loc} < d_N$. 
In other words, if we denote by $N$ ($N^{glob}$) the number of reduced basis vectors required
for a localized (global) RB approximation for given $\varepsilon$ and denoting by
$N^{loc}$ the maximum dimension of the local RB spaces $V^m_N$, $V^\gamma_N$, $V^e_N$, $V^v_N$,
we have
\begin{equation}\label{eq:apriori_size}
    N \leq N^{loc}M^{tot} \leq N^{glob} M^{tot} \leq C_{\Grid} M N^{glob}.
\end{equation}
where the constant $C_{\Grid}$ only depends on the topology of the domain decomposition $\Grid$.
Whether or not estimate \eqref{eq:apriori_size} is sharp largely depends on the dependence of the
solution $u(\param)$ on the parameter $\param$.
When a change in $\param$ equally affects the solution in all subdomains $\Omega_m$, we expect that
optimal local RB spaces will be of similar dimension $N^{loc}$ and that $N^{loc} \approx N^{glob}$.
On the other hand, it may be the case that the influence of $\param$ on $u(\param)$ is weak in many
$\Omega_m$, in which case $N \ll N^{loc}M^{tot}$, or that each of the $p$ components of
$\param \in \mathbb{R}^p$ affects $u(\param)$ on different subdomains, in which case $N^{loc}M^{tot}
\ll N^{glob}M^{tot}$.
Thus, the actual loss in online efficiency due to localization will strongly depend on the type of
problem to be solved.

More importantly though, note that the localization of $V_N$ results in a change of the structure of
the reduced system matrices $\mathbb{A}^q$.
While these matrices are dense for global RB approximations, localized RB schemes yield
$\mathbb{A}^q$ with a sparse block structure of $M^{tot} \times M^{tot}$ blocks of maximum dimension
$N^{loc} \times N^{loc}$ and a maximum of $C_{cup}$ blocks per row.
$C_{cup}$ depends on the specific localization method and on the topology of $\Grid$.
For instance, for non-conforming methods, $C_{cup} - 1$ is given by the maximum number of interfaces
of a given subdomain $\Omega_m$, whereas for the ArbiLoMod with a quadrilateral mesh $C_{cup} =
25$.

Thus, estimate \eqref{eq:apriori_size} has to be interpreted in relation to the fact that the computational
complexity for solving \eqref{eq:rom_in_matrices} can be vastly reduced in comparison to
\eqref{eq:rom_effort} by exploiting the structure of the $\mathbb{A}^a$.
In particular, the costs for assembling \eqref{eq:rom_in_matrices} can be reduced to
$\mathcal{O}(C_{\Grid}C_{cup}(N^{loc})^2M)$.
For the solution of \eqref{eq:rom_in_matrices} direct or block-preconditioned iterative solvers can
be used.
For the latter, the computational effort can be expected to increase sub-quadratically in the number
of subdomains $M$.
In the scRBE method, the volume degrees of freedom associated with the spaces $V^m_N$ are eliminated
from \eqref{eq:rom_in_matrices} using static condensation to improve computational efficiency.


\subsection{Offline costs and parallelization}

While the local RB spaces $V^m_N$, $V^\gamma_N$, $V^e_N$, $V^v_N$ can be initialized by decomposing
global solution snapshots $u_h(\param)$ w.r.t.\ \eqref{eq:localizing_space_decomposition} (see
\cite{AlbrechtHaasdonkEtAl2012}), the core element of localized RB methods is the construction of
local RB spaces from local problems associated with the subdomains $\Omega_m$ as described in
Section~\ref{sec:prep}.
This has various computational benefits:

First, we can expect a reduction of computational complexity as for most
linear solvers we expect a super-linear increase in computational complexity for an increasing
dimension of $V_h$, whereas the ratio of the dimensions of $V_h$ and the local subspaces in
\eqref{eq:localizing_space_decomposition} remains constant of order $1/M$ (for the volume spaces
$V_h^m$ and smaller for the spaces $V_h^\gamma$, $V_h^e$, $V_h^v$) as $h \to 0$.
Thus, solving $\mathcal{O}(MN^{loc})$ training problems of size $\dim V_h / M$ is expected to be
faster than solving $N$ global problems.
At the same time, we expect $N^{loc}$ to decrease for $H \to h$.
Thus smaller subdomains $\Omega_m$
will generally lead to shorter offline times at the expense of less-optimal spaces $V_N$.
In particular, for the non-conforming schemes in Subsection~\ref{subsec:ipdg} it is readily seen
that for $H = h$ we have $V_N = V_h$ and that \eqref{eq:rom_solution} and \eqref{eq:fom_solution}
are equivalent.

Even more important than a potential reduction of complexity is the possibility to choose $H$ small
enough such that each local training problem can be solved communication-free on a single compute node
without the need for a high-performance interconnect.
Also the problem setup and the computation of the reduced system \eqref{eq:rom_in_matrices} can be
performed mostly communication-free:
Instead of instantiating a global fine-scale compute mesh, each compute node can generate a local
mesh from a geometry definition, solve training problems for a given local RB space and all coupling
spaces to obtain corresponding block-entries in $\mathbb{A}^q$, $\mathbb{F}^q$.
Only the local geometry and the resulting reduced-order quantities are communicated (see
\cite[Section 8]{Betal17}).
This makes localized RB methods attractive for cloud-based environments, where large computational
resources can be dynamically made available, but communication speed is limited.

Depending on the problem structure, the use of online enrichment (Section~\ref{sec:adaptivity}) can
yield smaller, problem adapted reduced spaces $V_N$.
Similar to the training of $V_N$, online enrichment is based on the solution of small independent
local problems, that can easily be parallelized.
As typically only some fraction of the subspaces of $V_N$ undergo enrichment, less computational
resources need to be allocated during an online-enrichment phase.
It has to be noted, however, that online enrichment leads to a propagation of snapshot data through
the computational domain as the value of the current solution $u_N(\mu)$ at the boundary of the
enrichment-problem domain enters the problem definition.
Thus, to perform online-enrichment, (boundary values of) reduced basis vectors have to be
communicated between compute nodes and the entire reduced basis has to be kept available.


\section{Applications and numerical experiments}\label{sec:examples}

\subsection{Multiscale problems}\label{sec:multiscale_experiments}

We demonstrate the IP localized RB methods from Section \ref{sec:ip-nonconforming} in the context of parametric multiscale problems, such as Example \ref{ex:multiscale_problem}, with a focus on online adaptivity as in Section \ref{sec:adaptivity} (using the a posteriori error estimate from Section \ref{sec:flux-reconstruction}), rather than offline training.
These experiments were first published in the context of the online adaptive LRBMS in \cite{OS2015}.
We consider a multiplicative splitting of the parameter dependency and the multiscale nature of the data functions, in the sense that $\kappa(\param) := \lambda(\param) \kappa_\varepsilon$, with a parametric total mobility $\lambda: \Params \to L^\infty(\Omega)$ and a highly heterogeneous permeability field $\kappa_\varepsilon \in L^\infty(\Omega)^{d \times d}$.
To be more precise, we
consider \eqref{eq:parametric-multiscale} on $\Omega = [0, 5] \times [0, 1]$ with $f(x,y) = 2\cdot 10^3 $ if $(x,y) \in [0.95, 1.10] \times [0.30, 0.45]$, $f(x,y) = -1\cdot 10^3$ if $(x,y) \in [3.00, 3.15] \times [0.75, 0.90]$ or $(x,y) \in [4.25, 4.40] \times [0.25, 0.40]$ and $0$ everywhere else, $\lambda(x,y; \param) = 1 + (1 - \param) \lambda_c(x,y)$, homogeneous Dirichlet boundary values and a parameter space $\Params = [0.1, 1]$.
On each $t \in \grid$, $\kappa_\varepsilon|_t$ is the corresponding $0$th entry of the permeability tensor used in the first model of the 10th SPE Comparative Solution Project (which is given by $100 \times 20$ constant tensors, see \cite{spe10}) and $\lambda_c$ models a channel, as depicted in Fig. \ref{figure::experiments::estimator::spe10::datafunctions}, top left.

\begin{figure}[t]
  \centering%
  \footnotesize%
  \begin{tikzpicture}
    \node[anchor=north east] (left_column)
      at (-0.1,1.6)                              
      {\includegraphics[width=0.475\textwidth,trim=50 68 58 68, clip]{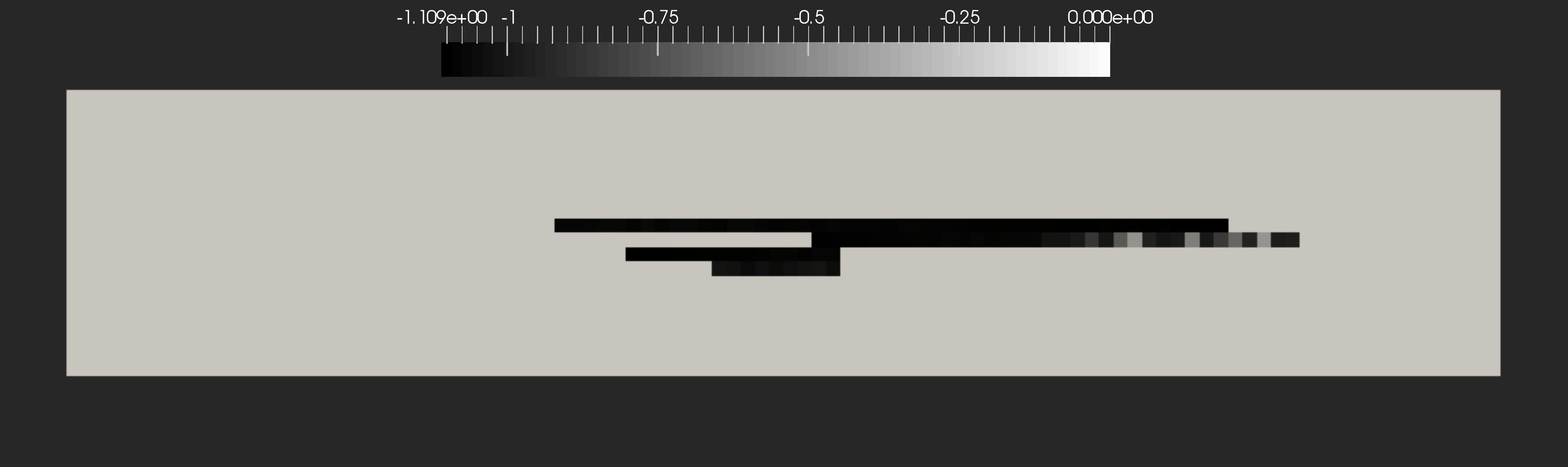}};
    \node[anchor=north east]
      at (-0.1,0)
      {\includegraphics[width=0.475\textwidth,trim=50 68 58 68, clip]{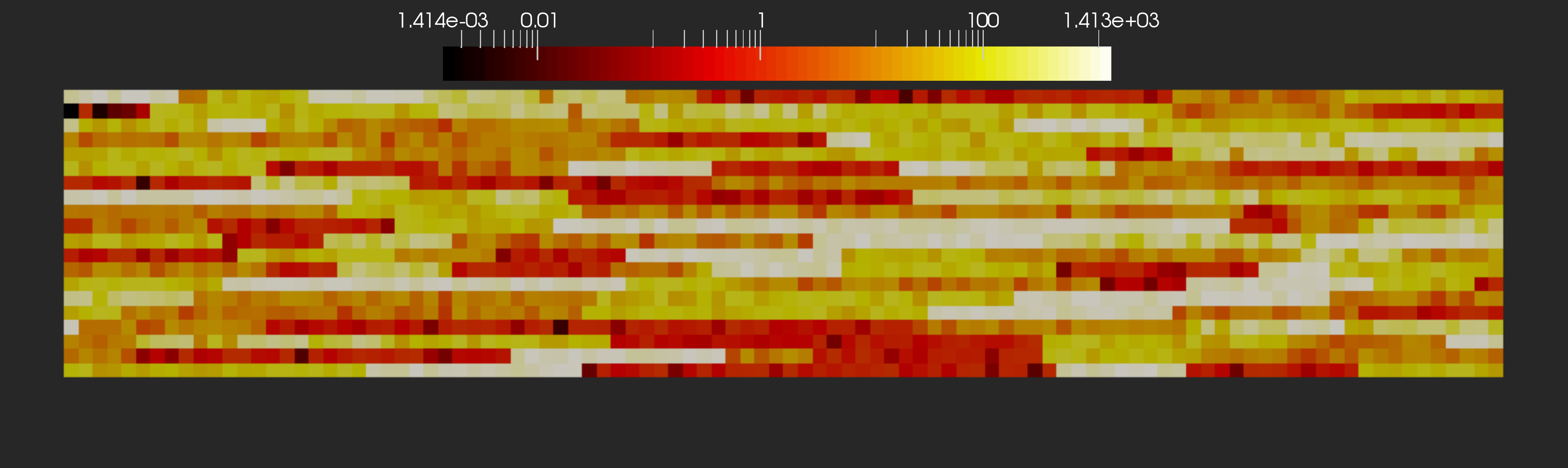}};
    \node[anchor=north east]
      at (-0.1,-1.5)
      {\includegraphics[width=0.475\textwidth,trim=50 69 58 68, clip]{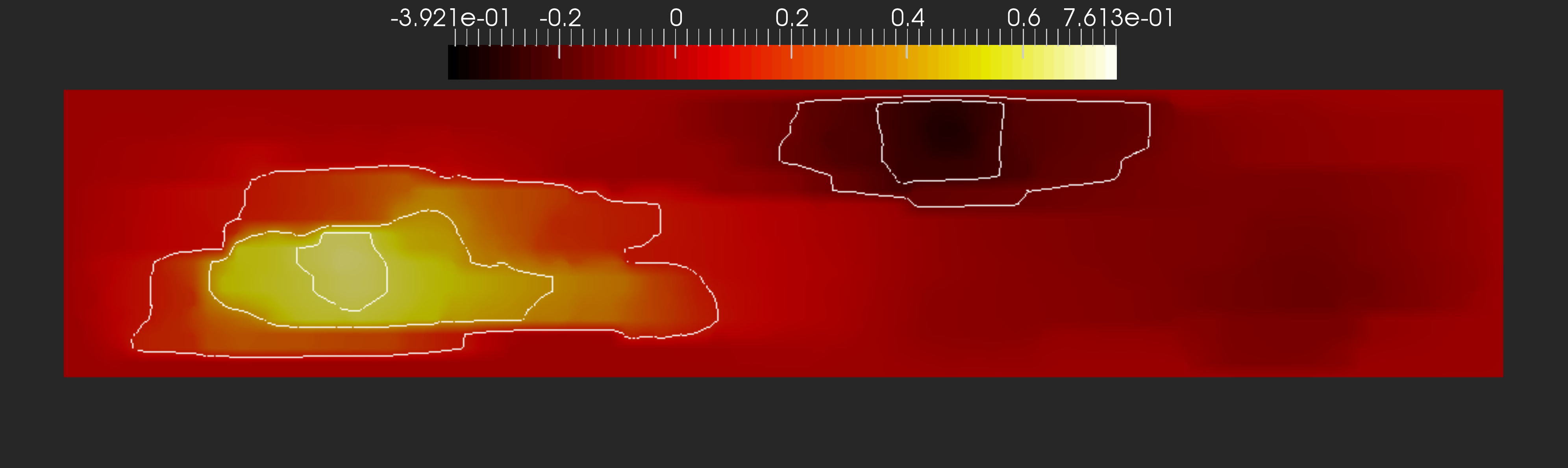}};
    \node[anchor=north east]
      at (-0.1,-3.0)
      {\includegraphics[width=0.475\textwidth,trim=50 69 58 68, clip]{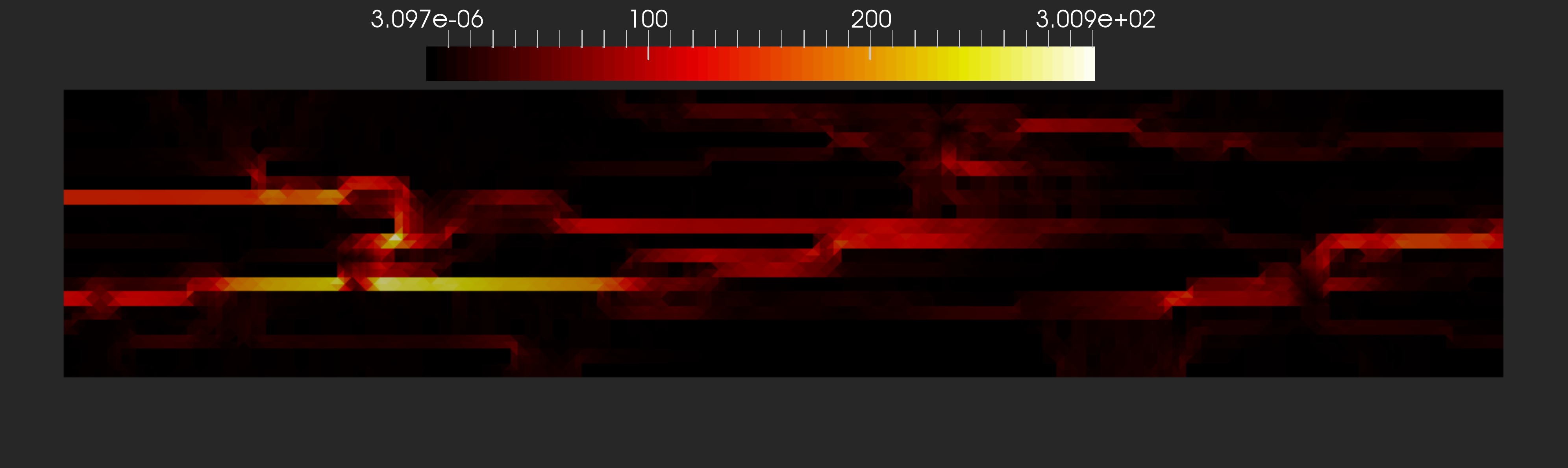}};
    \node[anchor=north west] (right_column)
      at (0.1,1.6)
      {\includegraphics[width=0.475\textwidth,trim=50 69 58 68, clip]{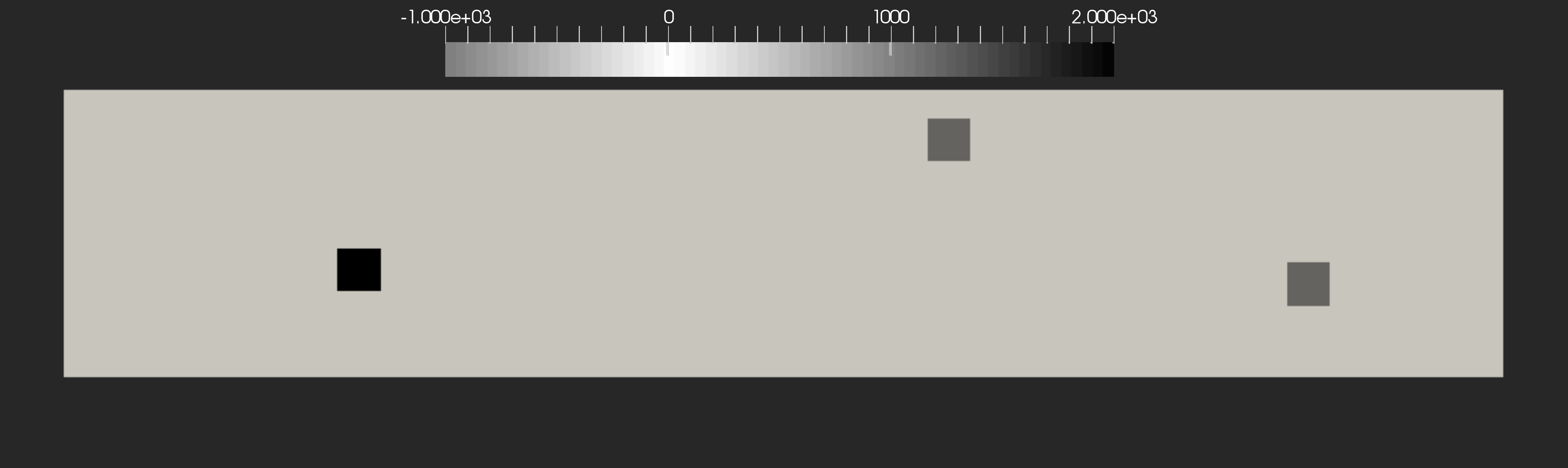}};
    \node[anchor=north west]
      at (0.1,0)
      {\includegraphics[width=0.475\textwidth,trim=50 68 58 68, clip]{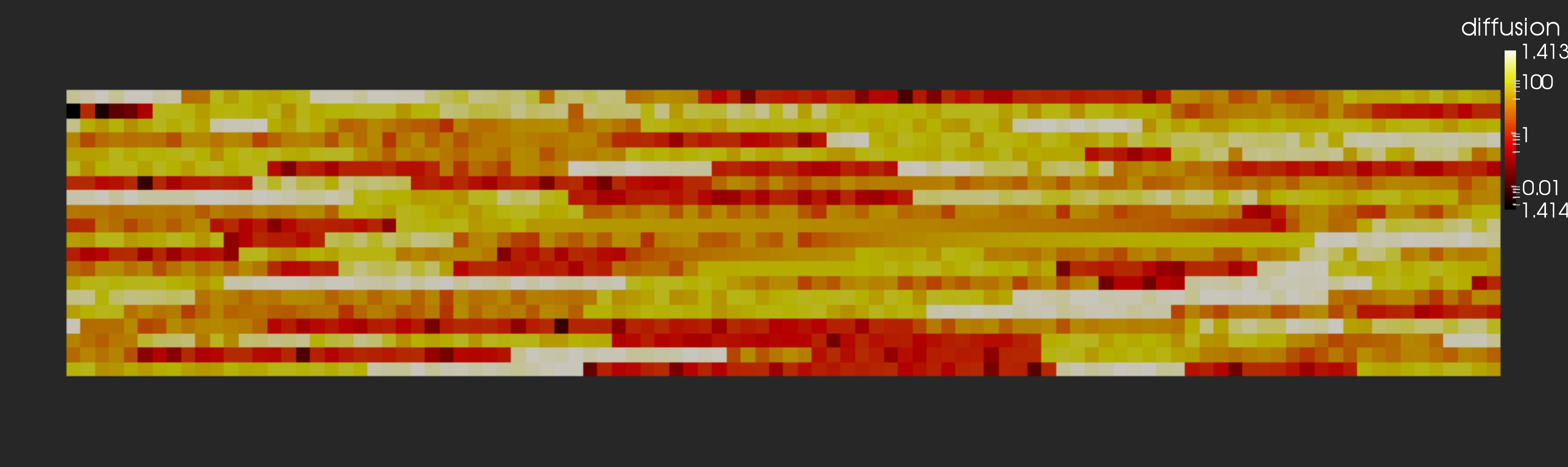}};
    \node[anchor=north west]
      at (0.1,-1.5)
      {\includegraphics[width=0.475\textwidth,trim=50 68 58 68, clip]{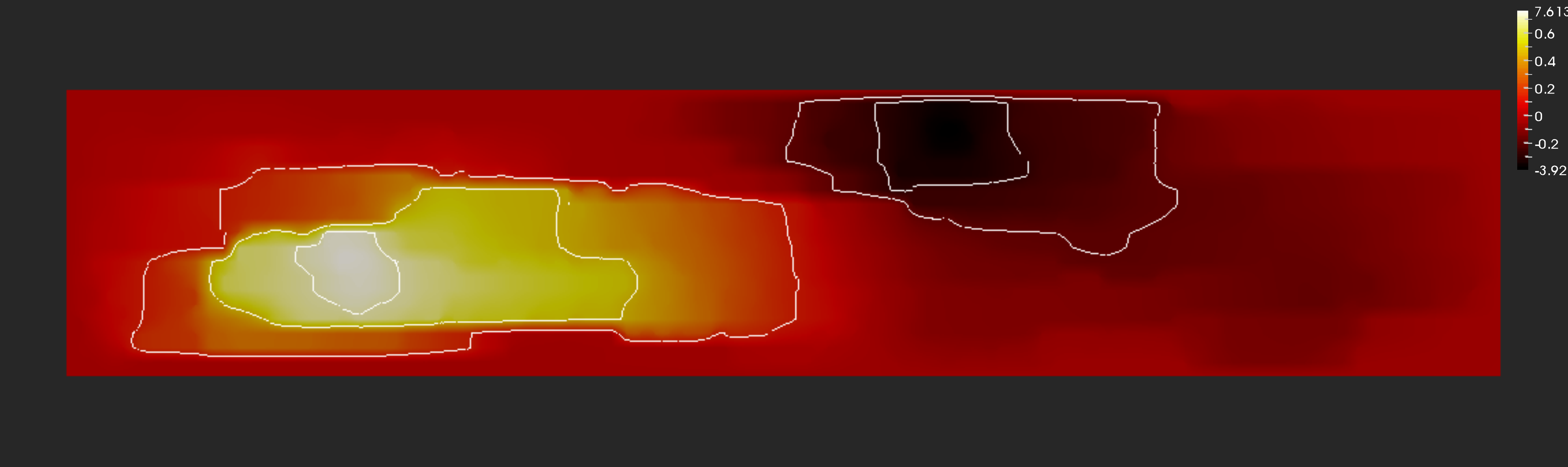}};
    \node[anchor=north west]
      at (0.1,-3.0)
      {\includegraphics[width=0.475\textwidth,trim=50 68 58 68, clip]{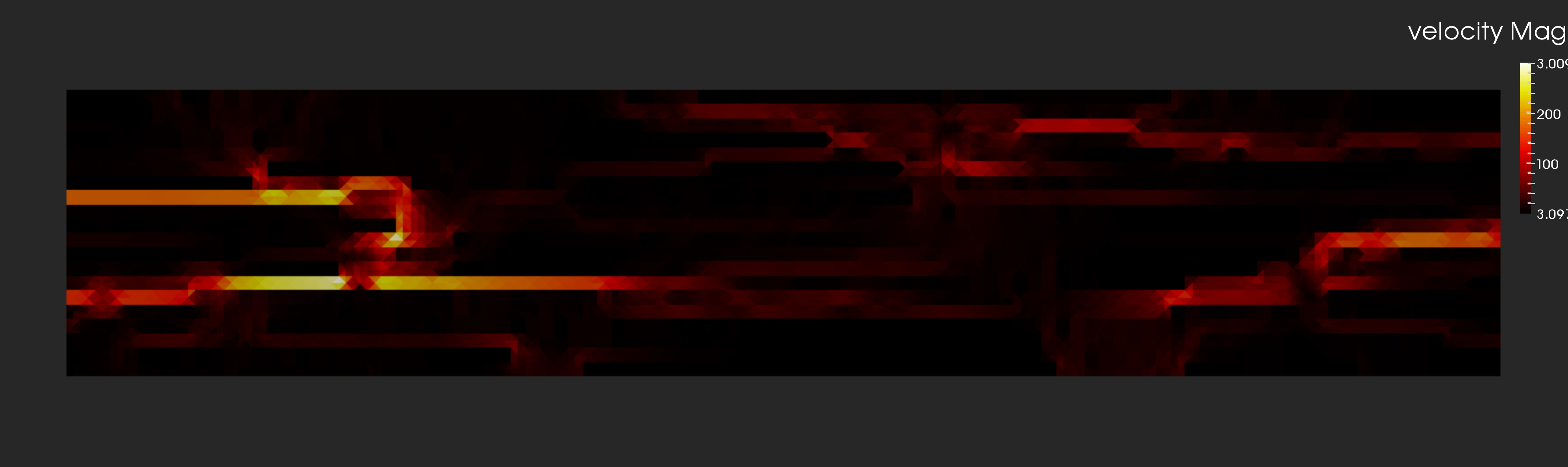}};
    \node[anchor=north]
      at (0,0)
      {\includegraphics[height=40.8mm,width=2mm,trim=1102 188 41 40, clip]{diffusion_mu_min_log}};
    \draw (left_column.south west) -- (right_column.south east);
    \node[anchor=north,scale=1.0] [above=of left_column,above=-50pt] {$\param = 1$};
    \node[anchor=north,scale=1.0] [above=of right_column,above=-50pt] {$\param = 0.1$};
  \end{tikzpicture}
  \caption{%
    Data functions and sample solutions of the experiment in Section \ref{sec:multiscale_experiments}.
    First row: location of the channel function $\lambda_c$ (left) and plot of the force $f$ (right) modeling one source (black: $2\mydot 10^3$) and two sinks (dark gray: $-1\mydot 10^3$, zero elsewhere).
    Second to fourth row: both plots in each row share the same color map (middle) with different ranges per row, for parameters $\param = 1$ (left column) and $\param = 0.1$ (right column).
    From top to bottom: logarithmic plot of $\lambda(\param)\kappa_\varepsilon$ (dark: $1.41\mydot 10^{-3}$, light: $1.41\mydot 10^3$), plot of the pressure $u_h(\param)$ (IP localized FOM solution of \eqref{eq:parametric-multiscale}, dark: $-3.92\mydot 10^{-1}$, light: $7.61\mydot 10{-1}$, isolines at 10\%, 20\%, 45\%, 75\% and 95\%) and plot of the magnitude of the reconstructed diffusive flux $R_h^0[u_h(\param); \param]$ (defined in \eqref{eq::flux_recostruction::1} and \eqref{eq::flux_recostruction::2}, dark: $3.10\mydot 10^{-6}$, light: $3.01\mydot 10^2$).
    Note the presence of high-conductivity channels in the permeability (second row left, light regions) throughout large parts of the domain.
    The parameter dependency models a removal of one such channel in the middle right of the domain (second row right), well visible in the reconstructed Darcy velocity fields (bottom).
  \label{figure::experiments::estimator::spe10::solutions}\label{figure::experiments::estimator::spe10::datafunctions}}
\end{figure}

The right hand side $f$ models a strong source in the middle left of the domain and two sinks in the top and right middle of the domain, as is visible in the structure of the solutions (see Fig. \ref{figure::experiments::estimator::spe10::datafunctions}, third row).
The role of the parameter $\param$ is to toggle the existence of the channel $\lambda_c$.
Thus $\lambda(1)\kappa_\varepsilon = \kappa_\varepsilon$ while $\param = 0.1$ models the removal of a large conductivity region near the center of the domain (see the second row in Fig. \ref{figure::experiments::estimator::spe10::solutions}).
This missing channel has a visible impact on the structure of the pressure distribution as well as the reconstructed velocities, as we observe in the last two rows of Fig. \ref{figure::experiments::estimator::spe10::solutions}.
With a contrast of $10^6$ in the diffusion tensor and an $\varepsilon$ of about $|\Omega|/2,\hnS000$ this setup is a challenging heterogeneous multi-scale problem.

We used several software packages for this numerical experiment and refer to \cite{OS2015} for a full list and instructions on how to reproduce these results.
We would like to mention that all grid-related structures (such as data functions, operators, functionals, products, norms) were implemented in a DUNE-based C++ discretization (which is by now contained in the DUNE extension modules\footnote{\url{https://github.com/dune-community/dune-xt-common/}} and the generic discretization toolbox \texttt{dune-gdt}\footnote{\url{https://github.com/dune-community/dune-gdt/}}), while we used \texttt{pyMOR} \cite{pymor} for everything related to model reduction (such as Gram-Schmidt, Greedy).
We consider a domain decomposition of $|\Grid| = 25 \times 5$ squares, each refined such that the full global grid would consist of $|\grid| = 1,\hnS014,\hnS000$ elements.
For the IP localized FOM, following Section \ref{subsec:ipdg}, we choose on each subdomain $\Omega_m \in \Grid$ the DG space (1st order), product and bilinear form from example \ref{ex:ipdg}.
For error estimation, we employed the flux reconstruction ansatz from Section \ref{sec:flux-reconstruction} using a zero order diffusive flux reconstruction (compare Theorem \ref{thm:flux_reconstruction}).

The sole purpose of these experiments is to demonstrate the capabilities of localized RB methods regarding online enrichment.
We thus initialize the local reduced spaces $V_N^m$ on each subdomain a priori by orthonormalized Lagrangian shape functions of order up to one, thus obtaining a reduced space with poor approximation properties (comparable to a standard DG space w.r.t.~the domain decomposition).
Since we employ the a posteriori error estimate $\eta$ on the full approximation error (including the discretization as well as the model reduction error) from Theorem \ref{thm:flux_reconstruction}, and since we omit grid-refinement in these experiments, the estimated discretization error over all parameters of $1.66$ is a lower bound for the overall approximation error, and we thus choose a tolerance of $\Delta_\text{online} = 2$ for the online enrichment in Algorithm \ref{alg:online_adaptive_LRBMS}.

We compare two different strategies, corresponding to the two plots in Fig.~\ref{figure::experiments::adaptive::spe10::error_evolution}.
In both cases, we simulate an outer-loop application in the online part by randomly choosing ten parameters $\Params_\textnormal{online} \subset \Params$ which are subsequently processed.
For each parameter, the local reduced spaces are enriched according to Algorithm \ref{alg:online_adaptive_LRBMS} and the respective marking strategy, until the estimated error is below the specified tolerance.
Note that the evaluation of the localizable a posteriori error estimate can be fully offline/online decomposed and that after each enrichment only information from a subdomain and its neighbors are required to locally update the offline/online decomposed data.

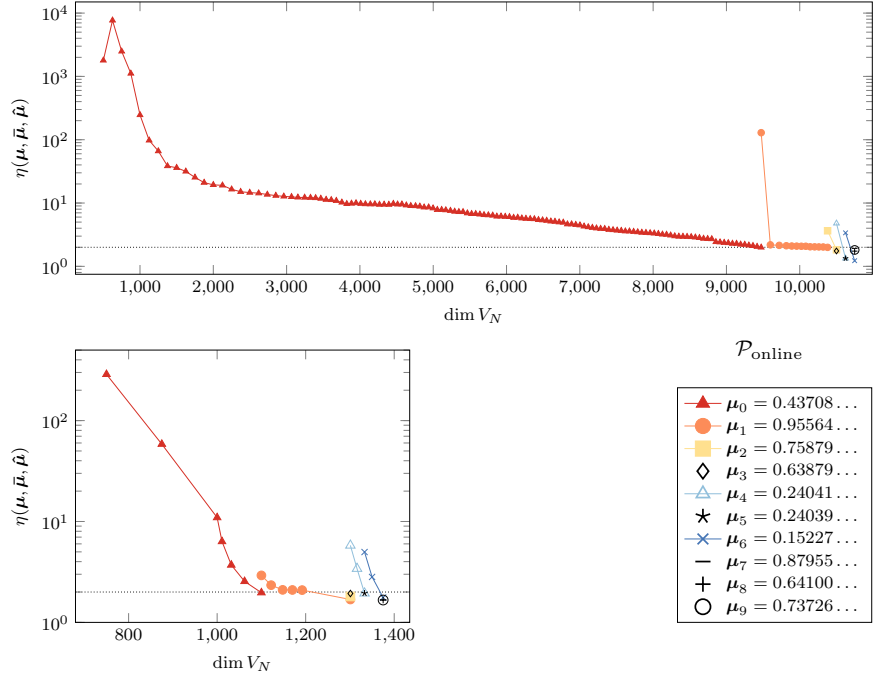
\begin{figure}
  \centering%
  \footnotesize%
  \centering%
  \footnotesize%
  \begin{tikzpicture}
    [scale=0.82,
     spy using outlines={magnification=2, connect spies}]
    \begin{semilogyaxis}[name=top,
                         xlabel={$\dim V_N$},
                         ylabel={$\eta(\param, \paramFixed, \paramHat)$},
                         scaled x ticks = false,
                         enlargelimits=false,
                         width=145mm,
                         height=60mm]
      \addplot [forget plot, black, mark=none] coordinates {
        (120, 15000)
      };
      \addplot [forget plot, black, mark=none] coordinates {
        (10990, 0.75)
      };
      \addplot[forget plot, densely dotted, mark=none, black] coordinates {
        (120, 2)
        (10990, 2)
      };
      \addplot [sixclassRdYlBu1, mark=triangle*, mark options={scale=0.75}] coordinates {
        (500, 1791.33754537)
        (625, 7680.88289396)
        (750, 2486.68081057)
        (875, 1113.40280704)
        (1000, 246.142330746)
        (1125, 97.7376430743)
        (1250, 66.1598863803)
        (1375, 38.4457735398)
        (1500, 36.0896435072)
        (1625, 31.6080475542)
        (1750, 25.4053893157)
        (1875, 21.0221225586)
        (2000, 19.3686800691)
        (2125, 19.0299602266)
        (2250, 16.4754609606)
        (2373, 15.0751096061)
        (2497, 14.6346835105)
        (2619, 14.2930538657)
        (2737, 13.6360137414)
        (2853, 13.1001761902)
        (2962, 12.7360724405)
        (3067, 12.5143946031)
        (3154, 12.2925061639)
        (3246, 12.2319733022)
        (3326, 12.1503742837)
        (3396, 12.0306830511)
        (3467, 11.7639209704)
        (3539, 11.3167997134)
        (3609, 11.2015457826)
        (3678, 10.8684284934)
        (3754, 10.2310399198)
        (3821, 9.69867145456)
        (3885, 9.84190285591)
        (3947, 9.99019406714)
        (4011, 9.77713183329)
        (4075, 9.63987857181)
        (4142, 9.59603509278)
        (4204, 9.6172497561)
        (4267, 9.46266980342)
        (4330, 9.48404190123)
        (4395, 9.44094640865)
        (4456, 9.7501336933)
        (4514, 9.53445596185)
        (4573, 9.5524919607)
        (4633, 9.22341758146)
        (4694, 9.02086928211)
        (4759, 8.98452211813)
        (4818, 8.83891997914)
        (4876, 8.56017418985)
        (4938, 8.58395879656)
        (4998, 8.28091753242)
        (5056, 7.83378955133)
        (5115, 7.82825205162)
        (5172, 7.73319573703)
        (5235, 7.56711245449)
        (5293, 7.3484328586)
        (5350, 7.24490794977)
        (5407, 7.23901681332)
        (5466, 6.91139409544)
        (5524, 6.74616469785)
        (5580, 6.74187306353)
        (5637, 6.62655446227)
        (5693, 6.49064695744)
        (5752, 6.44896733844)
        (5812, 6.37011335393)
        (5868, 6.14824615449)
        (5925, 6.16732247086)
        (5982, 6.09371497642)
        (6040, 6.04357859008)
        (6097, 5.8753347338)
        (6153, 5.89987213515)
        (6209, 5.74662580294)
        (6265, 5.67444494845)
        (6323, 5.61830065415)
        (6380, 5.59523968756)
        (6436, 5.42035242685)
        (6492, 5.36218039354)
        (6548, 5.28052988122)
        (6605, 5.13870159087)
        (6661, 5.08995694036)
        (6718, 4.97873606545)
        (6774, 4.93092564148)
        (6831, 4.6318192624)
        (6888, 4.63697995154)
        (6945, 4.54133952715)
        (7002, 4.46343748876)
        (7058, 4.28066509982)
        (7115, 4.17732748453)
        (7171, 4.08180440305)
        (7227, 3.99674347526)
        (7284, 3.93623176737)
        (7340, 3.91404078941)
        (7396, 3.7973865107)
        (7453, 3.75374591268)
        (7509, 3.69548712328)
        (7565, 3.6545222513)
        (7621, 3.58590971586)
        (7677, 3.58066226866)
        (7734, 3.52611290621)
        (7790, 3.47936039259)
        (7846, 3.41782751792)
        (7902, 3.39040551531)
        (7959, 3.36184178768)
        (8015, 3.32151424963)
        (8071, 3.26454323671)
        (8127, 3.19703158207)
        (8183, 3.17618529978)
        (8239, 3.09046901813)
        (8295, 3.00433162643)
        (8351, 2.96330669963)
        (8407, 2.95178227761)
        (8463, 2.93578057499)
        (8519, 2.91572844739)
        (8576, 2.87190669306)
        (8632, 2.81760496568)
        (8688, 2.74591956037)
        (8744, 2.71872484325)
        (8800, 2.70409469122)
        (8856, 2.43076217239)
        (8913, 2.40408897679)
        (8969, 2.36232664098)
        (9025, 2.33983772045)
        (9082, 2.28882802309)
        (9138, 2.28302701847)
        (9194, 2.21781293469)
        (9250, 2.20618332884)
        (9306, 2.14616833056)
        (9362, 2.10805542878)
        (9418, 2.04617357353)
        (9474, 1.99157821553)
      };
      \addplot [sixclassRdYlBu2, mark=*, mark options={scale=0.75}] coordinates {
        (9474, 129.119684223)
        (9599, 2.18601110121)
        (9720, 2.14538581493)
        (9817, 2.11443927376)
        (9891, 2.09320280548)
        (9958, 2.08534682757)
        (10022, 2.07658852726)
        (10082, 2.07591484604)
        (10142, 2.03557934695)
        (10203, 2.03139596606)
        (10261, 2.02575991265)
        (10320, 2.02156770925)
        (10380, 1.99476356954)
      };
      \addplot [sixclassRdYlBu3, mark=square*, mark options={scale=0.75}] coordinates {
        (10380, 3.6480883673)
        (10501, 1.82774849066)
      };
      \addplot [mark=diamond, mark options={scale=0.66}] coordinates {
        (10501, 1.7504215426)
      };
      \addplot [sixclassRdYlBu5, mark=triangle, mark options={scale=0.75}] coordinates {
        (10501, 4.77498709236)
        (10625, 1.33570172312)
      };
      \addplot [mark=star, mark options={scale=0.66}] coordinates {
        (10625, 1.33567497384)
      };
      \addplot [sixclassRdYlBu6, mark=x, mark options={scale=0.75}] coordinates {
        (10625, 3.38234887057)
        (10749, 1.22776786594)
      };
      \addplot [mark=-, mark options={scale=0.66}] coordinates {
        (10749, 1.91879753945)
      };
      \addplot [mark=+, mark options={scale=0.66}] coordinates {
        (10749, 1.73865463879)
      };
      \addplot [mark=o, mark options={scale=1}] coordinates {
        (10749, 1.81404519339)
      };
    \end{semilogyaxis}
    \begin{semilogyaxis}[name=bottom_left,
                         at=(top.south west),
                         anchor=north west,
                         xlabel={$\dim V_N$},
                         ylabel={$\eta(\param, \paramFixed, \paramHat)$},
                         yshift=-35pt,
                         enlargelimits=false,
                         width=70mm,
                         height=60mm,
                         legend style={at={(1.8,-0.0015)},anchor=south west},
                         legend entries={$\param_0 = 0.43708\dots$,
                                         $\param_1 = 0.95564\dots$,
                                         $\param_2 = 0.75879\dots$,
                                         $\param_3 = 0.63879\dots$,
                                         $\param_4 = 0.24041\dots$,
                                         $\param_5 = 0.24039\dots$,
                                         $\param_6 = 0.15227\dots$,
                                         $\param_7 = 0.87955\dots$,
                                         $\param_8 = 0.64100\dots$,
                                         $\param_9 = 0.73726\dots$}]
      \addlegendimage{sixclassRdYlBu1, mark=triangle*, mark options={scale=1.75}};
      \addlegendimage{sixclassRdYlBu2, mark=*, mark options={scale=1.75}};
      \addlegendimage{sixclassRdYlBu3, mark=square*, mark options={scale=1.75}};
      \addlegendimage{mark=diamond, mark options={thick,scale=1.75},only marks};
      \addlegendimage{sixclassRdYlBu5,mark=triangle, mark options={thick,scale=1.75}};
      \addlegendimage{mark=star, mark options={thick,scale=1.75},only marks};
      \addlegendimage{sixclassRdYlBu6,mark=x, mark options={thick,scale=1.75}};
      \addlegendimage{mark=-, mark options={thick,scale=1.75},only marks};
      \addlegendimage{mark=+, mark options={thick,scale=1.75},only marks};
      \addlegendimage{mark=o, mark options={thick,scale=1.75},only marks};
      \addplot [forget plot, black, mark=none] coordinates {
        (680, 500)
      };
      \addplot [forget plot, black, mark=none] coordinates {
        (1435, 1)
      };
      \addplot[forget plot, densely dotted, mark=none, black] coordinates {
        (680, 2)
        (1435, 2)
      };
      \addplot [forget plot, sixclassRdYlBu1,mark=triangle*] coordinates {
        (750, 287.440079916)
        (875, 58.4401549639)
        (1000, 10.9396349868)
        (1011, 6.33602720962)
        (1032, 3.70553345532)
        (1062, 2.54925802882)
        (1100, 1.96955930977)
      };
      \addplot [forget plot, sixclassRdYlBu2,mark=*] coordinates {
        (1100, 2.92291126423)
        (1122, 2.33908821635)
        (1148, 2.10069106593)
        (1170, 2.09592729593)
        (1192, 2.08985608326)
        (1301, 1.68664296569)
      };
      \addplot [forget plot, sixclassRdYlBu3,mark=square*] coordinates {
        (1301, 1.80585990862)
      };
      \addplot [forget plot, mark=diamond, mark options={scale=0.75}] coordinates {
        (1301, 1.93346080208)
      };
      \addplot [forget plot, sixclassRdYlBu5,mark=triangle, mark options={scale=1.25}] coordinates {
        (1301, 5.79430703102)
        (1316, 3.40990349211)
        (1333, 1.94339058214)
      };
      \addplot [forget plot, mark=star, mark options={scale=0.75}] coordinates {
        (1333, 1.94354276301)
      };
      \addplot [forget plot, sixclassRdYlBu6,mark=x] coordinates {
        (1333, 4.98104151044)
        (1350, 2.822429994 )
        (1375, 1.76565436835)
      };
      \addplot [forget plot, mark=-, mark options={scale=0.75}] coordinates {
        (1375, 1.67276869395)
      };
      \addplot [forget plot, mark=+, mark options={scale=0.75}] coordinates {
        (1375, 1.66465635735)
      };
      \addplot [forget plot, mark=o, mark options={scale=1.125}] coordinates {
        (1375, 1.66598114674)
      };
    \end{semilogyaxis}
    \node [anchor=south] at (11.25,-1.5) {$\Params_\text{online}$};
  \end{tikzpicture}
  \caption{%
    Estimated error evolution during the adaptive online phase for the experiment in Section \ref{sec:multiscale_experiments} with $|\Grid| = 125$, $k_H = 1$, $\Delta_\text{online} = 2$ (dotted line), $\paramFixed = \paramHat = 0.1$, for different on-line and offline strategies:
    no global snapshot (greedy search disabled, $N_\text{greedy} = 0$) during the offline phase, uniform marking during the online phase (top) and two global snapshots (greedy search on $\Params_\text{train} = \{0.1, 1\}$, $N_\text{greedy} = 2$) and combined uniform marking while $\eta(\param, \paramFixed, \paramHat) > \theta_\text{uni} \Delta_\text{online}$ with $\theta_\text{uni} = 10$, D\"orfler marking with $\theta_\text{doerf} = 0.85$ and age-based marking with $N_\text{age} = 4$ (bottom left); note the different scales.
    With each strategy the local reduced bases are enriched according to Algorithm \ref{alg:online_adaptive_LRBMS} while subsequently processing the online parameters $\param_0, \dots, \param_9$ (bottom right).
  \label{figure::experiments::adaptive::spe10::error_evolution}}
\end{figure}
\begin{figure}
  \footnotesize%
  \centering%
  \begin{tikzpicture}
    \node[anchor=north east]
      at (-0.1,0)                                
     {\includegraphics[width=0.475\textwidth,trim=50 68 58 68, clip]{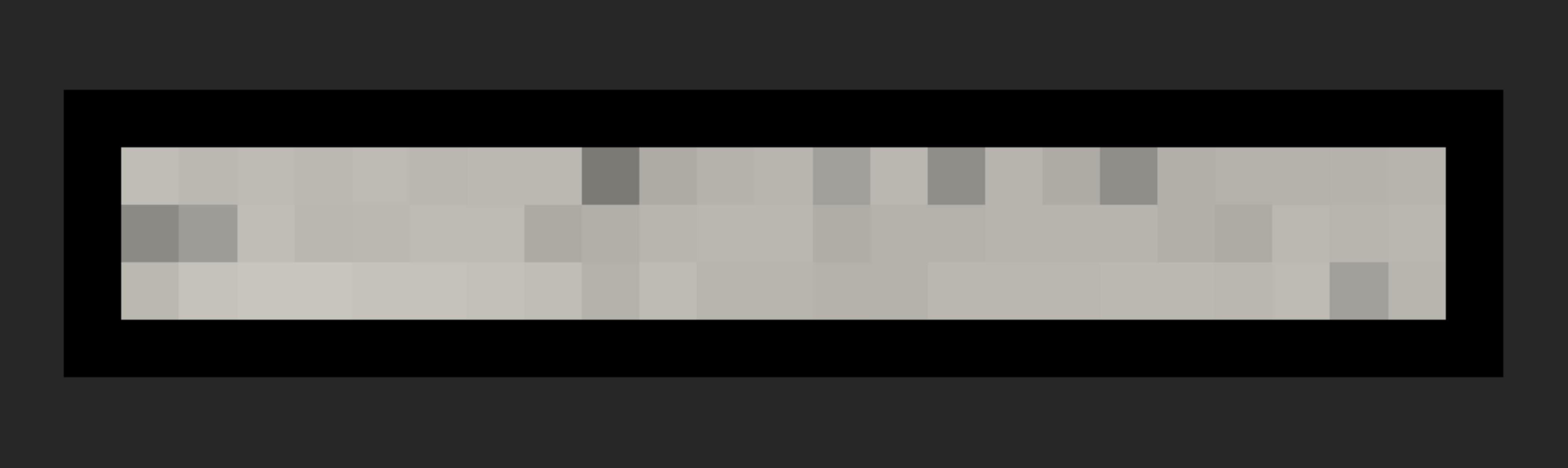}};
    \node[anchor=north west]
      at (0.1,0)
      {\includegraphics[width=0.475\textwidth,trim=50 68 58 68, clip]{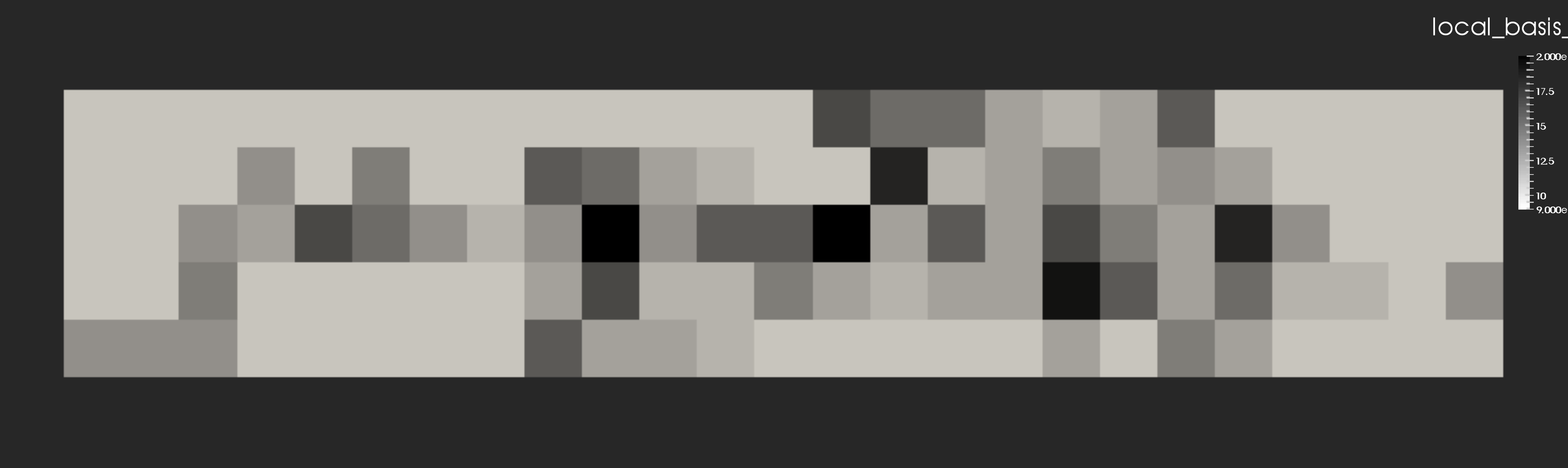}};
    \node[anchor=north]
      at (0,0)
      {\includegraphics[height=11.0mm,width=2mm,trim=1106 210 38 35, clip]{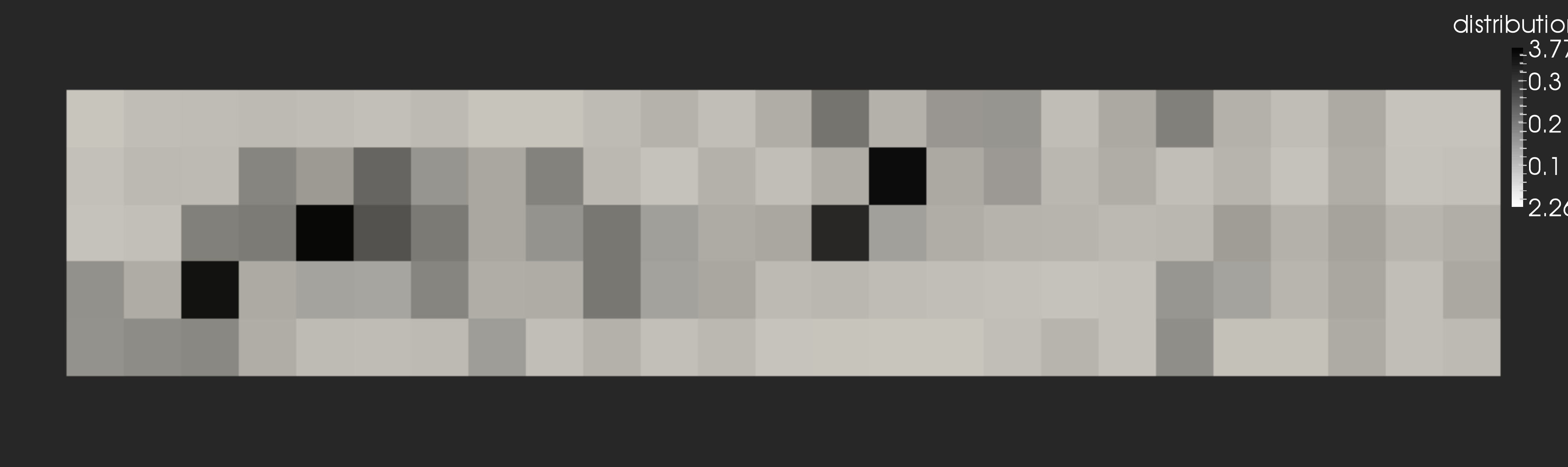}};
  \end{tikzpicture}
  \caption{%
    Spatial distribution of the final sizes of the local reduced bases on each subdomain, after the adaptive online phase for the experiment in Section \ref{sec:multiscale_experiments} with $\Omega = [0, 5]\times[0, 1]$, $|\Grid| = 25 \times 5$ for the two strategies shown in Fig. \ref{figure::experiments::adaptive::spe10::error_evolution}: no global snapshot with uniform enrichment (left, light: 24, dark: 148) and two global snapshots with adaptive enrichment (right, light: 9, dark: 20).
    Note the pronounced structure (right) reflecting the spatial structure of the data functions (compare Fig. \ref{figure::experiments::estimator::spe10::solutions}).}
  \label{figure::experiments::adaptive::spe10::final_local_basis_sizes}
\end{figure}

In the first experiment, we use a uniform marking strategy, which results in an unconditional enrichment on each subdomain (comparable to domain decomposition methods).
As we observe in Fig. \ref{figure::experiments::adaptive::spe10::error_evolution} (top), however, it takes 129 enrichment steps to lower the estimated error below the desired tolerance for the first online parameter $\param_0$.
After this extensive enrichment it takes 12 steps for $\param_1$ and none or one enrichment steps to reach the desired tolerance for the other online parameters.
The resulting coarse reduced space is of size $10,\hnS749$ (with an average of 86 basis functions per subdomain), which is clearly not optimal.
Although each subdomain was marked for enrichment, the sizes of the final local reduced bases differ since the local Gram Schmidt basis extension may reject updates (if the added basis function is locally not linearly independent).
As we observe in Fig. \ref{figure::experiments::adaptive::spe10::final_local_basis_sizes} (left) this is indeed the case with local basis sizes ranging between 24 and 148.
Obviously, a straightforward domain decomposition ansatz without suitable training is not feasible for this setup.
This is not surprising since the data functions exhibit strong multiscale features and non-local high-conductivity channels connecting domain boundaries, see Fig. \ref{figure::experiments::estimator::spe10::solutions}.

To remedy the situation we allow for two global snapshots during the offline phase (for parameters $\param \in \{0.1, 1\}$) and use an adaptive marking strategy which combines uniform marking, D\"orfler marking and age-based marking (see the caption of Fig. \ref{figure::experiments::adaptive::spe10::error_evolution}) in the online phase.
This strategy employs uniform marking until a saturation condition is reached, and afterwards uses a D\"orfler marking combined with a marking based on counting how often a subdomain has not been marked.
With two global solution snapshots incorporated in the basis the situation improves significantly, as we observe in Fig. \ref{figure::experiments::adaptive::spe10::error_evolution} (bottom left).
In total we observe only two enrichment steps with uniform marking (see the first two steps for $\param_0$), which indicates that further offline training would be beneficial.
The number of elements marked range between 11 and 110 (over all online parameters and all but the first two enrichment steps) with a mean of 29 and a median of 22.
Of these marked elements only once have 87 out of 110 elements been marked due to their age (see the last step for $\param_1$).
Overall we could reach a significantly lower overall basis size than in the previous setup ($1,\hnS375$ vs. $10,\hnS749$) and the sizes of the final local bases range between only nine and 20 (compared to 24 to 148 above).
We also observe in Fig. \ref{figure::experiments::adaptive::spe10::final_local_basis_sizes} (right) that the spatial distribution of the basis sizes follows the spatial structure of the data functions (compare Fig. \ref{figure::experiments::estimator::spe10::solutions}), which nicely shows the localization qualities of our error estimator.

\subsection{Fluid dynamics} 
\begin{figure} [th]
\centering
\includegraphics[scale=0.4]{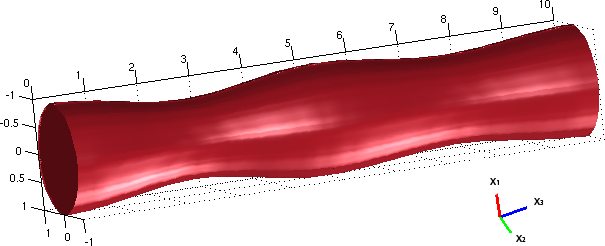}
\caption{\footnotesize{Computational domain ($\mu_1=7, \mu_2=10)$. }}
\label{cilindri4}
\end{figure} 

Flow simulations in pipelined channels  have a growing interest in many biological and industrial applications.
The  localized model order reduction approaches presented in this chapter are  suitable for the study of internal flows in hierarchical parametrized geometries. In particular, the non-conforming approach introduced in Section  \ref{subsec:non-conforming} has applications  in the analysis of the blood flow in specific compartments of the circulatory system  that can be represented as a  combination of few deformed vessels from a reference one.
\begin{figure} [th]
\centering
\includegraphics[width=0.32\linewidth,height=2cm]{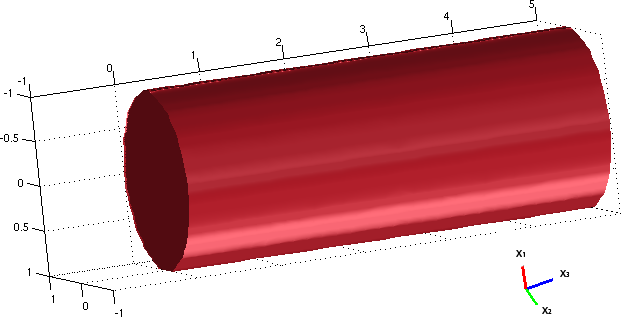}\includegraphics[width=0.32\linewidth,height=2cm]{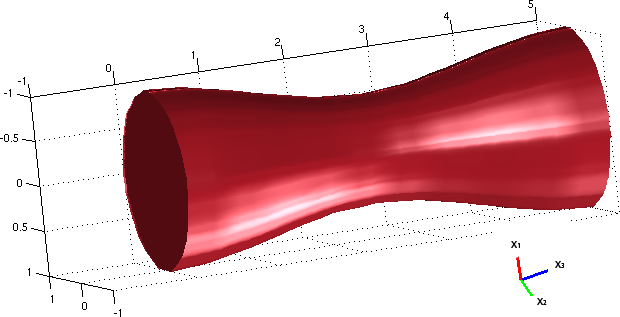}\includegraphics[width=0.32\linewidth,height=2cm]{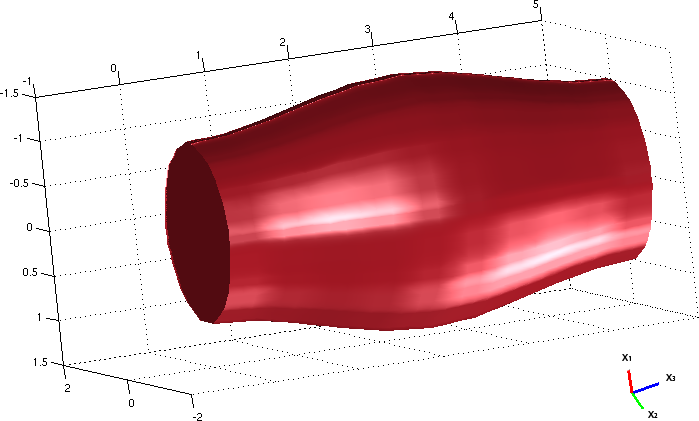}
\caption{\footnotesize{Reference pipe and two deformed pipes ($\mu=-5, \mu=5$): stenosis and aneurysm configuration. }}
\label{cilindri2}
\end{figure}

We want to solve the Stokes equation defined in \eqref{NSequation}, with $\delta=0$, in a computational domain  $\Omega$   composed by two stenosed blocks $\Omega_{\mu_1}$ and $\Omega_{\mu_2}$ (Fig. \ref{cilindri4}), by imposing non-homogeneous BCs $\sigma_n^{in}=[0,5]^{T}$ in the inlet surface ($x_1=10$), non-homogeneous BCs $\sigma_n^{in}=[0,-1]^{T}$ in the outlet surface  ($x_1=0$) and homogeneous Dirichlet BC on the remaining boundaries of the domain.
Here, the Taylor-Hood Finite Element Method has been used to compute the basis functions, $\mathbb{P}_2$ elements for velocity and supremizer,  $\mathbb{P}_1$ for pressure, respectively  and consequently $\mathbb{P}_1(\Gamma_{m,m'})$ for the Lagrange multipliers space.
\begin{figure} [t]
\centering
\includegraphics[width=7 cm,height=4cm]{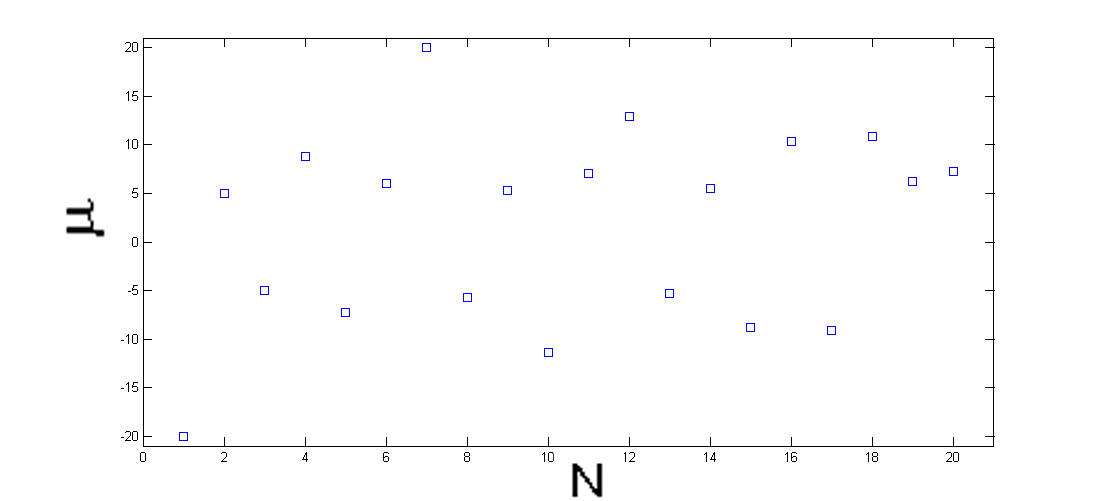}
\caption{\footnotesize{Distribution of the selected parameter values by the greedy algorithm used to generate the basis functions in a single block. }}
\label{sampling3d}
\end{figure}
\noindent Fig. $\ref{sampling3d}$ shows  the distribution of the parameter values selected by the greedy algorithm, by applying the offline stage of the reduced basis method to the single stenosis block. By taking into account that the range $[-5, 5]$ is not admitted, we can see that the higher concentration of values is in the intervals $[-10, -5]$ and $[5, 10]$ in correspondence to larger deformation of the pipe. 

\begin{figure} [b]
\centering
\includegraphics[scale=0.27]{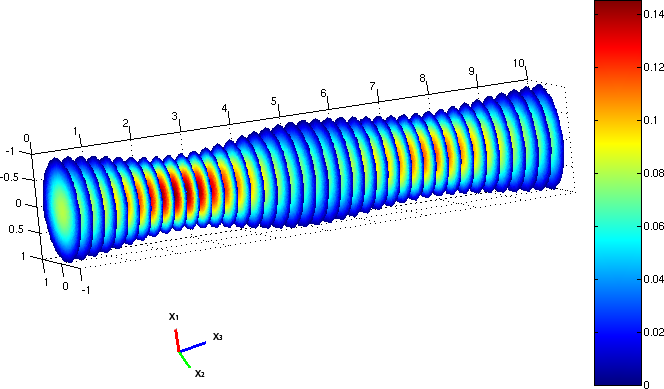}\includegraphics[scale=0.27]{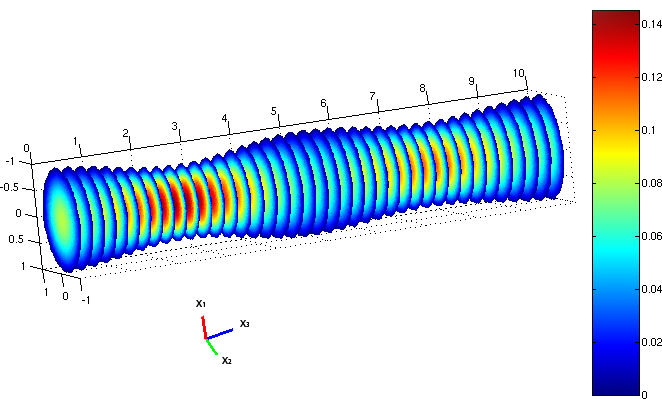}
\caption{\footnotesize{ Representative solutions of velocity    using RBHM  (with $N_1=N_2=19$) (left) and using FEM as a global solution (right), $\mu_1=7, \mu_2=10$.}}
\label{solRB3d}
\end{figure}
\begin{figure} [t]
\centering
\includegraphics[scale=0.27]{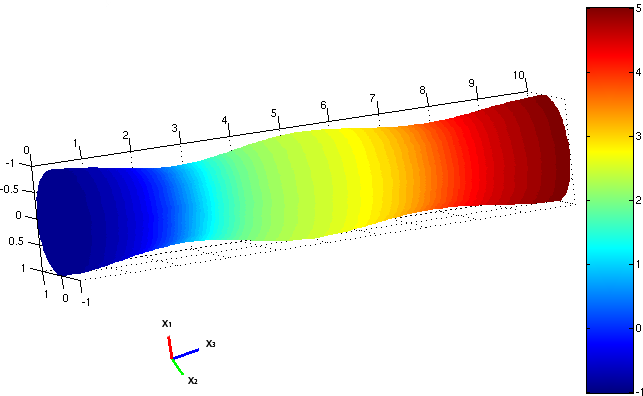}\includegraphics[scale=0.27]{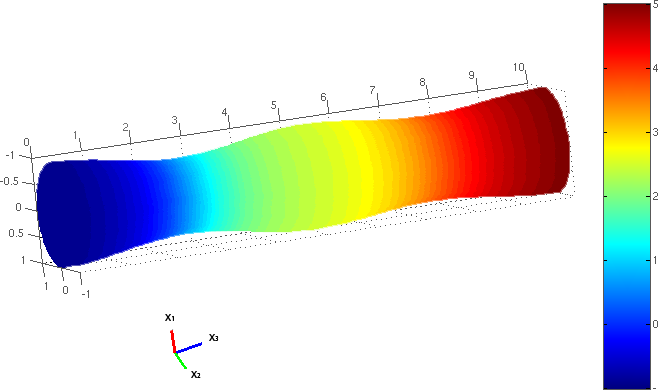}
\caption{\footnotesize{Representative solutions pressure   using RBHM  (with $N_1=N_2=19$) (left) and using FEM as a global solution (right), $\mu_1=7, \mu_2=10$. }}
\label{solFEM3d}
\end{figure}
The geometry of a single stenosis is obtained by the deformation of a reference pipe through a parameter that represents the contraction in the middle of the pipe.
The deformed domain $\Omega_{\mu}$ is mapped from the straight reference pipe $\hat \Omega$ of length $L=5$ and radius $r=1$ through the following coordinate transformation $T_{\mu}:\hat\Omega\to\Omega_{\mu}$ such as $\boldsymbol x=T_{\mu}(\hat {\boldsymbol x})$ and 
$x_1=\hat x_1+\frac{\hat x_1}{\mu}(\cos(\frac{2\pi \hat x_3}{L})-1),
x_2=\hat x_2+\frac{\hat x_2}{\mu}(\cos(\frac{2\pi \hat x_3}{L})-1),
x_3=\hat x_3.$
The range of the parameter $\mu$ is $[-20, -5]\cup[5, 20]$, Fig. \ref{cilindri2} shows the reference pipe and some representative deformations of the geometry.
In order to compute the basis functions, we consider a parametrized Stokes problem  for each subdomain. For the first subdomain, we compute the reduced basis imposing zero Dirichlet condition on the wall, Neumann boundary conditions given by imposing $\sigma_n=\sigma\cdot\mathbf n=\nu \displaystyle \frac{\partial \boldsymbol{u}}{\partial \mathbf{n}}  -p \mathbf{n}$  to be $\sigma_n^{in}=[0,5]^{T}$ on $\Gamma_{in}$ and $\sigma_n^{out}=\boldsymbol{0}$ on the internal interface. For the second subdomain, we compute the reduced basis imposing zero Dirichlet condition on the wall, Neumann boundary conditions imposing $\sigma_n^{in}=\boldsymbol 0$ on the internal interface  and $\sigma_n^{out}=[0,-1]^{T}$ on the outflow interface $\Gamma_{out}$. 

Moreover, we enrich the local RB spaces by a coarse finite element solution of the problem computed in the global domain. This strategy ensures not only the continuity of the velocity, but also the one of the normal stress along the internal interface. For this reason this method is called reduced basis hybrid method.
 Coarse and fine grids have been chosen in order to deal with respectively 155 and 2714 nodes in a single block domain.
\noindent Fig. \ref{solRB3d} shows a representative flow solution in $\Omega$, found with the reduced basis hybrid method, to be compared with the finite element solution.  The same comparison, regarding the pressure solutions, is shown in Fig. \ref{solFEM3d}.


\section{Further perspectives}\label{sec:perspectives}

\subsection{Parabolic problems} 

Most of the techniques presented in this chapter so far can be extended or even directly applied to parabolic problems. For instance, local approximation spaces that are optimal in the sense of Kolmogorov are proposed in \cite{Sch19} and the LRBMS for parabolic problems is presented in \cite{ORS2017,OR17}.
To facilitate an adaptive construction of the local reduced space or online-adaptivity, a suitable, localized a posteriori error estimator is key. Therefore, we present in this subsection an abstract framework for a posteriori error estimation 
for approximations of scalar parabolic evolution equations, based on elliptic reconstruction techniques. For further reading and the application to 
localized model reduction we refer to \cite{GeorgoulisLakkisEtAl2011,ORS2017}. 

\begin{definition}[Parameterized parabolic problem in variational form]
  \label{def:parametrized_problem_parabolic}
Let a Gelfand triple of suitable Hilbert spaces $V \subset H = H'  \subset V'$, an end time $T_\text{end} > 0$, initial data $u_0 \in V$ and right hand side $f \in H$ be given. 
For a parameter $\param \in \Params$ find $u(\cdot; \param) \in L^2(0, T_\text{end}; V)$ with $\partial_t u(\cdot; \param) \in L^2(0, T_\text{end}; V')$, such that $u(0; \param) = u_0$ and
\begin{align}
  \dualpair{\partial_t u(t; \param)}{q}  + a\big(u(t;\param), v; \param\big) &= f(v; \param) &&\text{for all } v \in V.
  \label{eq:parabolic_weak_solution}
\end{align}
\end{definition}

Depending on the error we want to quantify, the space $V$ in\ \eqref{eq:parabolic_weak_solution} can be either a analytical function space as in\ \eqref{eq:weak_solution} or an already discretized function space $V_h$.
We drop the parameter dependency in this section to simplify the notation.

\begin{definition}[Approximations of the parabolic problem]
Let $\widetilde{V} \subseteq H$ be a finite dimensional approximation space for $V$,
not necessarily contained in $V$. Potential candidates for $\widetilde{V}$ are conforming or non-conforming 
localized model reduction spaces $V_N$ as discussed above, but also finite element or finite volume spaces fit into this setting. 
Denote by $\SP{\cdot}{\cdot}$, $\|\cdot\|$ the $H$-inner product and the norm induced by it.

Let $f \in H$, and let $a_h: (V + \widetilde{V}) \times (V + \widetilde{V}) \to \R$ be a discrete bilinear
form which coincides with $a$ on $V \times V$ and is thus continuous and coercive on $V$.
Let further $\enorm{\cdot}$ be a norm over $V + \widetilde{V}$, which coincides with the square root of the symmetric part of $a_h$ over $V$.

Our goal is to bound the error $e(t) := u(t) - \tilde{u}(t)$ between the analytical (or discrete) solution $u \in L^2(0, T_\text{end}; V)$,
$\partial_t u \in L^2(0, T_\text{end}; V^\prime)$ of \eqref{eq:parabolic_weak_solution},
where the duality pairing $\dualpair{\partial_t u(t)}{v}$ is induced by the $H$-scalar product
via the Gelfand triple 
and the $\widetilde{V}$-Galerkin approximation $\tilde{u} \in L^2(0, T_\text{end}, \widetilde{V})$, $\partial_t \tilde{u} \in L^2(0, T_\text{end}, \widetilde{V})$, solution of
\begin{equation}\label{eq:reduced_problem}
	\SP{\partial_t \tilde{u}(t)}{ \tilde{v}} + a_h(\tilde{u}(t), \tilde{v}) = \SP{f}{ \tilde{v}} \qquad\text{for all } \tilde{v} \in \widetilde{V}.
\end{equation}
\end{definition}

\begin{definition}[Elliptic reconstruction]
  \label{def:elliptic_reconstruction}
	Denote by $\widetilde{\Pi}$ the $H$-orthogonal projection onto $\widetilde{V}$.
	For $\tilde{v} \in \widetilde{V}$, define the elliptic reconstruction $\Ell(\tilde{v}) \in V$ of $\tilde{v}$ to be the unique solution of the
	variational problem
	\begin{equation}\label{eq:elliptic_reconstruction}
		a_h(\Ell(\tilde{v}), v) = \SP{A_h(\tilde{v}) - \widetilde{\Pi}(f) + f}{v} \qquad\text{for all } v \in V,
	\end{equation}
	where $A_h(\tilde{v}) \in \widetilde{V}$ is the $H$-inner product Riesz representative  of the functional $a_h(\tilde{v}, \cdot)$, i.e.,
	$\SP{A_h(\tilde{v})}{\tilde{v}^\prime} = a_h(\tilde{v}, \tilde{v}^\prime)$ for all $\tilde{v}^\prime \in \widetilde{V}$.
	Note that $\Ell(\tilde{v})$ is well-defined, due to the coercivity of $a_h$ on $V$.
\end{definition}

From the definition it is clear that $\tilde{v}$ is the $\widetilde{V}$-Galerkin approximation of the elliptic reconstruction $\Ell(\tilde{v})$. 

Let us assume that for each $t$ we have a decomposition 
$\tilde{u}(t) =: \tilde{u}^c(t) + \tilde{u}^d(t)$ (not necessarily unique) where $\tilde{u}^c(t) \in V$,
$\tilde{u}^d(t) \in \widetilde{V}$ are the conforming and non-conforming parts of $\tilde{u}(t)$.
	We consider the following error quantities:
	\begin{align*}
		\rho(t)&:=u(t) - \Ell(\tilde{u}(t)), & \varepsilon(t) &:= \Ell(\tilde{u}(t)) - \tilde{u}(t), \\
		e^c(t) &:= u(t) - \tilde{u}^c(t), & \varepsilon^c(t)&:= \Ell(\tilde{u}(t)) - \tilde{u}^c(t).
	\end{align*}

\begin{theorem}[Abstract semi-discrete error estimate]\label{thm:abstract_estimate}
Let $C:=(2\gamma_h^2 + 1)^{1/2}$, where $\gamma_h$ denotes the continuity constant of $a_h$ on $V$ w.r.t.~$\enorm{\cdot}$, then
\begin{align*}\label{eq:abstract_estimate}
  \|e\|_{L^2(0, T_\text{end}; \enorm{\cdot})}\leq &\|e^c(0)\| + \sqrt{3} \|\partial_t \tilde{u}^d\|_{L^2(0, T_\text{end}; \enorm{\cdot}_{V,-1})} \\
  & + (C+1)\cdot\|\varepsilon\|_{L^2(0, T_\text{end}; \enorm{\cdot})} + C\cdot\|\tilde{u}^d\|_{L^2(0, T_\text{end}; \enorm{\cdot})}.
\end{align*}
\end{theorem}

Note that $\varepsilon(t)$ denotes the approximation error of the coercive variational problem \eqref{eq:elliptic_reconstruction}. 
Hence, this error contribution can be controlled by invoking any (localized) a posteriori error estimate for coercive variational problem
as e.g. presented in Section \ref{sec:apost}.

It is straightforward to modify the estimate in Theorem \ref{thm:abstract_estimate} for semi-discrete solutions $\tilde{u}(t)$ to
take the time discretization error into account:

\begin{corollary}\label{thm:fully_discrete_estimate}
  Let $\tilde{u} \in L^2(0, T_\text{end}, \widetilde{V})$, $\partial_t \tilde{u} \in L^2(0, T_\text{end}, \widetilde{V})$ be an arbitrary discrete approximation of
	$u(t)$, not necessarily satisfying \eqref{eq:reduced_problem}.
	Let $\mathcal{R}_T[\tilde{p}](t) \in \widetilde{V}$ denote the $\widetilde{V}$-Riesz representative w.r.t.~the $H$-inner product
	of the time-stepping residual of $\tilde{u}(t)$, i.e.
	\begin{equation*}
		(\mathcal{R}_T[\tilde{u}](t), \tilde{v}) = \SP{\partial_t \tilde{u}(t)}{\tilde{v}} + a_h(\tilde{u}(t), \tilde{v}) - \SP{f}{\tilde{q}} \qquad\forall \tilde{v}
		\in \widetilde{v}.
	\end{equation*}
	Then, with $C:=(3\gamma_h^2 + 2)^{1/2}$, the following error estimate holds:
	\begin{equation*}
		\begin{aligned}
      \|e\|_{L^2(0, T_\text{end}; \enorm{\cdot})}
      & \leq \|e^c(0)\| + 2 \|\partial_t \tilde{u}^d\|_{L^2(0, T_\text{end}; \enorm{\cdot}_{V,-1})} \\
      & \qquad + (C+1)\cdot\|\varepsilon\|_{L^2(0, T_\text{end}; \enorm{\cdot})} + C\cdot\|\tilde{u}^d\|_{L^2(0, T_\text{end}; \enorm{\cdot})} \\
      & \qquad + 2C_{H,V}^b\cdot\|\mathcal{R}_T[\tilde{u}]\|_{L^2(0, T_\text{end}; H)}.
		\end{aligned}
	\end{equation*}
\end{corollary}

\subsection{Non-affine parameter dependence and non-linear problems} 
A key ingredient towards model order reduction for nonlinear problems is the empirical interpolation method (EIM) introduced in \cite{BMNP04} and further developed in \cite{DHO12}, \cite{CS10,MMPY15}. 

In the context of localized model order reduction empirical interpolation has been employed in, e.g., \cite{CEGG14,Presho2015,OR17}.
Based on the concept of empirical operator interpolation from \cite{DHO12} localization strategies can be employed as follows.
To present the main ideas, let us assume the simple situation that 
$$V_h = \bigoplus_{m=1}^M V_h^{m}$$
and that we have a localized decomposition as follows 
$$
a_h\big(u_h(\param), v_h; \param\big)  =  \textstyle\sum_{m=1}^M a_h^m(u_h^m(\param), v_h^m; \param),
$$
with $a_h\big(u_h(\param), \cdot; \param\big) \in (V_h)'$.
The strategy will then rely on an empirical operator interpolation of the local volume operators $a_h^m(u_h^m(\param), \cdot; \param) \in (V_h^{m})'$
 and will thus only involve localized computations in the construction of the interpolation operator.
As an example, the interpolation of the local volume operator will be of the form
$$
{\mathcal I}_L^m \big[a_h^m(u_h^m(\param), \cdot; \param)\big] = \sum_{l=1}^L {\text{\small \calligra S \hspace*{.001em}}_l}^m \big(a_h^m(u_h^m(\param), \cdot; \param)\big)\, q^m_l 
$$
for a local collateral basis $\{q^m_l\}_{l=1}^L \subset {(V_h^m)}'$
and corresponding interpolation functionals $\{{\text{\small \calligra S \hspace*{.001em}}_l}^m\}_{l=1}^L \subset {\Sigma_h^m}''$ from a suitable local dictionary ${\Sigma_h^m}'' \subset (V_h^m)''$, the choice of which is crucial to ensure the accuracy as well as an online-efficient evaluation of the interpolant.
Note that due to the isomorphism between $V_h^m$ and its bi-dual, the local dictionary of interpolation functionals ${\Sigma_h^m}''$ can be identified with a dictionary of functions 
$\Sigma_h^m \subset V_h^m$, such that ${\text{\small \calligra S \hspace*{.001em}}_l}^m\big(  a_h^m(u_h^m(\param), \cdot; \param) \big) = a_h^m(u_h^m(\param), \sigma_l^m; \param)$, where $\sigma_l^m \in \Sigma_h^m$ corresponds to ${\text{\small \calligra S \hspace*{.001em}}_l}^m \in {\Sigma_h^m}''$.
An online-efficient evaluation of the interpolated operator ${\mathcal I}_L^m   \big[a_h^m(u_h^m(\param), \cdot; \param)\big]$ can be ensured by choosing the local dictionary $\Sigma_h^m$ such that the computational complexity of the evaluation $a_h^m(u_h^m(\param),  \sigma^m; \param)$ for $\sigma^m \in \Sigma_h^m$ does not depend on the dimension of $V_h^m$.
The choice of $\Sigma_h^m$ thus depends on the underlying discretization: possible choices in the context of finite element schemes include the finite element basis of $V_h^m$.
Other choices of $\Sigma_h^m$ are conceivable and could improve the interpolation quality, which is subject to further investigation.



\bibliographystyle{alpha}

\end{document}